\documentclass[graybox]{svmult}
\usepackage{amssymb}
\def\isfullversion{true}

\usepackage{helvet}         
\usepackage{courier}        
\usepackage{type1cm}        
\usepackage{makeidx}         
\usepackage{graphicx}        
\usepackage{multicol}        
\usepackage[bottom]{footmisc}
\usepackage{amsmath, latexsym}
\usepackage{xspace}
\usepackage{array}
\usepackage{tabularx}
\usepackage{url}
\usepackage{relsize}
\usepackage{ifthen}
\usepackage{xcolor}

\newcommand{\correction}{} 
\newcommand{\Thm}[1]{Thm.~\ref{#1}}

\newcommand{\fullversion}{\cite{gunnFull2010}\xspace}
\newcommand{\fullversioncolor}{\color{black!50!blue}}
\newcommand{\versions}[2]{\ifthenelse{\equal{\isfullversion}{true}}{{\fullversioncolor#1}}{#2}}

\renewcommand{\vec}[1]{\mathbf{#1}}
\newcommand{\quot}[1]{``#1''}
\newcommand{\R}[1]{\mathbb{R}^{#1}\xspace}
\newcommand\proj[1]{\mathbf{P}(#1)}
\newcommand\RP[1]{\mathbb{R}{P^{#1}}\xspace}
\newcommand{\e}[1]{\vec{e}_{#1}}
\newcommand{\EE}[1]{\vec{E}_{#1}}
\newcommand{\one}{\vec{1}}
\newcommand{\eye}{\vec{I}}
\newcommand{\inert}{\mathbf{A}}
\newcommand{\stripeye}[1]{\ensuremath{\mathbf{S}(#1)}}
\newcolumntype{Y}{X}

\newcommand{\inpro}{\cdot}
\newcommand{\beam}{spear}
\newcommand{\bivo}[1]{{\vec{#1}}}
\newcommand{\pip}{\bivo{\Xi}}

\newcommand{\pvelo}{\bivo{\Gamma}}
\newcommand{\velo}{\bivo{\Omega}}
\newcommand{\momo}{\bivo{\Pi}}
\newcommand{\foro}{\bivo{\Delta}}
\newcommand{\fio}{\bivo{\Phi}}
\newcommand{\thio}{\bivo{\Theta}}
\newcommand{\hard}{${}^*$}
\newcommand{\grade}[2]{\langle #1 \rangle_{#2}}

\newcommand{\pdgrass}[1]{\proj{\bigwedge{\mathbb{R}^{#1*}}}}

\newcommand{\pdgrassgr}[2]{\proj{\bigwedge^{#2}{\mathbb{R}^{#1*}}}}

\newcommand{\pclal}[3]{\proj{\mathbb{R}_{#1,#2,#3}}\xspace}
\newcommand{\pdclal}[3]{\proj{\mathbb{R}^*_{#1,#2,#3}}\xspace}
\newcommand{\pdclplus}[3]{\proj{\mathbb{R}^{*+}_{#1,#2,#3}}\xspace}
\newcommand{\spin}[3]{\mathbf{Spin({#1,#2,#3})}}
\newcommand{\exerc}{\subsubsection{Exercises}}
\newcommand{\linespace}{\ensuremath{\mathbf{L}^{3,3}}\xspace}
\newcommand{\complexspace}{\ensuremath{\mathfrak{B}}\xspace}

\newcommand{\mysubsubsection}[1]{\hspace{-.0in}\textbf{#1.}}
\newcommand{\mymarginpar}[1]{}

\newcommand{\Eq}[1]{(\ref{#1})}
\newcommand{\Fig}[1]{Fig.~\ref{#1}}

\newcommand{\Sec}[1]{Sect.~\ref{#1}}

\newcommand{\Theorem}[1]{Theorem~\ref{#1}}

\begin{document}

\title*{On the Homogeneous Model of Euclidean Geometry}
\titlerunning{On the Homogeneous Model of Euclidean Geometry}
\author{Charles Gunn}
\institute{Charles Gunn \at DFG-Forschungszentrum Matheon,  MA 8-3, Technisches Universit\"{a}t Berlin, Str. des 17. Juni 136, D-10623 Berlin \email{gunn@math.tu-berlin.de}
}
\maketitle
\abstract{
We attach the degenerate signature $(n,0,1)$  to the projectivized dual Grassmann algebra $\pdgrass{(n+1)}$ to obtain the Clifford algebra $\pdclal{n}{0}{1}$ and explore its use as a model for  euclidean geometry.  We avoid problems with the degenerate metric  by constructing an algebra isomorphism $\mathbf{J}$ between the Grassmann algebra and its dual that yields non-metric  \emph{meet} and  \emph{join} operators. 
We focus on the cases of $n=2$ and $n=3$ in detail, 
enumerating the geometric products between $k$- and $l$-blades. We establish that sandwich operators of the form $\vec{X} \rightarrow \vec{g} \vec{X} \widetilde{\vec{g}}$ provide all euclidean isometries, both direct and indirect.  
We locate the spin group, a double cover of the direct euclidean group, inside the even subalgebra of the Clifford algebra, and provide a simple algorithm for calculating the logarithm of such elements. We conclude with an elementary account of euclidean rigid body motion within this framework.}

\section{Introduction}
The work presented here was motivated by the desire to integrate the work of Study on dual quaternions (\cite{study03}) into a Clifford algebra setting.   
The following exposition introduces the modern mathematical structures -- projective space, exterior algebra, and Cayley-Klein metrics -- required to imbed the dual quaternions as the even subalgebra of a particular Clifford algebra, and shows how the result can be applied to euclidean geometry, kinematics, and dynamics.  \versions{This paper is an extended version of a paper with the same name which was published in \cite{gunn2011}.  To help identify extensions, text not contained in the original article appears dark blue.  The original publication is available at \url{www.springerlink.com}.}{Those interested in more details, exercises, and background material are referred to an extended on-line version  (\fullversion).}
\section{The Grassmann Algebra(s) of Projective Space}
\label{sec:gaps}

\textbf{Real projective n-space} $\RP{n}$
is obtained from the $(n+1)$-dimensional euclidean vector space $\R{n+1}$ by introducing an equivalence relation  on vectors $\vec{x}, \vec{y} \in \R{n+1} \setminus \{\mathbf{0}\}$ defined by: $\vec{x}  \sim \vec{y} \iff \vec{x} = \lambda\vec{y}$ for some $\lambda \neq 0$.  That is, points in $\RP{n}$ correspond to lines through the origin in $\R{n+1}$.  

\mysubsubsection{Grassmann algebra} The Grassmann, or exterior, algebra $\bigwedge (\R{n})$, is generated by the outer (or exterior) product $\wedge$ applied to the vectors of $\R{n}$.  The outer product is an alternating, bilinear operation. The product of a $k$- and $m$-vector is a $(k+m)$-vector when the operands are linearly independent subspaces.  An element that can be represented as a wedge product of $k$ 1-vectors is called a simple $k$-vector, or $k$-\emph{blade}. The $k$-blades generate the vector subspace $\bigwedge^k (\R{n})$, whose elements are said to have grade $k$. This subspace has dimension  $n \choose k$, hence the total dimension of the exterior algebra is $2^n$.   $\bigwedge^n (\R{n})$ is one-dimensional, generated by a single element $\eye$ sometimes called the pseudo-scalar.

\mysubsubsection{Simple and non-simple vectors} A $k$-blade represents the subspace of $\R{n}$ spanned by the $k$ vectors which define it.  Hence, the exterior algebra contains within it a representation of the subspace lattice of $\R{n}$.  For $n>3$ there are also $k$-vectors which are not blades and do not represent a subspace of $\R{n}$.  Such vectors occur as bivectors when $n=4$ and play an important role in the discussion of kinematics and dynamics see \Sec{sec:rbm}. 

\mysubsubsection{Projectivized exterior algebra} The exterior algebra can be projectivized using the same process defined above for the construction of $\RP{n}$ from $\R{n+1}$, but applied to the vector spaces $\bigwedge^k (\R{n+1})$. This  yields the projectivized exterior algebra $W := \proj{\bigwedge (\R{n+1})}$. 
The operations of $\bigwedge (\R{n+1})$ carry over to $W$, since, roughly speaking:  \quot{Projectivization commutes with outer product.}  The difference lies in how the elements and operations are projectively interpreted. 
The $k$-blades of $W$ correspond to $(k-1)$-dimensional subspaces of $\RP{n}$.  All multiples of the same $k$-blade  represent the same projective subspace, and differ only by intensity (\cite{whitehead98}, \S 16-17).  1-blades correspond to points; 2-blades to lines; 3-blades to planes, etc.  

\mysubsubsection{Dual exterior algebra} The dual algebra $W^* := \pdgrass{(n+1)}$  is formed by projectivizing the exterior algebra of the dual vector space $(\R{n+1})^*$.
Details can be found in the excellent Wikipedia article \cite{wikiextalg}, based on \cite{bourbaki89}. $W^*$ is the alternating algebra of $k$-multilinear forms, and is naturally isomorphic to $W$; again, the difference lies in how the elements and operations are interpreted.  Like $W$, $W^*$ represents the subspace structure of $\RP{n}$, but turned on its head: 1-blades represent projective hyperplanes, while $n$-blades represent projective points. The outer product $\vec{a} \wedge \vec{b}$ corresponds to the \emph{meet} rather than \emph{join} operator.   In order to distinguish the two outer products of $W$ and $W^*$, we write the outer product in $W$ as $\vee$, and leave the outer product in $W^*$ as $\wedge$.  These symbols match closely the affiliated operations of join (union $\cup$) and meet (intersection $\cap$), resp. 

\subsection{Remarks on homogeneous coordinates}
\label{sec:homcoord}

We use the terms \emph{homogeneous} model and \emph{projective} model interchangeably, to denote the projectivized version of Grassmann (and, later, Clifford) algebra. 

The projective model allows a certain freedom in specifying results within the algebra. In particular, when the calculated quantity is a subspace, then the answer is only defined up to a non-zero scalar multiple.  In some literature, this fact is represented by always surrounding an expression $x$ in square brackets $[x]$ when one means \quot{the projective element  corresponding to the vector space element $x$}. We do not adhere to this level of rigor here, since in most cases the intention is clear. 

Some of the formulas introduced below take on a simpler form which take advantage of this freedom, but they may appear unfamiliar to those used to working in the more strict vector-space environment.  On the other hand, when the discussion later turns to kinematics and dynamics, then this projective equivalence is no longer strictly valid. \mymarginpar{Make sure this is picked up again in the RBM section.} Different representatives of the same subspace represent weaker or stronger instances of a velocity or momentum (to mention two possibilities). In such situations terms such as \emph{weighted} point or \quot{point with intensity} will be used.  See  \cite{whitehead98}, Book III, Ch. 4.

\subsection{Equal rights for $W$ and $W^*$}
\label{sec:equalrights}
From the point of view of representing $\RP{n}$,  $W$ and $W^*$ are  equivalent.  There is no \emph{a priori} reason to prefer one to the other.  Every geometric element in one algebra occurs in the other, and any configuration in one algebra has a dual configuration in the other obtained by applying the Principle of Duality \cite{cox:pg}, to the configuration. We refer to $W$ as a \emph{point-based}, and $W^*$ as a \emph{plane-based}, algebra.\versions{\footnote{We prefer the dimension-dependent formulation \emph{plane-based} to the more precise \emph{hyperplane-based}.  We also prefer not to refer to the plane-based algebra as the \emph{dual} algebra, since this formulation depends on the accident that the original algebra is interpreted as point-based.}}{}

Depending on the context, one or the other of the two algebras may be more useful.   
Here are some examples:
\begin{itemize}
\item  \textbf{Joins and meets.} $W$ is the natural choice to calculate subspace joins, and $W^*$, to calculate subspace meets. See \Sec{sec:joinmeet}.
\item \textbf{Spears and axes.} Lines appear in two aspects: as {\beam}s (bivectors in $W$) and axes (bivectors in $W^*$).  See \Sec{sec:nolines}.
\item \textbf{Euclidean geometry.} $W^*$ is the correct choice to use for modeling euclidean geometry.  See \Sec{sec:rightmodel}.
\item \textbf{Reflections in planes.} $W^*$ has advantages for kinematics, since it naturally allows building up rotations as products of reflections in planes. See \Sec{sec:eucsand2d}.
\end{itemize}

\mysubsubsection{Bases and isomorphisms for $W$ and  $W^*$}
\label{sec:wwstar}
Our treatment differs from other approaches (for example, Grassmann-Cayley algebras) in explicitly maintaining both algebras on an equal footing rather than expressing the wedge product in one in terms of the wedge product of the other (as in the Grassman-Cayley \emph{shuffle} product) (\cite{selig05}, \cite{perwass09}). To switch back and forth between the two algebras, we construct an algebra isomorphism that, given an element of one algebra, produces the element of the second algebra which corresponds to the same geometric entity of $\RP{n}$.  We show how this works for the case of interest $n=3$.  

\begin{figure}[b]
\sidecaption
\includegraphics[width=.62\columnwidth]{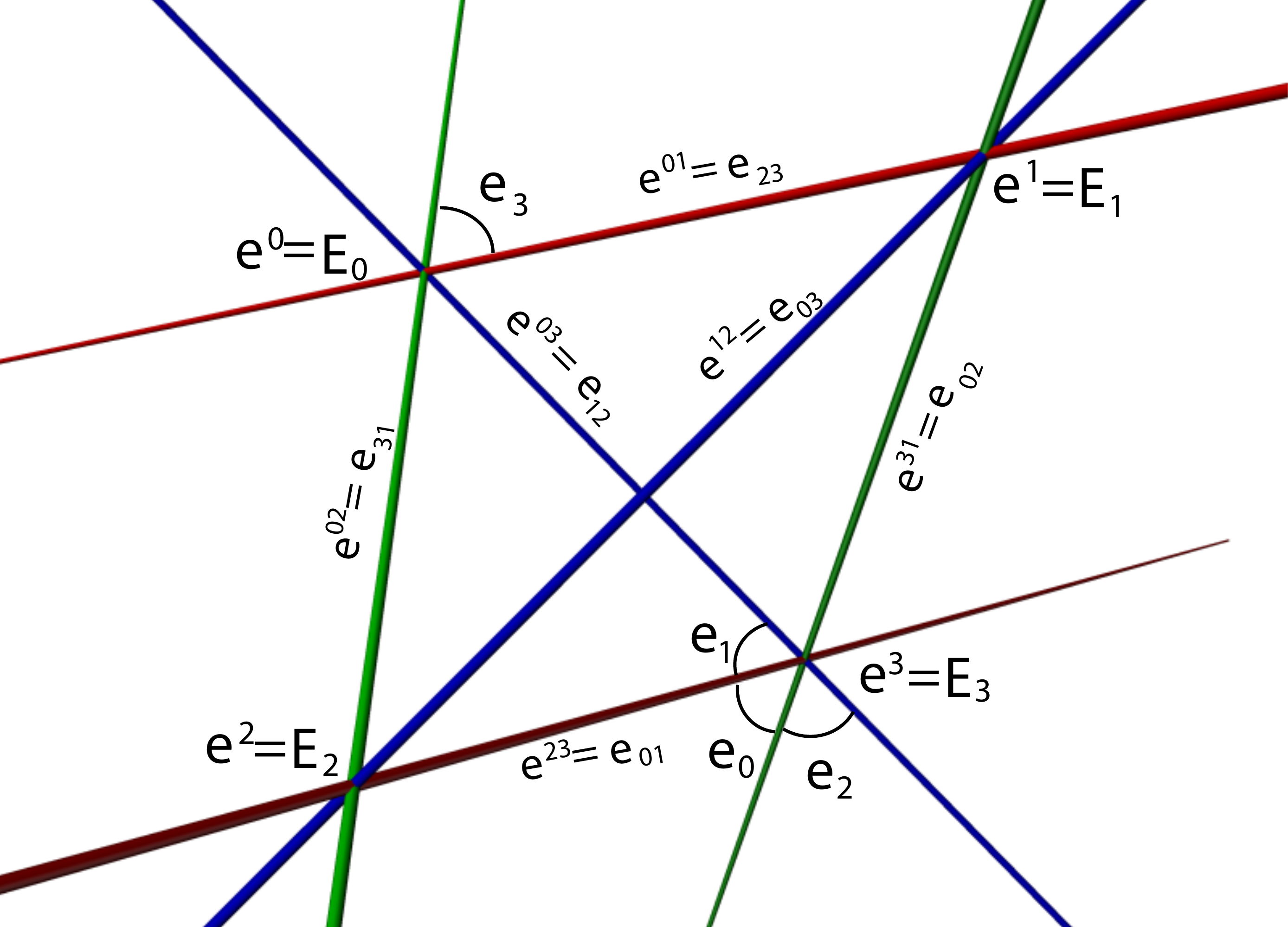}
\caption{Fundamental tetrahedron with dual labeling.  Entities in $W$ have superscripts; entities in $W^*$ have subscripts. Planes are identified by labeled angles of two spanning lines. A representative sampling of equivalent elements is shown.}
\label{fig:fundTetra}
\end{figure}

\mysubsubsection{The isomorphism $\mathbf{J}$} Each weighted subspace $S$ of $\RP{3}$ corresponds to a unique element $S_W$ of $W$ and to a unique element  $S_{W^*}$ of $W^*$.   We seek a bijection $\mathbf{J}: W \rightarrow W^*$ such that $J(S_W) = S_{W^*}$.   If we have found $\mathbf{J}$ for the basis $k$-blades, then it extends by linearity to multivectors.   To that end,  we introduce a basis for $\R{4}$ and extend it to a basis for $W$ and $W^*$ so that $\mathbf{J}$ takes a particularly simple form. Refer to  \Fig{fig:fundTetra}. 

\mysubsubsection{The canonical basis} A basis $\{\vec{e}^0, \vec{e}^1, \vec{e}^2, \vec{e}^3\}$ of $\R{4}$ corresponds to a coordinate tetrahedron for $\RP{3}$, with corners occupied by the basis elements\footnote{We use superscripts for $W$ and subscripts for $W^*$ since $W^*$ will be the more important algebra for our purposes.}.  Use the same names to identify the elements of $P(\bigwedge^1 (\R{4}))$ which correspond to these projective points.  Further, let  $\vec{I}^0 := \vec{e}^0 \wedge \vec{e}^1 \wedge \vec{e}^2 \wedge \vec{e}^3$ be the basis element of $P(\bigwedge^4 (\R{4}))$, and $\vec{1}^0$ be the basis element for $P(\bigwedge^0 (\R{4}))$.   Let the basis for  $P(\bigwedge^2 (\R{4}))$ be given by the six edges of the tetrahedron:
\begin{eqnarray*}
\{\vec{e}^{01}, \vec{e}^{02}, \vec{e}^{03}, \vec{e}^{12}, \vec{e}^{31}, \vec{e}^{23}\}
\end{eqnarray*}
where $\vec{e}^{ij} := \vec{e}^i \wedge \vec{e}^j$ represents the oriented line joining $\vec{e}^i$ and $\vec{e}^j$.\footnote{Note that the orientation of $\vec{e}^{31}$ is reversed; this is traditional since Pl\"{u}cker introduced these line coordinates.} 
Finally, choose a basis $\{\vec{E}^0, \vec{E}^1, \vec{E}^2, \vec{E}^3\}$ for $P(\bigwedge^3 (\R{4}))$  satisfying the condition that $\vec{e}^i \vee \vec{E}^i = \vec{I}^0$. This corresponds to choosing the $i^{th}$ basis 3-vector to be the plane opposite the $i^{th}$ basis 1-vector in the fundamental tetrahedron, oriented in a consistent way.  

We repeat the process for the algebra $W^*$,  writing indices as subscripts.  Choose the basis 1-vector $\vec{e}_i$ of $W^*$ to represent the same plane as $\vec{E}^i$.  That is, $\mathbf{J}(\vec{E}^i) = \vec{e}_i$.   Let $\vec{I}_0 := \vec{e}_0 \wedge \vec{e}_1 \wedge \vec{e}_2 \wedge \vec{e}_3$ be the pseudoscalar of the algebra.   Construct bases for grade-0, grade-2, and grade-3 using the same rules as above for $W$ (i. e., replacing subscripts by superscripts).   The results are represented in Table \ref{tab:wwstar}.

\begin{table}
\begin{center}
\begin{tabular} {|  l  |  c | c |} \hline 
feature & $W$ & $W^*$ \\ \hline \hline
0-vector & scalar $ \mathbf{1}^0$ & scalar $\mathbf{1}_0$ \\ \hline
vector & point $ \{e^i\}$ & plane $\{e_i\}$ \\ \hline
bivector & \quot{{\beam}} $\{e^{ij}\}$ & \quot{axis} $\{e_{ij}\}$ \\ \hline
trivector & plane $\{E^i\}$ & point $\{E_i\}$ \\ \hline
4-vector & $\mathbf{I}^0$ & $\mathbf{I}_0$ \\ \hline
outer product & {\color{black} join $\vee$} & {\color{black} meet $\wedge$} \\ \hline
\end{tabular}
\end{center}
\label{tab:wwstar}
\caption{Comparison of $W$ and $W^*$.}
\end{table}
Given this choice of bases for $W$ and $W^*$, examination of \Fig{fig:fundTetra} makes clear that, on the basis elements, $\mathbf{J}$ takes the following simple form: 
\begin{equation} \label{eqn:jdef}
\mathbf{J}(e^i) := E_i, ~~~~\mathbf{J}(E^i) :=  e_i, ~~~~\mathbf{J}(e^{ij}) := e_{kl}   
\end{equation}
where in the last equation, $(ijkl)$ an even permutation of $(0123)$.

\begin{figure}
\sidecaption
\includegraphics[width = .6\columnwidth]{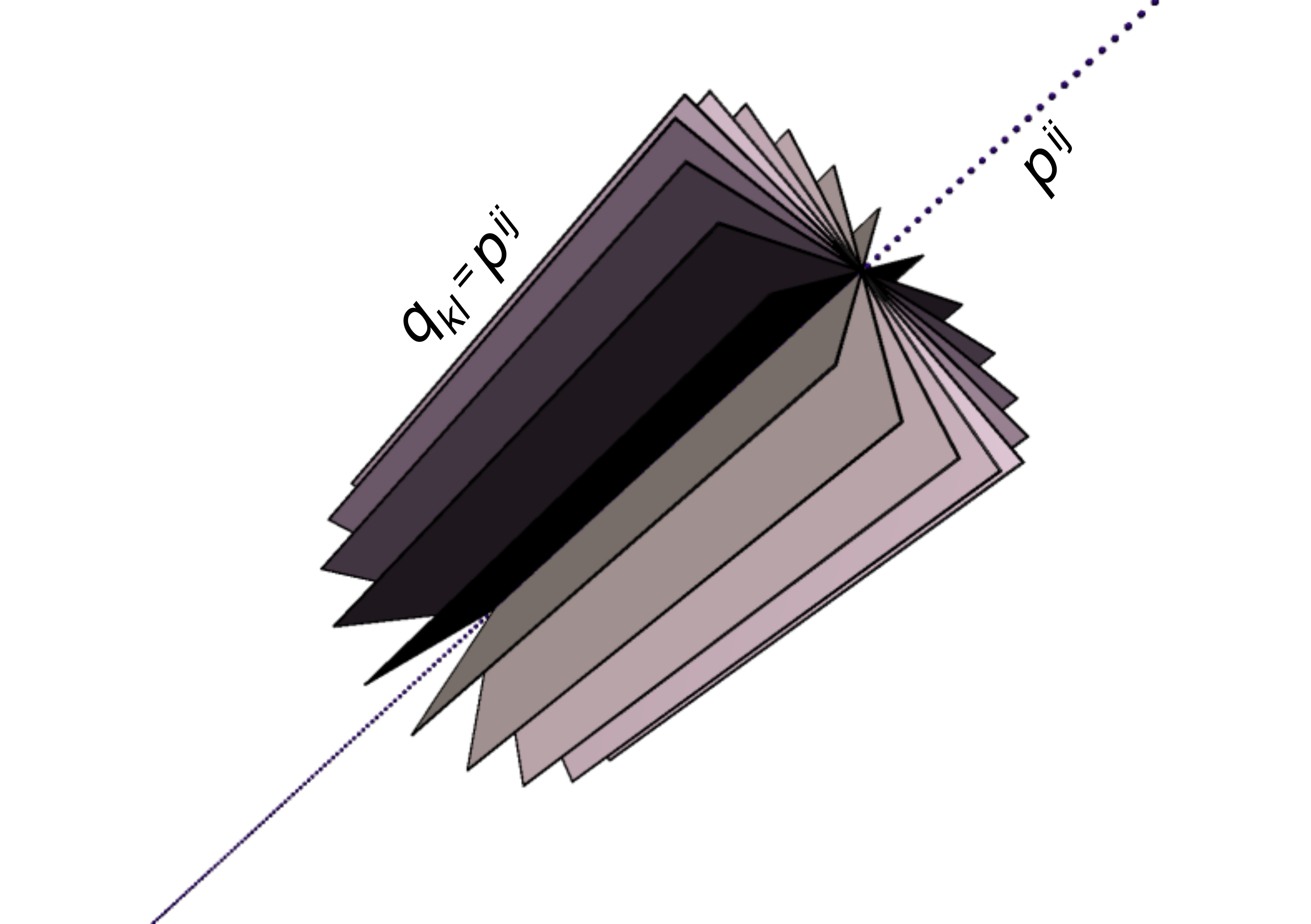}
\caption{A line in its dual nature as {\beam}, or point range; and as axis, or plane pencil.}
\label{fig:spearaxis}
\end{figure}

Furthermore, $\mathbf{J}(\mathbf{1}^0) = \mathbf{I}_0$ and $\mathbf{J}(\mathbf{I}^0) = \mathbf{1}_0$ since these grades are one-dimensional. To sum up: the map $\mathbf{J}$ is grade-reversing and, considered as a map of coordinate-tuples, it is the identity map on all grades except for bivectors.  What happens for bivectors?  In $W$, consider $\vec{e}^{01}$, the joining line of points $\vec{e}^0$ and $\vec{e}^1$ (refer to \Fig{fig:fundTetra}). In $W^*$, the same line is $\vec{e}_{23}$, the intersection of the only two planes which contain both of these points, $\vec{e}_2$ and $\vec{e}_3$.
\versions{On a general bivector $\mathbf{J}$ takes the form:
\begin{eqnarray*}
\mathbf{J}(a_{01} e^{01} + a_{02} e^{02} + a_{03} e^{03} + a_{12} e^{12} + a_{31} e^{31} + a_{23} e^{23}) = \\
a_{23} e_{01} + a_{31} e_{02} + a_{12} e_{03} + a_{03} e_{12} + a_{02} e_{31} + a_{01} e_{23}
\end{eqnarray*}
The coordinate-tuple is reversed! }{}
See \Fig{fig:spearaxis}. 
 Since $\mathbf{J}^{-1}$ is obtained from the definition of $\mathbf{J}$ by swapping superscripts and subscripts, we can consider $\mathbf{J}: W \leftrightarrow W^*$ as a defined on both algebras,  with $\mathbf{J}^2$ the identity.  The full significance of $\mathbf{J}$ will only become evident after metrics are introduced (\Sec{sec:jandmp}).   
 \versions{
Appendix 1 contains a detailed description of how $\mathbf{J}$ is constructed in the $n$-dimensional case, and its close relation to the use of a nondegenerate metric to achieve the same goal. 
 }{} 
 We now show how to use $\mathbf{J}$ to define meet and join operators valid for both $W$ and $W^*$.

\mysubsubsection{Projective join and meet}
\label{sec:joinmeet}
Knowledge of $\mathbf{J}$ allows equal access to join and meet operations.  We define a meet operation $\wedge$ for two blades $A, B \in W$:
\begin{equation}  \label{eqn:joinmeet}
A \wedge B  =  \mathbf{J}( \mathbf{J}(A) \wedge \mathbf{J}(B))
\end{equation}
and extend by linearity to the whole algebra.
There is a similar expression for the join $\vee$ operation for two blades $A, B \in W^*$:
\begin{equation}
A \vee B :=  \mathbf{J}( \mathbf{J}(A) \vee \mathbf{J}(B))
\end{equation}

We turn now to another feature highlighting the importance of maintaining $W$ and $W^*$ as equal citizens.

\textbf{There are no lines, only {\beam}s and axes!}
\label{sec:nolines}
Given two points $\vec{x}$ and $\vec{y}$ $\in W$,  the condition that a third point $\vec{z}$ lies in the subspace spanned by the 2-blade $\vec{l} := \vec{x} \vee \vec{y}$ is that $\vec{x} \vee \vec{y} \vee \vec{z} = 0$, which implies that $\vec{z} = \alpha \vec{x} + \beta \vec{y}$ for some $\alpha, \beta$ not both zero.  In projective geometry, such a set is called a \emph{point range}.  We prefer the more colorful term \emph{\beam}. Dually, given two planes $\vec{x}$ and $\vec{y}$ $\in W^*$,  the condition that a third plane $\vec{z}$ passes through the subspace spanned by the 2-blade $\vec{l} := \vec{x} \wedge \vec{y}$ is that $\vec{z} = \alpha \vec{x} + \beta \vec{y}$.  In projective geometry, such a set is called a \emph{plane pencil}.  We prefer the more colorful term \emph{axis}.  

Within the context of $W$ and $W*$, lines exist only in one of these two aspects:  of {\beam} -- as bivector in $W$ -- and axis -- as bivector in $W^*$. 
This naturally generalizes to non-simple bivectors: there are point-wise bivectors (in $W$), and plane-wise bivectors (in $W^*$.)  Many of the important operators of geometry and dynamics we will meet below, such as the polarity on the metric quadric (\Sec{sec:ckc}), and the inertia tensor of a rigid body (\Sec{sec:rbm2}), map $\langle W\rangle_2$ to $\langle W^*\rangle_2$ and hence map {\beam}s to axes and vice-versa.  
Having both algebras on hand preserves the qualitative difference between these dual aspects of the generic term \quot{line}.

We now proceed to describe how to introduce metric relations.
\section{Clifford algebra for euclidean geometry}
\label{sec:cafeg}

The outer product is anti-symmetric, so $\vec{x} \wedge \vec{x} =0$. However, in geometry there are important bilinear products which are symmetric. We introduce a  real-valued \emph{inner product} on pairs of vectors $\vec{x} \cdot \vec{y}$ which is a real-valued symmetric bilinear map.  Then, the geometric product on 1-vectors is defined as the sum of the inner and outer products:
\[
\vec{x} \vec{y} := \vec{x} \cdot \vec{y} + \vec{x} \wedge \vec{y}
\]
How this definition can be extended to the full exterior algebra is described elsewhere (\cite{dfm07}, \cite{hessob87}).  The resulting algebraic structure is called a \emph{real Clifford algebra}.  It is fully determined by its signature, which describes the inner product structure.  The signature is a  triple of integers $(p,n,z)$ where $p+n+z$ is the dimension of the underlying vector space, and $p$, $n$, and $z$ are the numbers of positive, negative, and zero entries along the diagonal of the quadratic form representing the inner product.  We denote the corresponding Clifford algebra constructed on the point-based Grassmann algebra as $\pclal{p}{n}{z}$; that based on the plane-based Grassmann algebra, as $\pdclal{p}{n}{z}$.\mymarginpar{Is there a standard notation for this distinction?}

The discovery and application of signatures to create different sorts of metric spaces within projective space goes back to a technique invented by Arthur Cayley and developed by Felix Klein \cite{klein:neg}.  The so-called \emph{Cayley-Klein construction}  provides models of the three standard metric geometries (hyperbolic, elliptic, and euclidean) -- along with many others! -- within projective space.  This work provides the mathematical foundation for the inner product as it appears within the homogeneous model of Clifford algebra. Since the Cayley-Klein construction for euclidean space has some subtle points, is relatively sparsely represented in the current literature, and is crucial to what follows, we describe it below.

\subsection{The Cayley-Klein Construction}
\label{sec:ckc}


For simplicity we focus on the case $n=3$. To obtain metric spaces inside $\RP{3}$ begin with a symmetric bilinear form $Q$ on $\R{4}$. The quadric surface associated to $Q$ is then defined to be the points  $\{\mathbf{x}\mid Q(\vec{x}, \vec{x}) = 0\}$. For nondegenerate $Q$,  a distance between points $A$ and $B$ can be defined by considering the cross ratio of the four points $A$, $B$, and the two intersections of the line $AB$ with $Q$.  Such a $Q$ is characterized by its signature, there are two cases of interest for $n=3$: $(4,0,0)$ yielding elliptic geometry and $(3,1,0)$ yielding hyperbolic geometry. These are point-based metrics; they induce an inner product on planes, which one can show is identical to the original signatures.  By interpolating between these two cases, one is led to the degenerate case in which the quadric surface collapses to a plane, or to a point.  In the first case, one obtains euclidean geometry; the plane is called the \emph{ideal} plane.  The signature breaks into two parts: for points, it's $(1,0,3)$ and for planes it's $(3,0,1)$.  The distance function for euclidean geometry is based on a related limiting process. For details see \cite{klein:neg} or \fullversion. \textbf{Warning}: in the projective model, the signature $(n,0,0)$ is called the \emph{elliptic} metric, and \emph{euclidean} metric refers to  these degenerate signatures.

\versions{

\mysubsubsection{Details of the limiting process}
We restrict attention to a special family of $Q$ parametrized by
a real parameter $\epsilon$, and define an inner product between points $\vec{X} = (x_0,x_1,x_2,x_3$ and $\vec{Y} = (y_0, y_1,y_2,y_3)$ as follows:
\vspace{-.1in}
\begin{eqnarray}
\langle \vec{X}, \vec{Z}  \rangle_\epsilon := Q_\epsilon(\vec{X} , \vec{Y})  :=  \epsilon x_0 y_0 + x_1 y_1 +x_2 y_2 + x_3 y_3
\end{eqnarray}

The inner product is not well-defined on $\RP{3}$ since different choices of representatives for the arguments will yield multiples of the inner product.  But the expression
\begin{equation} \label{eqn:hommetric}
\frac{\langle \vec{X}, \vec{Y}  \rangle_\epsilon^2}{\langle \vec{X}, \vec{X}  \rangle_\epsilon \langle \vec{Y}, \vec{Y}  \rangle_\epsilon}
\end{equation}
is  homogeneous in its arguments, hence well-defined (assuming neither argument is a null vector).  It is this normalized inner product which appears in the metric formulae below.

$\langle ,  \rangle_1$ gives the inner product for the signature $(4,0,0)$ (elliptic space), and $\langle ,  \rangle_{-1}$ for $(3,1,0)$ (hyperbolic space).  For brevity we write these two inner products as $\langle , \rangle_+$ and $\langle , \rangle_-$, resp.  

The inner products above are defined on the \emph{points} of space; there is an induced inner product on the \emph{planes} of space formally defined as the adjoint operator of the operator $Q_\epsilon$.  One can easily show that for two planes $\vec{u}$ and $\vec{v}$ this takes the form : 
\begin{eqnarray}
\langle \vec{u}, \vec{v}  \rangle_\epsilon := Q^*_\epsilon(\vec{u} , \vec{y})  &:=& \epsilon^{-1} u_0 v_0 + u_1 v_1 +u_2 v_2 + u_3 v_3 \\
&\simeq& u_0 v_0 +  \epsilon (u_1 v_1 +u_2 v_2 + u_3 v_3)
\end{eqnarray}
where the $\simeq$ symbol denotes projective equivalence obtained by multiplying the first equation by $\epsilon \neq 0$. 
For $\epsilon \in \{1,-1\}$ the resulting signatures are equivalent to the point signatures; we can use the same notation for both points and planes in the formulae below.

The metric quadric for elliptic space is  $\{\mathbf{x} \mid \langle \mathbf{x}, \mathbf{x} \rangle_+ = 0\}$.  There are no real solutions; this is called a totally imaginary quadric.  The points $\mathbf{x}  \mid \langle \mathbf{x}, \mathbf{x} \rangle_+ > 0 $ constitute the projective model of elliptic space: all of $\RP{3}$.  
The metric quadric for hyperbolic space is $\{\mathbf{x}  \mid \langle \mathbf{x}, \mathbf{x} \rangle_- = 0\}$, the unit sphere in euclidean space.  The points $\mathbf{x} \mid \langle \mathbf{x}, \mathbf{x} \rangle_- < 0 $ constitute the projective model of hyperbolic space: the interior of the euclidean unit ball. The sphere itself is sometimes called the \emph{ideal sphere} of hyperbolic space.

The euclidean metric is a limiting case of the above family of metrics as $\epsilon$ grows larger and larger, from either the positive or negative side.   For $\lim_{\epsilon \rightarrow \infty}$, we arrive at the euclidean metric.   Due to the homogeneity of  \Eq{eqn:hommetric} we allow ourselves to apply arbitrary non-zero scaling factors to the inner product.  For points the limiting process yields: 
\begin{eqnarray*}
\lim_{\epsilon \rightarrow \infty}\langle \vec{X}, \vec{Y}  \rangle_\epsilon &=& \lim_{\epsilon \rightarrow \infty} (\epsilon x_0 y_0 + x_1 y_1 +x_2 y_2 + x_3 y_3)\\
&\simeq& \lim_{\epsilon \rightarrow \infty}(x_0 y_0 + \frac{1}{\epsilon}(x_1 y_1 +x_2 y_2 + x_3 y_3))\\
&=& x_0 y_0
\end{eqnarray*}
where the second equation results from the first by multiplication by $\frac1{\epsilon}$.
This is the signature $(1,0,3)$.   For planes $\vec{u}$ and $\vec{v}$, one works with the adjoint form:
\begin{eqnarray*}
\lim_{\epsilon \rightarrow \infty}(\langle \vec{u}, \vec{v}  \rangle_\epsilon)  &=& \lim_{\epsilon \rightarrow \infty} (u_0 v_0 + \epsilon(u_1 v_1 +u_2 v_2 +u_3 v_3)) \\
&\simeq&  \lim_{\epsilon \rightarrow \infty}( \frac{1}{\epsilon} u_0 v_0 + u_1 v_1 +u_2 v_2 + u_3 v_3)\\
&=& u_1 v_1 +u_2 v_2 + u_3 v_3
\end{eqnarray*}
which corresponds to the signature $(3,0,1)$.  
The point- and plane-signatures for euclidean space are complementary!

The projective model of euclidean space consists of $\RP{3}$ with the plane $x_0=0$ removed, the so-called \emph{ideal plane}.  For the points of this plane, the original inner product remains valid. These are euclidean free vectors, see discussion below.

\mysubsubsection{Distance and Angle Formulae}
Let $\vec{X}$ and $\vec{Y}$ be two points, and $l_{XY}$ their joining line.  In general $l_{XY}$ will have two (possibly imaginary) intersection points $\vec{Q}_1$ and $\vec{Q}_2$ with the metric quadric.
The original definition of the distance of two points $\vec{X}$ and $\vec{Y}$ in these noneuclidean spaces relied on logarithm of the cross ratio of the points $\vec{X}$, $\vec{Y}$, $\vec{Q}_1, $, and $ \vec{Q}_2$ (\cite{cox:neg}).  By straightforward functional identities these formulae can be brought into alternative form. 
The distance $d$ between two points $\mathbf{X}$ and $\mathbf{Y}$ in the  elliptic (resp. hyperbolic) metric is then given by:
\begin{eqnarray} \label{eq:elldist}
d& = & \cos^{-1}(\frac{\langle \mathbf{X}, \mathbf{Y} \rangle_+}{\sqrt(\langle \mathbf{X}, \mathbf{X} \rangle_+\langle \mathbf{Y}, \mathbf{Y} \rangle_+)}) \\
d & = & \cosh^{-1}(\frac{-\langle \mathbf{X}, \mathbf{Y} \rangle_-}{\sqrt(\langle \mathbf{X}, \mathbf{X} \rangle_- \langle \mathbf{Y}, \mathbf{Y} \rangle_-)} ) 
\end{eqnarray}

The familiar euclidean distance between two (non-homogeneous) points
\begin{eqnarray*}
d& = & \sqrt{(x_1-y_1)^2 + (x_2-y_2)^2 +  (x_3 - y_3)^2 }
\end{eqnarray*}
 can be derived by parametrizing the above formulas by $epsilon$, and evaluating a carefully chosen limit as as $\epsilon \rightarrow \infty$. 
 To simplify the treatment, we work with a projective line ($n=1$), and introduce the family of inner products parametrized by the real parameter $\epsilon$ as above.  We want to produce a euclidean metric on this line in which the basis vector $\vec{e}_1$ is the ideal point.  We consider two points $\vec{X} = (x_0, x_1)$ and $\vec{Y} = (y_0, y_1)$.  We can assume that $x_0 \neq 0$ and $y_0 \neq 0$, since $\vec{e}_1$ is the ideal point and doesn't belong to the euclidean line.  Choose the projective representative so that $x_0 = y_0 = 1$.  Then:
 \begin{eqnarray*}
 \langle \vec{X}, \vec{Y} \rangle_\epsilon := \epsilon x_0 y_0 + x_1 y_1 = \epsilon+ x_1 y_1
 \end{eqnarray*}
By \Eq{eq:elldist}, the distance function $d_\epsilon$ associated to $\langle,\rangle_\epsilon$ is determined by:
 \begin{equation} \label{eq:epsmetric}
 \cos{d_\epsilon(\vec{X},\vec{Y})} = \frac{\langle \vec{X},\vec{Y} \rangle_\epsilon}{\sqrt(\langle \mathbf{X}, \mathbf{X} \rangle_\epsilon\langle \mathbf{Y}, \mathbf{Y} \rangle_\epsilon)}
 \end{equation}
 Abbreviate $d_\epsilon(\vec{X},\vec{Y})$ as $d_\epsilon$.  Now consider the limit as $\epsilon\rightarrow \infty$. It's clear from  \Eq{eq:epsmetric} that $\lim_{\epsilon\rightarrow \infty}  \cos{d_\epsilon}= 1$, that is, the distance goes to zero.   Define a new distance function $\hat{d}_\epsilon := \sqrt{\epsilon}d_\epsilon$.  The scaling factor $ \sqrt{\epsilon}$ prevents the distance from going to zero in the limit. Instead one obtains:
\begin{align}
\lim_{\epsilon \rightarrow \infty} \hat{d}_\epsilon^2 &= (x_1 - y_1)^2 
 \end{align}
which is clearly equivalent to the euclidean distance between the two points $\vec{X}$ and $\vec{Y}$.  For the details consult Appendix 2.

In all three geometries the angle $\alpha$ between two oriented planes $\mathbf{u}$ and $\mathbf{v}$  is given by (where $\langle,\rangle$ represents the appropriate inner product):
\begin{eqnarray*}
\alpha & = & \cos^{-1}(\frac{\langle \mathbf{u}, \mathbf{v} \rangle}{\sqrt(\langle \mathbf{u}, \mathbf{u} \rangle\langle \mathbf{v}, \mathbf{v} \rangle)}) 
\end{eqnarray*}

}{}

\mysubsubsection{Polarity on the metric quadric}
For a $Q$ and a point $\vec{P}$, define the set $\mathbf{P}^\perp := \{\mathbf{X} \mid Q(\mathbf{X}, \mathbf{P}) = 0\}$.  When $\mathbf{P}^\perp$ is a plane, it's called the \emph{polar plane} of the point.  For a plane $\vec{a}$, there is also an associated \emph{polar point} defined analogously using the \quot{plane-based} metric. Points and planes with such polar partners are called \emph{regular}.  
 In the euclidean case,
the polar plane of every finite point is the ideal plane; the polar point of a finite plane is the ideal point in the normal direction to the plane. Ideal points and the ideal plane are not regular and have no polar partner. 
The polar plane of a point is important since it can be identified with the tangent space of the point when the metric space is considered as a differential manifold.  Many of the peculiarities of euclidean geometry may be elegantly explained due to the degenerate form of the polarity operator.  In the Clifford algebra setting, this polarity is implemented by multiplication by the pseudoscalar. 

\mysubsubsection{Free vectors and the euclidean metric} \label{sec:freevec}
As mentioned above, the tangent space at a point is the polar plane of the point. Every euclidean point shares the same polar plane, the ideal plane.  In fact, the ideal points (points of the ideal plane) can be identified with euclidean free vectors.  A model for euclidean geometry should handle both euclidean points and euclidean free vectors. This is complicated by the fact that free vectors have a natural signature $(3,0,0)$.  However, since the limiting process (in Cayley-Klein) that led to the degenerate point metric $(1,0,3)$ only effects the non-ideal points, it turns out that the original non-degenerate metric, restricted to the ideal plane, yields the desired signature $(3,0,0)$. As we'll see in \Sec{sec:enumpro} and \Sec{sec:enumpro3}, the model presented here is capable of mirroring this subtle fact.

\versions{ To learn more about the mathematics,  see \cite{klein:neg}, \cite{cox:pg}  \cite{cox:neg}, and \cite{kowol2009}.  For a modern formulation with a good collection of formulae see
 \cite{lt:gt3m}. 
}{}

\subsection{A model for euclidean geometry}
\label{sec:rightmodel}
As noted above, the euclidean inner product has signature $(1,0,3)$ on points and $(3,0,1)$ on planes.  If we attach the first signature to $W$, we have the following relations for the basis 1-vectors:
\[
(\vec{e}^0)^2 = 1; ~~(\vec{e}^1)^2 = (\vec{e}^2)^2 = (\vec{e}^3)^2 = 0
\]
It's easy to see that these relations imply that, for all basis trivectors $E_i$, $E_i^2=0$.  But the trivectors represent planes, and the signature for the plane-wise euclidean metric is $(3,0,1)$, not $(0,0,4)$.  Hence, we cannot use $W$ to arrive at euclidean space. If  instead, we begin with $W^*$, and attach the plane-wise signature $(3,0,1)$, we obtain:
\[
(\vec{e}_0)^2 = 0; ~~(\vec{e}_1)^2 = (\vec{e}_2)^2 = (\vec{e}_3)^2 = 1
\]
It is easy to check that this inner product, when extended to the higher grades, produces the proper behavior on the trivectors, since only $\vec{E}_0 = \vec{e}_1\vec{e}_2\vec{e}_3$ has non-zero square, producing the point-wise signature $(0,1,3)$ (equivalent to the signature $(1,0,3))$.  Hence, $W^*$ is the correct choice for constructing a model of euclidean geometry. 

\textbf{Counterspace.} What space \textbf{does} one obtain by attaching the signature $(3,0,1)$ to $W$? 
One obtains a different metric space, sometimes called polar-euclidean space or counterspace.  Its metric quadric is a \emph{point} along with all the \emph{planes passing through it} (dual to the euclidean ideal \emph{plane} and all the \emph{points lying in it})\footnote{Blurring the distinction between these two spaces may have led some authors to incorrect conclusions about the homogeneous model, see \cite{li08}, p. 11}. \versions{Like euclidean space, it arises as a limiting case of the Cayley-Klein construction, when one lets $\epsilon \rightarrow 0$.  Instead of flattening out the metric quadric into a plane, this limit contracts it to a point.}{} See \cite{conradt08}, pp. 71ff., for a related discussion.

We retain $W$, the point-based algebra, solely as a Grassmann algebra, primarily for calculating the join operator.   All euclidean metric operations are carried out in $W^*$. Or equivalently, we attach the metric $(0,0,4)$ to $W$, forcing all inner products to zero. Due to the more prominent role of $W^*$,   the basis element for scalar and pseudoscalar  in $W^*$ will be written without index as $\vec{1}$ and $\vec{I}$; we may even omit $\one$ when writing scalars, as is common in the literature.

\subsection{$\mathbf{J}$, metric polarity, and the regressive product }
\label{sec:jandmp}
We can now appreciate better the significance of $\mathbf{J}: W \rightarrow W^*$.  Consider the map $\mathbf{\Pi}: W \rightarrow W$ defined analogously to $\mathbf{J}$ in \Eq{eqn:jdef}:
\begin{equation} \label{eqn:pidef}
\mathbf{\Pi}(e^i) := E^i,~~~~ \mathbf{\Pi}(E^i) := e^i, ~~~~\mathbf{\Pi}(e^{ij}) := e^{kl}   
\end{equation}
$\mathbf{\Pi}$ is the same as  $\mathbf{J}$, but interpreted as a map to $W$ instead of $W^*$.  It's easy to see that $\mathbf{\Pi}$ is the polarity on the elliptic metric quadric with signature $(4,0,0)$.   Many authors (see \cite{hessob87}) define the meet operation between two blades $A, B \in W$ (also known since Grassmann as the \emph{regressive} product) via $ \vec{\Pi}(\vec{\Pi}(A) \wedge \vec{\Pi}(B))$, where $\wedge$ is the exterior product in $W$.  One can define a similar join operator in $W^*$.  We prefer to use $\mathbf{J} $ for this purpose (see \Eq{eqn:joinmeet}) since it provides a \emph{projective} solution for a \emph{projective} (incidence) problem,  and it is useful on its own (see for example \Sec{sec:rbm2}), while $\mathbf{\Pi}$, being a foreign entity, must always appear in the second power so that it has no side-effects.  To distinguish the two approaches, we suggest calling $\mathbf{\Pi}$ the \emph{metric polarity} and $\mathbf{J}$, the \emph{duality} operator, consistent with mathematical literature. For an $n$-dimensional discussion and proof, see  \versions{Appendix 1}{Appendix 1 of \fullversion}.

\section{The euclidean plane via $\pdclal{2}{0}{1}$} 
\label{sec:eucplane}
Due to the combination of unfamiliar concepts involved in the algebras $\pdclal{n}{0}{1}$ -- notably the dual construction and the degenerate metric -- we  begin our study with the Clifford algebra for the euclidean plane:  $\pdclal{2}{0}{1}$.  Then, when we turn to the 3D case, we can focus on the special challenges which it presents, notably the existence of non-simple bivectors.  A basis for the full algebra of  $\pdclal{2}{0}{1}$ is given by 
\[
\{\vec{1} := \vec{1}_0, \vec{e}_0, \vec{e}_1, \vec{e}_2,\vec{E}_0 := \vec{\e{1}\e{2}}, \vec{E}_1 := \vec{\e{2}\e{0}}, \vec{E}_2 := \vec{\e{0}\e{1}}, \vec{I} := \vec{e}_0 \vec{e}_1 \vec{e}_2\}
\]
with the relations $
\{\vec{e}_0^2 = 0; ~~\vec{e}_1^2 = \vec{e}_2^2 = 1\}
$.  
\versions{See \Fig{fig:projplane}.

\begin{figure}[t]
\sidecaption
\includegraphics[width=2.5in]{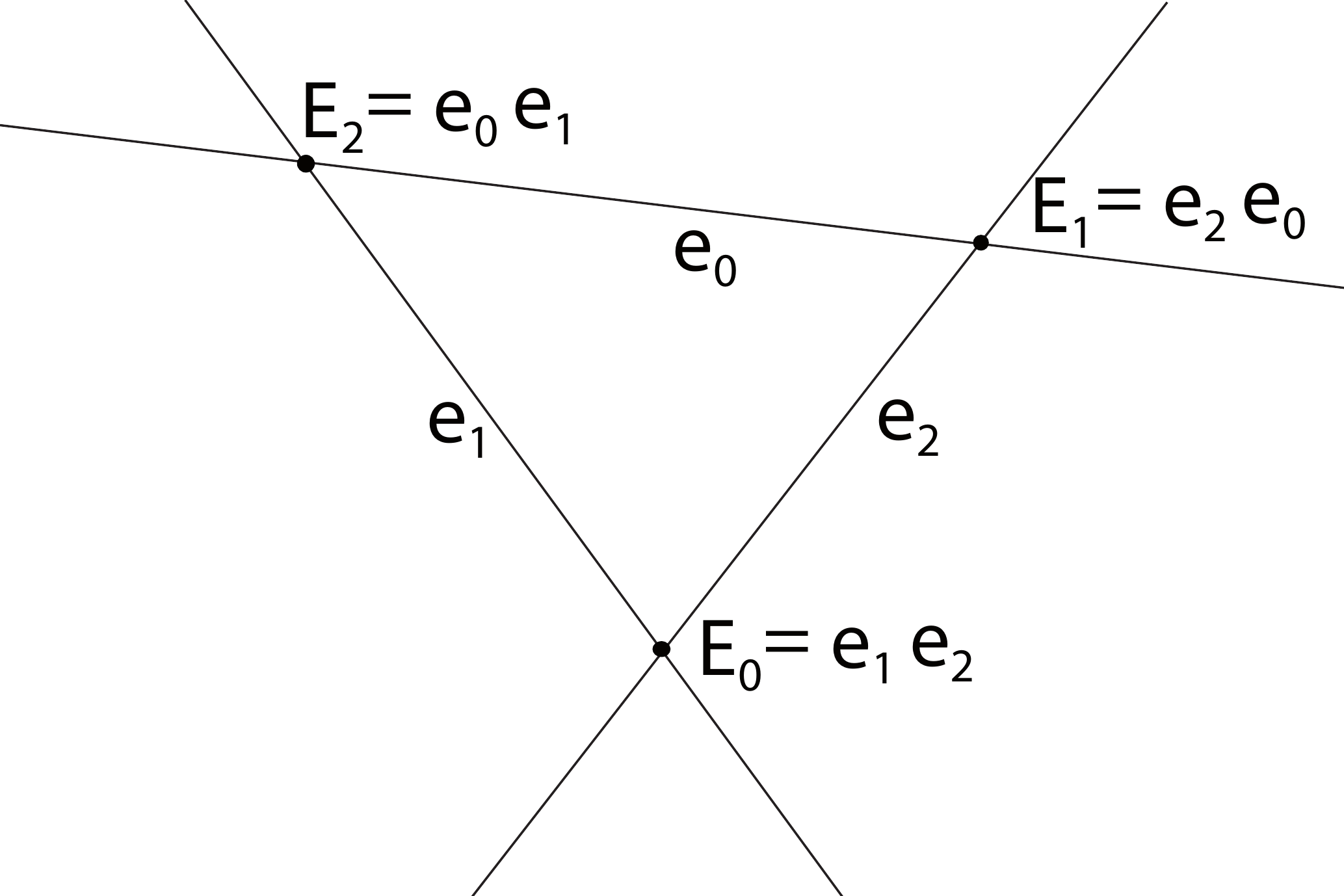}
\caption{The basis 1- and 2-vectors for  $\pdclal{2}{0}{1}$. $\vec{e}_0$ is the ideal line, $\EE{0}$ is the origin.}
\label{fig:projplane}
\end{figure}
}{}

\mysubsubsection{Consequences of degeneracy} The pseudoscalar $\eye$ satisfies $\eye^2=0$. Hence, $\eye^{-1}$ is not defined.  Many standard formulas of geometric algebra  are, however,  typically stated using $\eye^{-1}$ ( \cite{dfm07, hessob87}), since that can simplify things for nondegenerate metrics.   As explained in \Sec{sec:homcoord},  many formulas remain projectively valid when $\eye^{-1}$ is replaced by $\eye$; in such cases this is the solution we adopt. 

\mysubsubsection{Notation} We denote 1-vectors with bold small letters, and 2-vectors with bold capital letters. We will use the term {\emph{ideal} to refer to geometric elements contained in projective space but not in euclidean space.  Then $\e{0}$ is the \emph{ideal} line of the plane, $\e{1}$ is the line $x = 0$ and $\e{2}$, the line $y = 0$.  $\EE{0}$ is the origin $(1,0,0)$ while $\EE{1}$ and $\EE{2}$ are the \emph{ideal} points in the $x-$ and $y-$ direction, resp.  Points and lines which are not ideal, are called \emph{finite}, or \emph{euclidean}. 

The natural embedding of a euclidean \emph{position} $\mathbf{x} = (x,y)$ we write as $\mathbf{i}(\mathbf{x}) = \EE{0}+x\EE{1} + y\EE{2}$.  A euclidean \emph{vector} $\mathbf{v} = (x,y)$ corresponds to an ideal point (see \Sec{sec:freevec}); we denote its embedding with the same symbol $\mathbf{i}(\mathbf{v}) = x\EE{1} + y\EE{2}$.  We sometimes refer to such an element as 
a \emph{free} vector.  Conversely, a bivector $w\EE{0}+x\EE{1} + y\EE{2}$ with $w \neq 0$  corresponds to the euclidean point $(\frac{x}{w}, \frac{y}{w})$.     We refer to $w$ as the \emph{intensity} or \emph{weight} of the bivector.  And, we write $\underline{\vec{A}}$ to refer to $\mathbf{i}^{-1}(\vec{A})$.  The line $ax + by + c=0$ maps to the 1-vector $c\e{0} + a\e{1} + b\e{2}$.  A line is euclidean if and only if $a^2 + b^2 \neq 0$. 

The multiplication table is shown in Table \ref{tab:cl201}. Inspection of the table reveals that the geometric product of a $k$- and $l$-vector yields a product that involves at most two grades.  When these two grades are $|k-l|$ and $k+l$, we can write the geometric product for 2 arbitrary blades $\vec{A}$ and $\vec{B}$ as
\[
\vec{A} \vec{B} = \vec{A} \inpro \vec{B} + \vec{A} \wedge \vec{B}
\]
where $\inpro$ is the generalized inner product, defined to be $\langle \vec{A} \vec{B}\rangle_{|k-l|}$(\cite{hessob87}).  The only exception is $(l,k)=(2,2)$ where the grades 
$|k-l|=0$ and $|k-l|+2 = 2$ occur.  Following \cite{hessob87} we write the grade-2 part as:
\[ 
\vec{A} \times \vec{B} := \langle \vec{A} \vec{B} \rangle_2 = \frac12 ( \vec{A} \vec{B} -  \vec{B} \vec{A})
\]
\correction
where $\vec{A}$ and $\vec{B}$ are bivectors.  This is called the \emph{commutator} product.   
Since all vectors in the algebra are blades, the above decompositions are valid for the product of any two vectors in our algebra. 
\begin{table}[t]
\label{tab:cl201}
\begin{center}
\renewcommand{\arraystretch}{1.1}
\begin{tabularx}{.6\columnwidth} {| Y  || Y | Y  |  Y  | Y | Y | Y | Y | Y  |} \hline 
           & $\one$ & $\e{0}$ & $\e{1}$ & $\e{2}$ & $\EE{0}$ & $\EE{1}$ & $\EE{2}$ & $\eye$  \\ \hline \hline
$\one$        & $\one$ & $\e{0}$ & $\e{1}$ & $\e{2}$ & $\EE{0}$ & $\EE{1}$ & $\EE{2}$ & $\eye$  \\ \hline
$\e{0}$  & $\e{0}$ & $0$     & $\EE{2}$  & $-\EE{1}$ & $\eye$    & $0$      & $0$          & $0      $  \\ \hline
$\e{1}$  & $\e{1}$ & $-\EE{2}$ & $\one$& $\EE{0}$ & $\e{2}$ & $\eye$ & $-\e{0}$ & $\EE{1}$ \\ \hline
$\e{2}$  & $\e{2} $  & $\EE{1}$ & $-\EE{0}$ & $\one$ & $-\e{1}$ & $\e{0}$ & $\eye$ & $\EE{2}$ \\ \hline
$\EE{0}$  & $\EE{0}$ & $\eye$ & $-\e{2}$   & $\e{1}$   & $-\one$ & $-\EE{2}$ & $\EE{1}$ & $-\e{0}$  \\ \hline 
$\EE{1}$  & $\EE{1}$ & $0$     & $\eye$     & $-\e{0}$    & $\EE{2}$ & $0$ & $0$ & $0$  \\ \hline
$\EE{2}$  & $\EE{2}$ & $0$     & $\e{0}$   & $\eye$    & $-\EE{1}$ & $0$ & $0$ & $0$  \\ \hline
$\eye$    & $\eye$     & $0$     & $\EE{1}$ & $\EE{2}$ & $-\e{0}$ & $0$ & $0$ & $0$ \\ \hline
\end{tabularx}
\caption{Geometric product in $\pdclal{2}{0}{1}$}
\end{center}
\end{table}

\subsection{Enumeration of various products}
\label{sec:enumpro}

We want to spend a bit of time now investigating the various forms which the geometric product takes in this algebra.  
For this purpose, 
define two arbitrary 1-vectors $\vec{a}$ and $\vec{b}$ and two arbitrary bivectors $\vec{P}$ and $\vec{Q}$ with 
\[
\vec{a} =  a_0 \e{0} + a_1 \e{1} + a_2 \e{2}, ~~~\text{etc.}
\]
These coordinates are of course not instrinsic but they can be useful in understanding how the euclidean metric is working in the various products.  See the companion diagram in \Fig{fig:geomProducts}.

\begin{enumerate}
\item \textbf{Norms.} It is often useful to normalize vectors to have a particular intensity. There are different definitions for each grade:
\begin{itemize}
\item \textbf{1-vectors.} $\vec{a}^2 = \vec{a} \cdot \vec{a} = a_1^2 + a_2^2$.  Define the \emph{norm} of $\vec{a}$ to be $\|\vec{a}\| := \sqrt{\vec{a} \cdot \vec{a}}$.  Then $\dfrac{\vec{a}}{\| \vec{a} \|}$ is a vector with norm 1, defined for all vectors except $\e{0}$ and its multiples.  In particular, all euclidean lines can be normalized to have norm 1.  Note that  when $\vec{a}$ is normalized,  then so is $-\vec{a}$. These two lines represents opposite \emph{orientations} of the line\footnote{Orientation is an interesting topic which lies outside the scope of this article.}.
\item \textbf{2-vectors.} $\vec{P}^2 = \vec{P} \cdot \vec{P} = p_0^2 \EE{0}^2 = -p_0^2$.  Define the \emph{norm} of $\vec{P}$ to be $p_0$ and write it $\|\vec{P}\|$.  Note that this can take positive or negative values, in contrast to $\sqrt{\vec{P} \cdot \vec{P}}$. Then $\dfrac{\vec{P}}{\| \vec{P} \|}$ is a bivector with norm 1, defined for all bivectors except where $p_0 = 0$, that is, ideal points.  In particular, all euclidean points can be normalized to have norm 1. This is also known as \emph{dehomogenizing}.
\item \textbf{3-vectors.} Define  $\mathbf{S}: \pdgrassgr{3}{3} \rightarrow \pdgrassgr{3}{0}$ by $\stripeye{\alpha \eye} = \alpha\one$. This gives the scalar magnitude of a pseudoscalar in relation to the basis pseudoscalar $\eye$.  We sometimes write $\frac{1}{\eye}(\alpha \eye)$ for the same.  In a non-degenerate metric, the same can be achieved by multiplication by $\eye ^{-1}$.
\end{itemize}
\item \textbf{Inverses.} $\vec{a}^{-1} = \dfrac{\vec{a}}{\vec{a} \cdot \vec{a}}$ and $\vec{P}^{-1} = \dfrac{-\vec{P}}{\vec{P} \cdot \vec{P}}$, for euclidean $\vec{a}$ and $\vec{P}$.
\item \textbf{Euclidean distance.} For normalized $\vec{P}$ and $\vec{Q}$, $\| \vec{P} \vee \vec{Q}\| $ is the euclidean distance between $\vec{P}$ and $\vec{Q}$.
\item  \textbf{Free vectors.} For an ideal point $\vec{V}$ (that is, a free vector) and \emph{any} normalized euclidean point $\vec{P}$, $\|\vec{V}\|_\infty := \| \vec{V} \vee \vec{P}\| = \sqrt{v_1^2+v_2^2}$ is the length of $\vec{V}$.  Then $\dfrac{\vec{V}}{\| \vec{V}\|_\infty}$ is normalized to have length 1. 
\item $\vec{a} \wedge \vec{P} = (a_0 p_0 + a_1 p_2 + a_2 p_2)\eye$ vanishes only if $\vec{a}$ and $\vec{P}$ are incident.  Otherwise, when $\vec{a}$ and $\vec{P}$ are normalized, it is equal to the signed distance of the point to the line times the pseudoscalar $\eye$ .
\item $\vec{P} \cdot \vec{a} = (p_2 a_1 -  p_1 a_2)\e{0} + p_0 a_2 \e{1} - p_0 a_1 \e{2}$ is a line which passes through $\vec{P}$ and is perpendicular to $\vec{a}$. Reversing the order changes the orientation of the line.
\item $\vec{a} \wedge \vec{b} =: \vec{T}$ is the intersection point of the lines $\vec{a}$ and $\vec{b}$. For normalized $\vec{a}$ and $\vec{b}$,  $\| \vec{T} \| = \sin{\alpha}$ where $\alpha$ is the angle between the lines. Reversing the order reverses the orientation of the resulting point.
\item $\vec{a} \cdot \vec{b} = \cos{\alpha}$  for normalized vectors $\vec{a}$ and $\vec{b}$. Which of the two possible angles is being measured here depends on the orientation of the lines. 
\item $\vec{P} \vee \vec{Q}$ is the joining line of $\vec{P}$ and $\vec{Q}$..  
\item $\vec{P} \times \vec{Q} =: \vec{T}$ is the ideal point in the direction perpendicular to the direction of the line $\vec{P} \vee \vec{Q}$. 
\item $\vec{a} \eye = a_1\EE{1} + a_2\EE{2}$ is the polar point of the line $\vec{a}$: the ideal point in the perpendicular direction to the line $\vec{a}$.  All lines parallel to $\vec{a}$ have the same polar point.
\item $\vec{P} \eye = p_0\e{0}$ is the polar line of the point $\vec{P}$: for finite points, the ideal line, weighted by the intensity of $\vec{P}$.  Ideal points have no polar line.
\item $\eye^2= 0$.  This is equivalent to the degeneracy of the metric.  Notice that this fact has no effect on the validity of the above calculations.\footnote{In fact, the validity of most of the above calculations \emph{requires} that $\eye^2=0$.}  
\end{enumerate}

\begin{figure}[t]
\sidecaption
{\setlength\fboxsep{0pt}\fbox{\includegraphics[width=.62\columnwidth]{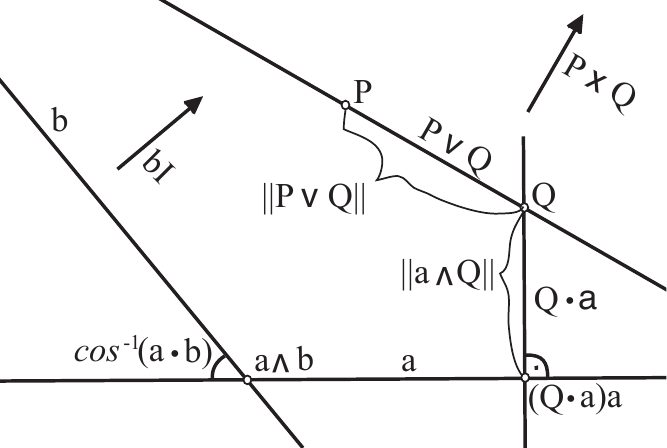}}}
\caption{A selection of the geometric products between various $k$-blades.  Points and lines are assumed to be normalized.  Ideal points are drawn as vectors, distances indicated by norms.}
\label{fig:geomProducts}
\end{figure}

\versions
{ 
\exerc
\label{sec:exercise1}
Rather than presenting finished results based on the above formulas, we present the material in the form of exercises for the reader to work through.  Readers who choose not to work through the results are still recommended to read through them, since what follows will build on the results of the exercises.  Exercises marked with an asterix are more challenging.
\begin{enumerate}
\item \textbf{Angles between vectors.} For normalized ideal points 
$\vec{U}$ and $\vec{V}$, 
\[(\vec{V} \vee \vec{P}) \cdot (\vec{U} \vee \vec{P}) = \cos{\alpha}\]
 where $\alpha$ is the angle between the vectors.  
\item \textbf{Projection onto a line.} For a euclidean line $\vec{a}$ and euclidean point $\vec{P}$, show that $(\vec{P} \cdot \vec{a})\vec{a}^{-1}$ is the orthogonal projection of the point $\vec{P}$ onto $\vec{a}$.  \label{ex:orthproj}
\begin{itemize}
\item What does $(\vec{a} \cdot \vec{P})\vec{P}^{-1}$ represent?  
\item Compare $(\vec{P} \cdot \vec{a})\vec{a}$ and  $(\vec{P} \cdot \vec{a})\vec{a}^{-1}$.
\end{itemize}
\item For normalized arguments, show that $\| \vec{a} \cdot \vec{P}\| = 1$.
\item \textbf{Parallel lines.} What is the condition that $\vec{a}$ and $\vec{b}$ are parallel?  In the expression for the distance of two lines above, what result is obtained when $\vec{a}$ and $\vec{b}$ are parallel?  
\item What changes have to be made to the above formulas when the arguments are not normalized?
\item  \textbf{Ideal elements.} What changes have to be made to the above formulas when the arguments are not finite?  In particular,  are the formulas (6)  and (7) above valid when one or both of $P$ and $Q$ are ideal points?
\item $\vec{P}\times\vec{Q} = (\vec{P} \vee \vec{Q})\eye$.
\item Show that for normalized $\vec{A}, \vec{B}$, and $\vec{C}$, the area of $\Delta \vec{A}\vec{B} \vec{C}$ is given by $(\vec{A} \vee \vec{B})\wedge \vec{C}$.

\begin{figure}[t]
\sidecaption
{\setlength\fboxsep{0pt}\fbox{\includegraphics[width=.62\columnwidth]{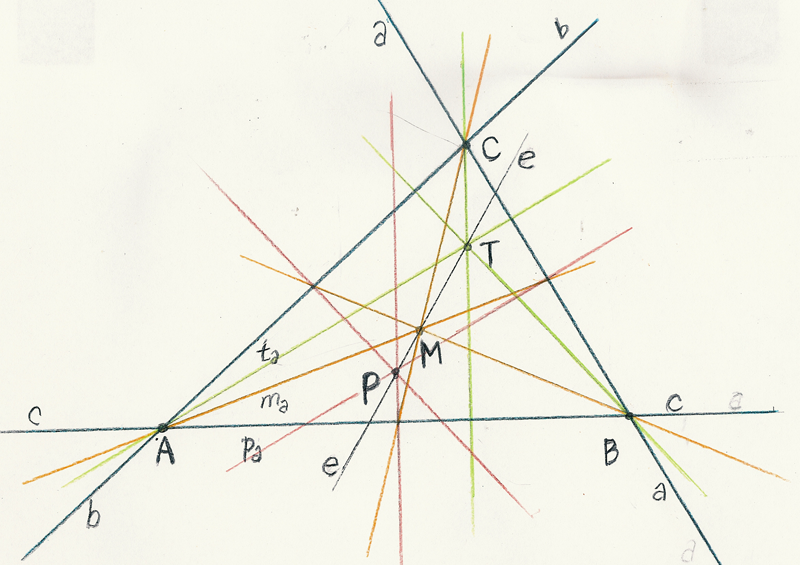}}}
\caption{The three triangle centers $\vec{M}$ (centroid), $\vec{T}$ (orthocenter), and $\vec{P}$ (circumcenter) lying on the Euler line $\vec{e}$ of the triangle $\vec{ABC}$.}
\label{fig:triangleCenters}
\end{figure}
\mymarginpar{Provisional image}
\item \textbf{Triangle centers.} Given a euclidean triangle $\vec{ABC}$, this exercise shows how to use $\pdclal{2}{0}{1}$ to calculate the four classical triangle centers, and prove three lie on the Euler line.  Consult \Fig{fig:triangleCenters}. 
\begin{enumerate} 
\item \label{en:one} Show $\Delta\vec{ABC}$ can be assumed to have normalized corners $A = \EE{0}$, $B=\EE{0} + b_1\EE{1}$, and $C= \EE{0}+c_1 \EE{1} + c_2 \EE{2}$. 
\item Define the edge-lines $\vec{a} := \vec{B}\vee \vec{C}$, $\vec{b} := \vec{C}\vee \vec{A}$, and $\vec{c} := \vec{A}\vee \vec{B}$.
\item Show that the median $\vec{m}_\inert$, perpendicular bisector $\vec{p}_a$, altitude $\vec{t}_a$ and angle bisector $\vec{n}_a$ associated to the pair $\vec{a}$ and $\vec{A}$ are given by the following formulas:
\begin{itemize}
\item $\vec{m}_a = (\vec{B} + \vec{C}) \vee \vec{A}$
\item $\vec{p}_a = (\vec{B} + \vec{C}) \cdot \vec{a}$
\item $\vec{t}_a = \vec{A} \cdot \vec{a}$
\item $\vec{n}_a = \| \vec{c} \| \vec{b}+  \| \vec{b} \|\vec{c}$
\end{itemize}
\item Derive analogous formulas for the lines associated to $\vec{b}$ and $\vec{c}$. 
\item Show that $\vec{m}_a$, $\vec{m}_b$, and $\vec{m}_c$ are co-punctual; their common point is the \emph{centroid} $\vec{M}$ of the triangle. [Hint: to show that the three lines go through the same point you can show the outer product of all three is 0.  If you work with the fully general forms for $\vec{A}$, $\vec{B}$, and $\vec{C}$ (don't use \ref{en:one}), you can also show that the expression for the intersection of two of the lines is symmetric in $\{\vec{a}$,$\vec{b}$,$\vec{c}\}$ and $\{\vec{A}$, $\vec{B}$,$\vec{C}\}$.]
\item Do the same for the triple of lines 
\begin{itemize}
\item ($\vec{p}_a$, $\vec{p}_b$, $\vec{p}_c$), to obtain the \emph{circumcenter} $\vec{P}$,
\item ($\vec{t}_a$, $\vec{t}_b$, $\vec{t}_c$), to obtain the \emph{orthocenter} $\vec{T}$, and
\item ($\vec{n}_a$, $\vec{n}_b$, $\vec{n}_c$), to obtain the \emph{incenter} $\vec{N}$.
\end{itemize}
\item Show that $\vec{M}$, $\vec{P}$, and $\vec{T}$ lie on a line $\vec{e}$ (the \emph{Euler line} of the triangle). [Hint: to show three points lie on a line, show that the join ($\vee$) of two of the points has vanishing outer product with the third.
\item Show that $\vec{M}$ lies between $\vec{P}$ and $\vec{T}$ on the Euler line, and is twice as far from $\vec{T}$ as from $\vec{P}$.
\end{enumerate}
\end{enumerate}

}{For a variety of exercises see \fullversion.}

\subsection{Euclidean isometries via sandwich operations}
\label{sec:eucsand2d}
One of the most powerful aspects of Clifford algebras for metric geometry is the ability to realize isometries as sandwich operations of the form 
\[\vec{X} \rightarrow \vec{g}\vec{X} \vec{g}^{-1}\]
where $\vec{X}$ is any geometric element of the algebra and $\vec{g}$ is a specific geometric element, unique to the isometry.  $\vec{g}$ is in general a \emph{versor}, that is, it can be written as the product of 1-vectors (\cite{hessob87}. 
Let's explore whether this works in $\pdclal{2}{0}{1}$.

\textbf{Reflections.} Let $\vec{a} := \e{0} - \e{1}$ (the line $x=1$), and $\vec{P}$ a normalized point $\EE{0} + x \EE{1} + y \EE{2}$.
Simple geometric reasoning shows  that reflection in the line $\vec{a}$ sends the point $(x,y)$ to the point $(2-x,y)$. 
\ifthenelse{\equal{\isfullversion}{false}}{Alternatively, the reader is encouraged to verify the missing steps of the following computation:
\begin{eqnarray*}
\vec{P}' &:=& \vec{a} \vec{P} \vec{a}^{-1} = \vec{a} \vec{P} \vec{a} = \dots \\
&=& (x-2)\EE{1} -\EE{0} - y \EE{2} \\
&=& \EE{0} + (2-x)\EE{1} + y\EE{2}
\end{eqnarray*}}
{Let's evaluate the versor operator:
\begin{eqnarray*}
\vec{P}' &:=& \vec{a} \vec{P} \vec{a}^{-1} = \vec{a} \vec{P} \vec{a} = \dots \\
&=& (\e{0} - \e{1})(\EE{0} + x \EE{1} + y \EE{2}) (\e{0} - \e{1}) \\
&=& (\e{0} - \e{1})(\e{1}\e{2} + x \e{2}\e{0}+ y \e{0}\e{1}) (\e{0} - \e{1}) \\
&=& (\e{0} - \e{1})(\e{1}\e{2}\e{0}-\e{1}\e{2}\e{1}- x \e{2}\e{0}\e{1} - y \e{0}\e{1}\e{1}) \\
&=& (\e{0} - \e{1})((1-x)\e{0}\e{1}\e{2}+\e{2} - y \e{0}) \\
&=& (\e{0}\e{2}  - (1-x)\e{1}\e{0}\e{1}\e{2}-\e{1}\e{2} + y \e{1}\e{0}) \\
&=& (x-2)\E{1} -\E{0} - y \E{2} \\
&=& \E{0} + (2-x)\E{1} + y\E{2}
\end{eqnarray*}}
where the final step is obtained by dehomogenizing.\ifthenelse{\equal{\isfullversion}{false}}{}{\footnote{Without homogenizing, the orientation of the point is reversed, as it probably should be.} }

This algebra element corresponds to the euclidean point $(2-x,y)$, so the sandwich operation \textbf{is} the desired reflection in the line $\vec{a}$.
We leave it as an exercise for the interested reader to carry out the same calculation for a general line.

 \textbf{Direct Isometries.}
By well-known results in plane geometry,
\versions{ \begin{quote}{The result of carrying out  reflections in two lines one after the other (the \emph{composition} of the reflections) is a rotation around the intersection point of the two lines, through an angle equal to twice the angle between the two lines, unless the two lines are parallel, in which case the composition is a translation in the direction perpendicular to the two lines, through a distance equal to twice the distance between the two lines.}\end{quote}}
{the composition of two reflections yields a rotation around the common point of the two lines.} Translating this into the language of the Clifford algebra,  the composition of reflections in lines $\vec{a}$ and $\vec{b}$ will look like:
\begin{eqnarray*} \label{eqn:tworefl}
\vec{P}' &=& \vec{b} (\vec{a} \vec{P} \vec{a}) \vec{b} \label{eqn:rotator}\\
&=& \vec{T} \vec{P} \vec{\widetilde{T}} 
\end{eqnarray*}
where we write $\vec{T} := \vec{b}\vec{a}$, and $\widetilde{\vec{T}}$ is the reversal of $\vec{T}$.

Note that the intersection of the two lines will be fixed by the resulting isometry.  There are two cases: the point is ideal, or it is euclidean. In the case of an ideal point, the two lines are parallel and the composition is a translation.  Let's look at an example.

Retaining $\vec{a}$ as above, define the normalized line $\vec{b} := 2 \e{0} - \e{1}$, the line $x=2$.  
By simple geometric reasoning, the composition "reflect first in $\vec{a}$, then in $\vec{b}$" \emph{should} be the translation $(x,y) \rightarrow (x+2,y)$.  Defining $\vec{T} := \vec{b} \vec{a}$,  the sandwich operator looks like: $\vec{T} \vec{P} \vec{\tilde{T}}$. 
Calculate the product $\vec{T} = \one - \EE{2}$, and 
 $\vec{T} \vec{P} \vec{\tilde{T}} = \EE{0} + (2+x) \EE{1} + \EE{2}$. 
This shows that $\vec{T} \vec{P} \vec{\widetilde{T}} $ is the desired translation operator.
One can generalize the above to show that a translation by the vector $(x_0, y_0)$  is given by the sandwich operation $  \vec{T}\vec{P}\vec{\widetilde{T}}$ where
\begin{equation} \label{eqn:translator}
\vec{T} := \one + \frac12 (y_0 \EE{1} -x_0 \EE{2})
\end{equation}
It's interesting to note that $\vec{T} \vec{P}$ and $\vec{P} \vec{\tilde{T}}$ are both translations of $(x,y) \rightarrow (x+1,y)$, so one doesn't need a sandwich to implement translations, but for simplicity of representation we continue to do so.

\textbf{Rotations.}
Similar remarks apply to rotations. A rotation around a normalized point $\vec{R}$ by an angle $\theta$ is given by
\versions{\[ \vec{T} = \cos{(\frac{\theta}{2})} + \sin{(\frac{\theta}{2})}\vec{R}\]}
{ $\vec{T} = \cos{(\frac{\theta}{2})} + \sin{(\frac{\theta}{2})}\vec{R}$.}
This can be checked by substituting into  \Eq{eqn:rotator} and multiplying out.  
We'll explore a method for constructing such rotators using the exponential function, in the next section.
See \fullversion for a more detailed discussion including constructions of glide reflections and point reflections. 

\versions{

%
%
\subsubsection{Reflections in points and in hyperplanes}
\label{sec:buildrefl}
When $\vec{a}$ is a 1-vector, then $\vec{a} \vec{X} \vec{a}$ is the reflection in $\vec{a}$.  In the dual algebra, $\vec{a}$ represents a hyperplane (in this case, a line); in the standard algebra, a point.  Because reflections in planes are often more practically useful than reflections in points, this can be seen as an advantage for the dual approach used here.  
  In  a non-degenerate metric this is not a significant advantage,  because in such metrics a reflection in a hyperplane is \emph{also} a reflection in the polar point of the plane; and vice-versa. 

This is also good place to point out that a common form for a reflection in a hyperplane, in the vector space model, includes a minus sign (\cite{dfm07}, p. 168).  Thus, the reflection in the plane whose normal vector is $\vec{P}$ is given by:
\[ \vec{X} \rightarrow -\vec{P} \vec{X} \vec{P}\]
This minus sign is due to the fact that the algebra is point-based but the desired reflection is in a hyperplane.  Without minus sign, the expression $\vec{P} \vec{X} \vec{P}$ represents a \emph{reflection} in the \emph{vector} $\vec{P}$.  In 3 dimensions, this is a rotation of 180 degrees around the  line $\lambda \vec{P}$.  To obtain the reflection in the plane $\vec{p}$ orthogonal to $\vec{P}$, one must compose the vector reflection with the point reflection in the origin, which is achieved by multiplying by $-1$.  This undoes the rotation around the line $\lambda\vec{P}$ and introduces a reflection in the orthogonal plane.  This yields the expression $-\vec{P} \vec{X} \vec{P}$ as the form for a reflection in the plane orthogonal to $\vec{P}$.  Compare this to the homogeneous model presented here, where reflections in points and in planes are represented without any extra minus signs. 

\exerc
\label{sec:exercise2}

\begin{enumerate}
\item A glide reflection is the product of a reflection in a line with a translation parallel to the line.  Describe how to realize a glide reflection using a sandwich operation.
\item A point reflection in a point $\vec{P}$ is an isometry that sends each point $\vec{Q}$ to its reflected image on the \quot{other side} of $\vec{P}$.  In particular, the image of $\vec{Q}$ lies on the line $\vec{P} \vee \vec{Q}$ on the other side of $\vec{P}$ an equal distance to $\vec{P}$.  Show that $\vec{Q} \rightarrow   \vec{P}  \vec{Q} \widetilde{\vec{P}}$ realizes this point reflection.  How is this consistent with the theory developed above that such rotors are rotations? [Hint: $\cos{\frac\pi{2}} = 0$.]
\item\hard Show that the above formulas for the action of reflections, translations, and rotations on \emph{points} are also valid when $\vec{P}$ is replaced by a euclidean \emph{line} $\vec{a}$.
\end{enumerate}

}{}

\subsection{Spin group, Exponentials and Logarithms}
\label{sec:spinexplog}

We have seen above in \Eq{eqn:rotator} that euclidean rotations and translations can be represented by sandwich operations in $\pdclal{2}{0}{1}$,  in fact, in the even subalgebra $\pdclplus{2}{0}{1}$.

\begin{definition}
The \emph{spin group} $\spin{2}{0}{1}$ consists of elements 
$\vec{g}$  of the even subalgebra $ \pdclplus{2}{0}{1}$ such that $\vec{g} \vec{\widetilde{g}} = 1$.   An element of the spin group is called a \emph{rotor}. \versions{Some authors (\cite{perwass09}) refer to an element of the spin group as a \emph{spinor} but this conflicts with the accepted definition of spinor in the mathematics community, so we avoid using the term here.}{}
\end{definition}
Write $\vec{g} = s \one + \vec{M}$ where $\vec{M} = \langle \vec{g} \rangle_2 = m_0 \EE{0} + m_1 \EE{1} + m_2 \EE{2}$. Then $\vec{g} \vec{\widetilde{g}} = s^2 + m_0^2 = 1$.  
There are two cases. 
\begin{itemize}
\item $m_0 \neq 0$, so $\vec{M}$ is a euclidean point.  Then there exists $\theta \neq 0$ such that $s = \cos{(\theta)}$ and $m_0 = \sin{(\theta)}$, yielding  $\vec{g} = \cos{(\theta)} + \sin{(\theta)}\vec{N}$ where $\vec{N} = \dfrac{\vec{M}}{\sin{\theta}}$ is a normalized point, hence $\vec{N}^2 = -1$.  Thus, the formal exponential $e^{t\vec{N}}$ can be evaluated to yield:
\begin{eqnarray*}
e^{t\vec{N}} &=& \sum_{i=0}^\infty \frac{(t\vec{N})^i}{i!} \\
&=& \cos{(t)} + \sin{(t)}\vec{N}
\end{eqnarray*}
Hence, the rotor $\vec{g}$ can be written as an exponential: $\vec{g} = e^{\theta \vec{N}}$.
\item $m_0 = 0$, so $\vec{M}$ is an ideal point. Then we can assume that $s = 1$ (if $s = -1$,  take the element $-\vec{g}$ with the same sandwich behavior as $\vec{g}$). Also, $\vec{M} = m_1 \EE{1} + m_2 \EE{2}$. Again, the formal exponential $e^{t\vec{M}}$ can be evaluated to yield:
\begin{eqnarray*}
e^{t\vec{M}} &=& \sum_{i=0}^\infty \frac{(t\vec{M})^i}{i!} \\
&=& 1 + t\vec{M}
\end{eqnarray*}
So in this case, too, $\vec{g} = e^{\vec{M}}$ has an exponential form. 
\end{itemize}
The above motivates the following definitions:
\begin{definition}
A \emph{rotator} is a rotor whose  bivector part is a euclidean point. A \emph{translator} is a rotor whose bivector part is a ideal point.  
\end{definition}
\begin{definition}
The \emph{logarithm} of a translator $\vec{g} = s\one + \vec{M} \in \pdclal{2}{0}{1}$ is $\vec{M}$, since $e^{\vec{M}}= \vec{g}$.
\end{definition}
\begin{definition}
Given a rotator $\vec{g} = s\one + m_0 \EE{0} + m_1 \EE{1} + m_2 \EE{2} \in \spin{2}{0}{1}$. Define $\theta :=  \text{tan}^{-1}{(m_0, s)}$ and $\vec{N} := \dfrac{\vec{M}}{\|\vec{M}\|}$. Then the \emph{logarithm} of $\vec{g}$ is $\theta \vec{N}$, since $e^{\theta \vec{N}}= \vec{g}$.
\end{definition}

\versions{
\exerc
\label{ex:exercise3}
\begin{enumerate}
\item Show that a translation of the point $\vec{P}$ by the vector $(x_0,y_0)$ can be written as 
$
e^{\vec{M}}\vec{P}e^{-\vec{M}}
$
for $\vec{M} = \one + \frac12(y_0 \EE{1} - x_0 \EE{2})$.
\item Show that a rotation of the point $\vec{P}$ around the point $(x_0,y_0)$ by angle $\theta$ can be written as 
$
e^{\vec{M}}\vec{P}e^{-\vec{M}}
$
for $\vec{M} = \sin{(\frac\theta2)}(\EE{0} + x_0\EE{1} + y_0\EE{2})$.
\item Deduce a condition on two rotors $g_1$ and $g_2$ so that their product is a translator.
\end{enumerate}
}{}
\mysubsubsection{Lie groups and Lie algebras}
The above remarks provide a realization of the 2-dimensional euclidean direct isometry group $SE(2)$ and its Lie algebra $se(2)$ within $\pdclplus{2}{0}{1}$.  The Spin group $
\spin{2}{0}{1}$ forms a double cover of $SE(2)$ since the rotors $g$ and $-g$ represent the same isometry.   Within $\pdclplus{2}{0}{1}$, the spin group consists of elements of unit norm; the Lie algebra consists of the pure bivectors plus the zero element.  
The exponential map $\vec{X} \rightarrow e^{\vec{X}}$ maps the latter bijectively onto the former.  This structure is completely analogous to the way the unit quaternions sit inside $\pdclplus{3}{0}{0}$ and form a double cover of $SO(3)$.  The full group including indirect isometries is also naturally represented in $\pdclal{2}{0}{1}$ as the group generated by reflections in lines,  sometimes called the \emph{Pin} group. 

\subsection{Guide to the literature}

There is a substantial literature on the four-dimensional even subalgebra $\pdclplus{2}{0}{1}$ with basis $\{\one, \EE{0}, \EE{1}, \EE{2}\}$.  In an ungraded setting, this structure is known as the \emph{planar quaternions}.  The original work appears to have been done by Study (\cite{study91}, \cite{study03}); this was subsequently expanded and refined by Blaschke (\cite{blaschke38}). Study's parametrization of the full planar euclidean group as \quot{quasi-elliptic} space is worthy of more attention.  Modern accounts include \cite{mcc}.
\versions{
\footnote{which however confuses the euclidean inner product on vectors with the inner product on points.}  }{}

\section{$\pdclal{3}{0}{1}$ and euclidean space}
\label{sec:cl301}
The extension of the results in the previous section to the three-dimensional case $\pdclal{3}{0}{1}$ is mostly straightforward.  Many of the results can be carried over virtually unchanged.   The main challenge is due to the existence of non-simple bivectors; in fact, most bivectors are \textbf{not} simple! (See \Sec{sec:propbiv} below.)
This means that the geometric interpretation of a bivector is usually \textbf{not} a simple geometric entity, such as a spear or an axis, but a more general object known in the classical literature as a \emph{linear line complex}, or  \emph{null system}.  Such entities are crucial in kinematics and dynamics; we'll discuss them below in more detail.

\textbf{Notation.} As a basis for the full algebra we adopt the terminology for the exterior algebra $W^*$ in \Sec{sec:wwstar}, interpreted as a \emph{plane-based} algebra.
We add an additional basis 1-vector satisfying $\e{3}^2 = 1$.  $\e{0}$ now represents the ideal \emph{plane} of space, the other basis vectors represent the coordinate planes. 
$\EE{0}$ is the origin of space while $\EE{1} = \e{0}\e{3}\e{2}$ is the ideal point in the $x$-direction, similarly for $\EE{2}$ and $\EE{3}$.   The bivector $\e{01}$ is the ideal line in the $x=0$ plane, and similarly for $\e{02}$ and $\e{03}$.  $\e{23}, \e{31}$, and $ \e{12}$ are the $x$-, $y$-, and $z$-axis, resp.   We use $\mathbf{i}$ again to denote the embedding of euclidean points, lines, and planes, from $\RP{3}$ into the Clifford algebra.  

We continue to denote 1-vectors with bold small Roman letters $\vec{a}$; trivectors will be denoted  with bold capital Roman letters $\vec{P}$; and bivectors will be represented with bold capital Greek letters $\pip$.\footnote{A convention apparently introduced by Klein, see \cite{klein72a}.}

We leave the construction of a multiplication table as an exercise.  Once again, most of the geometric products of two vectors obey the pattern $AB = A\cdot B + A \wedge B$. Two new exceptions involve the product of a bivector with another bivector, and with a trivector:
\begin{eqnarray} \label{eqn:gpbivector}
\pip\fio &=& \pip \cdot \fio + \pip \times \fio + \pip \wedge \fio \\
\pip \vec{P} &=& \pip \cdot \vec{P} + \pip \times \vec{P}
\end{eqnarray}

Here, as before, the commutator product $A \times B := \frac12(AB - BA)$.  

We now describe in more detail the nature of bivectors. We work in $W^*$, since that is the foundation of the metric.  As a result, even readers familiar with bivectors from a \emph{point-based} perspective will probably benefit from going through the following \emph{plane-based} development.

\subsection{Properties of Bivectors.}  
\label{sec:propbiv}

We begin with a simple bivector $\pip := \vec{a} \wedge \vec{b}$ where $\vec{a}$ and $\vec{b}$ are two planes with coefficients $\{a_i\}$ and $\{b_i\}$.  The resulting bivector has coefficients 
\begin{equation}\label{eqn:pluecker}
p_{ij}:=a_ib_j-a_j b_i ~~~(ij \in \{01, 02, 03, 12, 31, 23\})
\end{equation} These are the plane-based Pl\"{u}cker coordinates for the intersection line (\emph{axis}) of $\vec{a}$ and $\vec{b}$.  Clearly $\pip \wedge \pip = 0$.  Conversely, if \[\pip \wedge \pip = 2 (p_{01}p_{23} + p_{02} p_{31} + p_{03} p_{12})\eye = 0\] for a bivector $\pip$, the bivector is simple \cite{nh:pg}.

Given a second axis $\fio = \vec{c} \wedge \vec{d}$, the condition $\pip \wedge \fio = 0$ implies they have a plane in common, or, equivalently, they have a point in common. 
For general bivectors, 
\begin{equation}
\label{eqn:pluecker}
\pip \wedge \fio= (p_{01}q_{23} + p_{02} q_{31} + p_{03} q_{12} + p_{12}q_{03} + p_{31} q_{02} + p_{23} q_{01})\eye
\end{equation}
The parenthesized expression is called the Pl\"{u}cker inner product of the two lines, and is written $\langle \pip, \fio \rangle_P$. With this inner product, the space of bivectors $ P(\bigwedge^2(\R{4})^*)$ is the Cayley-Klein space $\complexspace :=\proj{\mathbb{R}^{3,3}}$, and the space of lines  is the quadric surface $\linespace = \{ \pip \mid \langle \pip, \pip \rangle_P = 0\} \subset \complexspace$.  \versions{$\complexspace$ is sometimes called \emph{Klein} quadric.}{}
When \Eq{eqn:pluecker} vanishes, the two bivectors are said to be \emph{in involution}.  

\versions{
\mysubsubsection{Pencils of line complexes} \label{eqn:linecomppencil}
Two bivectors $\pip$ and $\fio$ span a line in \complexspace, called a line complex pencil, or bivector pencil.  Points on this line are of the form  $\alpha \pip  + \beta \fio$ for $\alpha, \beta \in \mathbb{R}$, both not 0. Finding the simple bivectors on this line involves solving the following equation in the homogeneous coordinate $\lambda = \alpha : \beta$:
\begin{eqnarray*}
0 &=& \langle \alpha \pip  + \beta \fio, \alpha \pip  + \beta \fio\rangle_P \\
&=& \alpha^2 \langle \pip, \pip \rangle_P + 2\alpha \beta \langle \fio, \pip  \rangle_P + \beta^2\langle \fio,  \fio \rangle_P
\end{eqnarray*}
This equation can be undetermined, or have 0, 1, or 2 real homogeneous roots.  See \cite{pottmann01} for details.  Finding the intersection of a line in \complexspace with \linespace is a common procedure in line geometry.  It's used below to calculate 
the axis of a non-simple bivector. 
}{}

\versions{
\mysubsubsection{Null polarity associated to a non-simple bivector} For a fixed $\pip$, the orthogonal complement $\pip^\perp :=  \{ \fio \mid \langle \pip, \fio \rangle_P = 0\}$ is a 4-dimensional hyperplane of $\complexspace$ consisting of all bivectors in involution to $\pip$. The intersection with $\linespace$ is a 3-dimensional quadratic submanifold $\mathbf{M}^3_2$ called a \emph{linear line complex}.  When $\pip$ is simple, this consists of all lines which intersect $\pip$.  

Assume that $\pip$ is not simple.  Then $\pip$ determines a collineation, the \emph{harmonic homology} $\mathbf{H}_\pip: \RP{5} \rightarrow \RP{5}$, with center $\pip$ and axis $\pip^\perp$ given by the following formula:
\begin{align}
\mathbf{H}(\fio) &:= \langle \pip, \pip^\perp \rangle\fio - 2 \langle \fio, \pip^\perp\rangle \pip \\
&:= \langle \pip, \pip \rangle_P\fio - 2 \langle \fio, \pip\rangle_P \pip
\end{align}
Here $\langle , \rangle$ is the canonical pairing of a vector and a dual vector in the vector space $\R{6}$ underlying $\RP{5}$.  When expressed in terms of the Pl\"{u}cker inner product, this takes the second form given above.

$\mathbf{H}_\pip$ is a collineation of $\complexspace$ which fixes  $\pip$ and $\pip^\perp$ pointwise.  It's a kind of mirror operation when applied to a point $\fio$: Let $\vec{X}$ be the line in $\complexspace$ joining $\pip$ and $\fio$, and $\pip_\perp$ the intersection of $\vec{X}$ with $\pip^\perp$.  Then $\mathbf{H}_\pip(\fio)$ is the harmonic partner of $\fio$ with respect to the fixed points $\pip$ and $\pip_\perp$, the point $\thio$ on $\vec{X}$  satisfying the cross ratio condition $\mathbf{CR}(\thio, \fio; \pip, \pip_\perp) = -1$.  If $\fio = \alpha \pip + \beta \pip_\perp$, then $\thio = \alpha \pip - \beta \pip_\perp$ (Exercise).  When $\fio$ is  simple, then so is $\thio$ (Exercise).  To proceed, we need a theorem relating projectivities of $\RP{3}$ and projectivities of $\complexspace$ preserving $\linespace$.  See \cite{pottmann01}, p. 144, for a proof.

\begin{theorem}
The group of projectivities of $\RP{3}$ and the group of collineations of $\complexspace$ preserving $\linespace$ are isomorphic.
\end{theorem}

This result confirms that the harmonic homology $\mathbf{H}_\pip$ has an induced action on $\RP{3}$, the \emph{null polarity} associated to $\pip$.  An simple element of $\pip^\perp$, considered as a line in $\RP{3}$, is called a \emph{null line} of  the null polarity. First of all, how do the null lines of $\pip$ appear?  Consider a point $\vec{P}$, and define the null plane of $\vec{P}$ to be $N_\pip(\vec{P}) := \pip \vee \vec{P}$.   Then $N_\pip(\vec{P}) \vee \vec{P} = 0$, so $\vec{P}$ lies in its null plane.  Also, let $\vec{Q}$ be another point of $N_\pip(\vec{P})$. Then apply the associativity of the $\vee$ product to obtain:
\begin{align*}
0 = (N_\pip(\vec{P}) \vee \vec{Q} &= ( \pip \vee \vec{P}) \vee \vec{Q} \\
&= \pip \vee (\vec{P} \vee \vec{Q})
\end{align*}
This shows that the line $\vec{P} \vee \vec{Q}$ is a null line of $\pip$.  Conversely, any null line of $\pip$ passing through $\vec{P}$ can be written in the form $\vec{P} \vee \vec{Q}$, and one can reverse the reasoning to conclude it must lie in $N_\pip(\vec{P})$.  

One can also dualize the discussion to define the \emph{null point} of a plane $N_\pip(\vec{a}) := \pip \wedge \vec{a}$. $N_\pip(N_\pip(\vec{a})) = \vec{a}$ by \Eq{eqn:ex4-4} , showing that $N_\pip$ is an involution, hence in fact deserves the name null \emph{polarity}.  

\mysubsubsection{Comparison to versors}  \Sec{sec:spinexplog3d} shows that euclidean isometries can be represented in 3D also with versor operators of the form $\vec{g}\vec{X}\widetilde{\vec{g}}$.  These isometries are projective transformations of $\RP{3}$ that preserve the metric quadric.  The null polarity $N_\pip$ is also a projective transformation of $\RP{3}$ but it is realized within the Clifford algebra in a different way\footnote{Actually it's more accurate to say it is realized in the Grassman algebra since it doesn't involve the inner product}.  The action on planes, lines, and points is given by:
\begin{align}
N_\pip(\vec{a}) &= \pip \wedge \vec{a} \\
N_\pip(\fio) &= \stripeye{\pip \wedge \pip} \fio - 2 \stripeye{\fio \wedge \pip}\pip \\
N_\pip(\vec{P}) &= \pip \vee \vec{P}
\end{align}

Of course, once the action on points, on lines, or on planes is given, the action on the other two types of elements is determined. One can construct other projective transformations by composing such null polarities.  For example, if $\pip$ and $\fio$ are in involution, then the composition $N_\fio \circ N_\pip = N_\pip \circ N_\fio$ is an involution which preserves each null line common to both line complexes $\fio$ and $\pip$.  By composing three such null polarities, one arrives at the polarity with respect to a regulus (the common null lines of all 3 line complexes). For details see \cite{weiss35}.   

By taking all possible compositions of these null polarities, one arrives at a representation of the full projective group of $\RP{3}$ within the Grassman algebra. These representations do not have the simplicity and homogeneity of the versor representation of isometries. Still, this representation of the projective group seems worth investigating.  For example, it should be possible to implement uniform scaling around a given euclidean point in this way. 

With these remarks we close our discussion of the projective properties of bivectors and the induced null polarities. We'll meet the null system again in Section  \ref{sec:rbm} since it  is fundamental to understanding rigid body mechanics.
}{

\mysubsubsection{Null system} A line which is in involution with a given bivector $\pip$ is called a \emph{null} line of $\pip$.  Through every point and in every plane of space, lie a line pencil of null lines.  In the case of a non-simple $\pip$, this sets up an polarity\footnote{A \emph{polarity} is an involutive projectivity that swaps points and planes.} between the points and planes of space called the \emph{null polarity} determined by $\pip$:  the null plane of a point is the plane in which the null lines of the point lie, and vice-versa.  \Sec{sec:enumpro3} shows how the null plane and null point can be expressed in the Clifford algebra.  
We'll meet the null system again in Section  \ref{sec:rbm} since it  is fundamental to understanding rigid body mechanics.
}

\mysubsubsection{Metric properties of bivectors}
Write the bivector $\pip$ as the sum of two simple bivectors $\pip = \pip_\infty + \pip_o$:
\begin{eqnarray*}
\pip_\infty := p_{01} \e{01} + p_{02} \e{02}+p_{03} \e{03} \\
\pip_o := p_{12} \e{12} + p_{31} \e{31}+p_{23} \e{23}
\end{eqnarray*}
This is the unique decomposition of $\pip$ as the sum of a line lying in the ideal plane ($\vec{\pip}_\infty$) and a euclidean part ($\vec{\pip}_o$). We sometimes write $\pip = (\pip_{\infty}; \pip_o)$. 
\versions{
This decomposition is useful in characterizing the bivector. 
$\pip_o \wedge \pip_\infty = 0 \iff \pip_o \vee \pip_\infty = 0 \iff \pip$ is simple.
 We say a bivector is \emph{ideal} if $\pip_o = 0$, otherwise it is \emph{euclidean}.  
$\pip_o$ is a line through the origin, whose direction is given by the ideal point  $\vec{N}_\pip := \e{0} \pip_o =  p_{23} \EE{1}+p_{31} \EE{2} + p_{12} \EE{3}$.  We call $\vec{N}_\pip $ the \emph{direction vector} of the bivector.  The following identities are left as exercises for the reader:
\begin{align}
\vec{N}_\pip &= \e{0} \pip \\ \label{eqn:direction}
\EE{0} \vee \e{0}\pip_o &= \pip_o
\end{align}

$\pip_o$ is invariant under euclidean translations, $\pip_\infty$ is not.  This deserves a closer look. First, consider the case of a line passing through the origin.  Then $\pip = \EE{0} \vee \vec{N}$ where $\vec{N}$ is an ideal point.  Note that $\pip_\infty = 0$ in this case.  Let $\vec{T}$ be the ideal point representing a translation vector.  Then the image of $\pip$ under this translation is:
\begin{align}
\pip_T &=  (\EE{0}+\vec{T}) \vee \vec{N} \\
&= \pip + \vec{T} \vee \vec{N} \\ \label{eqn:transbiv}
&= (\vec{T} \vee \vec{N} ~;~\pip_o)
\end{align}
where we have used the fact that ideal points are invariant under translations, and that the join of two ideal points is an ideal line ($ \vec{T} \vee \vec{N} $).  This leads us to the following proposition:
\begin{theorem}\label{thm:transline}
A bivector of the form $(\vec{T} \vee \e{0} \pip_o; \pip_o)$ where $\vec{T}$ is ideal, is a simple bivector passing through the point $\EE{0} + \vec{T}$, and every such bivector can be so represented.
\end{theorem}
\begin{proof}
Set $\vec{N} := \e{0}\pip_o$ in \Eq{eqn:transbiv}.  Then $\pip_\infty =\vec{T}\vee  \e{0} \pip_o  $. The resulting bivector $\pip_\infty + \pip_o$ is simple since $\pip_\infty \vee \pip_o = 0$ and passes through $\EE{0} + \vec{T}$:
\begin{align}
\pip \vee (\EE{0} + \vec{T}) &=( \vec{T}\vee \e{0} \pip_o  + \pip_o) \vee (\EE{0} + \vec{T}) \\
&= (\vec{T} \vee \e{0} \pip_o )\vee\EE{0} + \pip_o \vee \vec{T} \\
&=\vec{T}  \vee (\EE{0} \vee \e{0}\pip_o) - \vec{T} \vee  \pip_o  \\
&=\vec{T} \vee (\EE{0} \vee \e{0}\pip_o - \pip_o) \\
&= 0
\end{align}
Here we have used the fact that $ (\vec{T}\vee \e{0} \pip_o)\vee \vec{T} = 0$, $\pip_o \vee \EE{0} = 0$, anti-symmetry and distributivity of $\vee$, and, finally, \Eq{eqn:direction}.
\end{proof}
When $\pip$ is not simple, we have a related result: 
\begin{theorem}\label{thm:transbiv}
A bivector of the form $(\pip_\infty+ \vec{T} \vee \e{0} \pip_o; \pip_o)$ where $\vec{T}$ is ideal, is the image of the bivector $\pip = \pip_\infty + \pip_o$ under the translation by $\vec{T}$, and every such bivector can be so represented.
\end{theorem}
\begin{proof}
We know from \Thm{thm:transline} that the result is true when $\pip_\infty = 0$.  But since $\pip_\infty$ is an ideal line, it is fixed by a euclidean translation.   The result follows immediately. \qed
\end{proof}
\Thm{thm:transline} can be so intepreted:  the set of all lines sharing the same finite part $\pip_o$ is the line bundle centered at the ideal point $\vec{N}_\pip$.  This bundle is mapped to itself by a translation $\vec{T}$, in such a way that  two elements in the bundle differ by the ideal line bundle element $\vec{T} \vee \e{0} \pip_o$  when one is the translated image of the other under $\vec{T}$.  This implies that a change of coordinate which moves $\EE{0}$ to $\EE{0} + \vec{T}$ will produce such a result on the coordinates of bivectors.  We'll return to this point below  in the discussion of the axis of a bivector (\Sec{sec:dualno}), where we show how to choose a canonical representative from this line bundle to represent a general bivector $\pip$.
}{
$\pip_o \wedge \pip_\infty = 0 \iff \pip$ is simple.
$\pip_o$ is invariant under euclidean translations, $\pip_\infty$ is not (Exercise).
We say a bivector is \emph{ideal} if $\pip_o = 0$, otherwise it is \emph{euclidean}.  
$\pip_o$ is a line through the origin, whose direction is given by the ideal point  $\vec{Q} := e_0 \pip_o = p_{23} \EE{1}+p_{31} \EE{2} + p_{12} \EE{3}$.  We call $\vec{Q}$ the \emph{direction vector} of the bivector. 
}.

\subsection{Enumeration of various products}
\label{sec:enumpro3}

All the products described in $\ref{sec:enumpro}$ have counter-parts here, obtained by leaving points alone and replacing lines  by planes.  We leave it as an exercise to the reader to enumerate them.
Here we focus on the task of enumerating the products that involve bivectors.
For that purpose, we extend the definition of $\vec{a}$,$\vec{b}$, $\vec{P}$, and $\vec{Q}$ to have an extra coordinate, and introduce two arbitrary bivectors, which may or may not be simple: 
$\pip := p_{01} \e{01} + ... $ and 
$ \fio :=  g_{01} \e{01} + ... $. 

\begin{enumerate}
\item \textbf{Inner product.} $\pip \cdot \fio = -(p_{12} g_{12} +  p_{31} g_{13} + p_{23} g_{23})  = -\cos{(\alpha)}$ where $\alpha$ is the angle between the direction vectors of the two bivectors (see \Sec{sec:propbiv} above). $\pip \cdot \fio$ is a symmetric bilinear form on bivectors, called the \emph{Killing} form. We sometimes write $\pip \cdot \fio = \langle \pip, \fio \rangle_k$.  Note that just as in the 2D case, the ideal elements play no role in this inner product. This angle formula is only valid for euclidean bivectors.
\item \textbf{Norm.} There are two cases:
\begin{enumerate}
\item \textbf{Euclidean bivectors.}  For euclidean $\pip$, define the norm $\| \pip \| = \sqrt{-\pip \cdot \pip}$.  Then $\dfrac{\pip}{\| \pip \|}$ has norm 1; we call it a normalized euclidean bivector. 
  \item \textbf{Ideal bivectors.} As in \Sec{sec:enumpro}, we get the desired norm on an ideal line $\pip$ by joining the line with \emph{any} euclidean point $\vec{P}$, and taking the norm of the plane: $\| \pip \|_\infty = \| \pip \vee \vec{P}\|$.  We normalize ideal bivectors with respect to this norm.
\end{enumerate}
\item \textbf{Distance.}  Verify that the euclidean distance of two normalized points $\vec{P}$ and $\vec{Q}$ is still given by $\| \vec{P} \vee \vec{Q} \|$, and the norm of an ideal point $\vec{V}$ (i. e., vector length) is given by $\| \vec{V} \vee \vec{P} \|$ where $\vec{P}$ is any normalized euclidean point.
\item \textbf{Inverses.} For euclidean $\pip$, define $\pip^{-1} = \dfrac{\pip}{\pip \cdot \pip}$.  Inverses are unique.
\item $\pip \wedge \fio = \pip \vee \fio =  \langle \pip, \fio \rangle_P \eye$ is the Pl\"{u}cker inner product  times $\eye$.  
When both bivectors are simple, this is proportional to the euclidean distance between the two lines they represent (Exercise).  
\item \textbf{Commutator.} $\pip \times \fio$ is a bivector which is in involution to both $\pip$ and $\fio$ (Exercise). We'll meet this later in the discussion of mechanics (\Sec{sec:rbm}) as the Lie bracket.  
\item  \textbf{Null point.} $\vec{a} \wedge \pip$, for simple $\pip$,  is the intersection point of $\pip$ with the plane $\vec{a}$; in general it's the \emph{null point }of the plane with respect  $\pip$.
\item \textbf{Null plane.} $\vec{P} \vee \pip$, for simple $\pip$, is the joining plane of $\vec{P}$ and $\pip$; in general it's the \emph{null plane} of the point with respect to  $\pip$.
\versions{
\item $\vec{P} \cdot \pip$, for simple $\pip$, is a plane passing through $\vec{P}$ perpendicular to $\pip$.
\item $\vec{P} \times \pip$, for simple $\pip$, is the normal direction to the plane through $\vec{P}$ and $\pip$.}{}
\item $\vec{a} \cdot \pip$, for simple $\pip$, is a plane containing $\pip$ whose intersection with $\vec{a}$ is perpendicular to $\pip$.
\item $\pip \eye = (p_{23} \e{01} + p_{31} \e{02} + p_{12} \e{03})\eye $ is the polar bivector of the bivector $\pip$.  It is an ideal line which is orthogonal (in the elliptic metric of the ideal plane, see \Sec{sec:ckc}) to the direction vector of $\pip$.  $(\pip_\infty; \pip_o)\eye = (\pip_o \eye; 0)$.
\end{enumerate}

\subsection{Dual Numbers}
\label{sec:dualno}
A number of the form $a + b\eye$ for $a, b \in \mathbf{R}$ we call a \emph{dual number}, after Study (\cite{study03}).  Dual numbers are similar to complex numbers, except $\eye^2=0$ rather than $i^2=-1$.  We'll need some results on dual numbers to calculate rotor logarithms below.  
\versions{
\begin{itemize}
\item Dual numbers commute with other elements of the Clifford algebra.
\item Given a dual number $\mathbf{z} = a+b\eye$, we say $\mathbf{z}$ is \emph{euclidean} if $a \neq 0$, otherwise $\mathbf{z}$ is \emph{ideal}.
\item \textbf{Conjugate.} Define the \emph{conjugate} $\overline{\mathbf{z}} = a-b\eye$. $\mathbf{z}\overline{\mathbf{z}} = a^2$. 
\item \textbf{Norm.} Define the \emph{norm} $\|\mathbf{z}\| := \sqrt{\mathbf{z}\overline{\mathbf{z}} }= a$.
\item \textbf{Inverse.} For euclidean $\mathbf{z}$, define the \emph{inverse} $\vec{z}^{-1} = \dfrac{\overline{\vec{z}}}{a^2}$.  The inverse is the unique dual number $\mathbf{w}$ such that $\mathbf{z}\mathbf{w} = 1$.
\correction \item \textbf{Square root.} Given a euclidean dual number $a+b\eye$, define $c=\sqrt{a}$ and $d = \dfrac{b}{2\sqrt{a}}$. Then $\mathbf{w} := c+d\eye$ satisfies $\mathbf{w}^2 = \mathbf{z}$ and we write $\mathbf{w} = \sqrt{\mathbf{z}}$.  
 \end{itemize}
}{
Dual numbers commute with other elements of the Clifford algebra.
Given a dual number $\mathbf{z} = a+b\eye$, we say $\mathbf{z}$ is \emph{euclidean} if $a \neq 0$, otherwise $\mathbf{z}$ is \emph{ideal}.
Define the \emph{conjugate} $\overline{\mathbf{z}} = a-b\eye$. Then $\mathbf{z}\overline{\mathbf{z}} = a^2$. Define the \emph{norm} $\| a+b\eye\| := \sqrt{\mathbf{z}\overline{\mathbf{z}} }= a$.
For euclidean $\mathbf{z}$, define the inverse $(a+b\eye)^{-1} = \dfrac{1}{a^2}({a-b\eye})$.  The inverse is the unique dual number $\mathbf{w}$ such that $\mathbf{z}\mathbf{w} = 1$.
Given a euclidean dual number $a+b\eye$, define $c=\sqrt{a}$ and $d = \dfrac{b}{2\sqrt{a}}$. Then $\mathbf{w} := c+d\eye$ satisfies $\mathbf{w}^2 = \mathbf{z}$ and we write $\mathbf{w} = \sqrt{\mathbf{z}}$.  
}

\textbf{Dual analysis.} Just as one can extend real power series to complex power series with reliable convergence properties, power series with a dual variable have well-behaved convergence properties.  See \cite{study03} for a proof. In particular, the power series for $\cos{(x+y\eye)}$ and $\sin{(x+y\eye)}$ have the same radii of convergence as their real counterparts.
One can use the addition formulae for $\cos$ and $\sin$ to show that (\cite{pottmann01}, p154):
\begin{eqnarray*}
\cos{(x+y\eye)} = \cos{x} - (y\eye) \sin{x}\\
\sin{(x+y\eye)} = \sin{x} + (y\eye)\cos{x}
\end{eqnarray*}

\textbf{The axis of a bivector.}
Working with euclidean bivectors is simplified by identifying a special line, the \emph{axis}, the unique euclidean line in the linear span of $\pip_\infty$ and $\pip_o$.  The axis $\pip_x$  is defined by  $\pip_x = (a+b\eye) \pip$ for a dual number $a+b\eye$.  In fact, one can easily check that the choice 
\versions{\[a: b = -2\langle \pip, \pip \eye  \rangle_P : \langle \pip, \pip \rangle_P\]}{\linebreak $a: b = -2\langle \pip, \pip \eye \rangle_P : \langle \pip, \pip \rangle_P$ }yields the desired simple bivector.  We usually normalize so that $\pip_x^{2} = -1$.
The axis appears later in the discussion of euclidean isometries in \Sec{sec:spinexplog3d}, since most isometries are characterized by a unique invariant axis.

\versions{
\exerc
\label{ex:exercise4}
\begin{enumerate}
\item Translate the products $\vec{a}\cdot \vec{b}, \vec{P}\cdot\vec{Q}, \vec{a}\wedge\vec{P}, \vec{a}\cdot\vec{P},  \vec{P} \vee \vec{Q}, \vec{P} \times \vec{Q}, \vec{a}\eye, \vec{P}\eye$ from \Sec{sec:enumpro}.  [Hint: a typical 2D item can be translated to 3D by replacing lines by planes, and 2-vectors with 3-vectors.  The grade of a product depends, as before, on the grades of the arguments.  For example, $\vec{a} \cdot \vec{P}$ is the line passing through $\vec{P}$ perpendicular to the plane $\vec{a}$.]   Also, translate the definition of norm to points and planes in 3D. 
\item \textbf{Relationships of lines.}
In this exercise, both $\pip$ and $\fio$ are normed euclidean simple bivectors.
\begin{enumerate}
\item \textbf{Distance between lines.} For normalized $\pip$ and $\fio$, $\pip \wedge \fio = \sin{(\alpha)} d_{\pip \fio}$, where $d_{\pip \fio}$ is the euclidean distance between the two lines, and $\alpha$ is the angle between their two direction vectors. [Hint: consider the tetrahedron spanned by unit vectors on $\pip$ and $\fio$.] 
\item \textbf{Dual norm and dual angle.} Define the \emph{dual norm} $\| \pip\|_d := \grade{\pip}{0} + \grade{\pip}{4}$. Then $\| \pip\fio\|_d = 0 \iff \pip$ and $\fio$ intersect at right angles.  In general,
\[   \| \pip\fio\|_d = \pm(\cos{(\alpha)} - \sin{(\alpha)} d_{\pip \fio} \eye \]
where $\alpha$ is the angle between the direction vectors of $\pip$ and $\fio$.  This is called the \emph{dual angle} of $\pip$ and $\fio$ and measures both the angle between the directions \textbf{and} the distance between the lines.   Give an algorithm to decide which choice of $\pm$ is correct.
\end{enumerate}
\item \textbf{The axis.} Let the axis of a non-simple $\pip = \pip_\infty + \pip_o$  be $\pip_x$.  Show that $\pip_x$ is the unique euclidean line such that $\pip_\infty$ is both the polar line of $\pip_x$ with respect to the euclidean metric (i.e., $\pip_\infty = \pip \eye$) \textbf{and} the conjugate line of $\pip_x$ with respect to the null polarity on $\pip$.
\item \textbf{Orthogonal projection.}  It's interesting to investigate orthogonal projectsions involving lines.  Not only can one project points onto lines, but lines can also be projected onto points.  To be more precise:  One can project points onto spears, and axes onto bundles.  Keep in mind in the following exercises that $\vec{X}^{-1}$ can in general be replaced by $\vec{X}$ without effecting the validity of the result, \emph{qua}  subspace.
\begin{enumerate}
\item \textbf{Projecting a point onto a line, and \emph{vice-versa}.} Show that for euclidean $\vec{P}$ and euclidean simple $\pip$, $(\vec{P} \cdot \vec{\pip})\pip^{-1}$ is the point of $\pip$ closest to $\vec{P}$.  What is $((\vec{P} \cdot \vec{\pip})\pip) \vee \vec{P}$?  Show that $(\vec{P} \cdot \vec{\pip})\vec{P}^{-1}$ is a line parallel to $\pip$ passing through $\vec{P}$.
\item \textbf{Projecting a line onto a plane, and \emph{vice-versa}..}  Show that for euclidean $\vec{a}$ and simple $\pip$, ($\vec{\pip} \cdot \vec{a})\vec{a}^{-1}$ is the orthogonal projection of $\pip$ onto $\vec{a}$. Show that  $(\vec{a} \cdot \vec{\pip})\vec{a}^{-1}$ represents a plane containing $\pip$,  parallel to $\vec{a}$.
\end{enumerate}
\item \textbf{The common normal of two lines.}
Let $\pip$ and $\fio$ be two simple euclidean bivectors. 
\begin{enumerate}
\item Show that $\thio := \pip \times \fio$ is a bivector which is in involution to both $\pip$ and $\fio$.  
\item Show that the axis $\thio_x$ of $\thio$ is a line perpendicular to both $\pip$ and $\fio$.
\item Calculate the points where $\thio_x$ intersects $\pip$ and $\fio$. [Hint: consider where the plane spanned by $\vec{P} := \e{0} \thio_x$ (the direction vector of $\thio_x$) and $\fio$ cuts $\pip$.]
\end{enumerate}
\item \hard Many of the products in \Sec{sec:enumpro3} are described only for simple $\pip$. Can you provide an interpretation for non-simple $\pip$?
\item Assume  $\pip \wedge \fio  = 0$ for simple bivectors $\pip$ and $\fio$, with $\pip \neq \fio$. Find the unique common point and unique common plane of $\pip$ and $\fio$. [Hint:  for a point $\vec{P} \notin \pip, \notin \fio$, consider $\pip \wedge (\fio \vee \vec{P})$.]
\item \textbf{The pitch of a bivector.} Let $\pip$ be a euclidean bivector and define the pitch $p$ to be the ratio: 
$ p := -\dfrac{\langle \pip, \pip \rangle_P}{\langle \pip, \pip \eye \rangle_P}$.
 To account for ideal $\pip$, one defines $p = \infty$.  Show that 
\begin{enumerate}
\item For simple euclidean $\pip$, $p=0$.
\item $p$ is a euclidean invariant. 
\end{enumerate}
\item For two general bivectors $\pip$ and $\fio$, use the ideal-finite decomposition $\pip = \pip_\infty + \pip_o$ and $\fio = \fio_\infty + \fio_o$ (\Sec{sec:propbiv}).   
Establish the following formulas:
\begin{enumerate}
\item $\pip \wedge \fio = \pip_\infty \wedge \fio_o + \pip_o \wedge \fio_\infty$.
\item $\pip \times \fio = (\pip_\infty \times \fio_o + \pip_o \times \fio_\infty; \pip_o \times \fio_o)$.
\item $\pip \cdot \fio = \pip_o \cdot \fio_o$.
\item $\pip$ is simple $\iff \pip_\infty \wedge \pip_o = 0$.
\end{enumerate}
\item Prove the following identities involving involving general bivectors, and an arbitrary 3-vector $\vec{P}$:
\begin{align} \label{eqn:ex4-1}
 (\pip \times \fio) \wedge \pip &= 0 \\ \label{eqn:ex4-2}
 \stripeye{(\pip\eye) \wedge \fio} =\stripeye{\pip \wedge (\fio\eye)} &= \pip \cdot \fio \\ \label{eqn:ex4-3}
 \pip \vee (\pip \wedge \vec{P}) &= \stripeye{\pip \wedge \pip}\vec{P} \\ \label{eqn:ex4-4}
 \pip \times \vec{P} &=(\pip \vee \vec{P})\eye \\ \label{eqn:ex4-5}
 \vec{P} \times (\pip \times \vec{P}) &= \vec{P} \pip \vec{P} + \pip 
\end{align}
\item \textbf{Incidence of point and line.}  Show that a euclidean point $\vec{P}$ lies on the simple bivector $\pip \iff \vec{P} \vee \pip = 0$.  Use \Eq{eqn:ex4-4} above to show this is equivalent to $\vec{P} \times \vec{\pip} = 0$. Show that in both cases the formulas are valid also for non-simple $\pip$.  (That is, for non-simple $\pip$, neither expression can vanish.)
\end{enumerate}
} {}

\subsection{Reflections, Translations, Rotations, and ... } 
The results of \Sec{sec:eucsand2d} can be carried over without significant change to 3D:
\begin{enumerate}
\item For a 1-vector  $\vec{a}$, the sandwich operation $\vec{P} \rightarrow \vec{a} \vec{P} \vec{a}$ is a euclidean reflection in the plane represented by $\vec{a}$.
\item For a pair of 1-vectors $\vec{a}$ and $\vec{b}$ such that $\vec{g} := \vec{a} \vec{b}$,   $\vec{P} \rightarrow  \vec{g}  \vec{P} \widetilde{\vec{g}}$ is a euclidean isometry.  There are two cases:
\begin{enumerate}
\item When $\langle \vec{g} \rangle_2$ is euclidean, it's a rotation around the line represented by $\langle \vec{g} \rangle_2$ by twice the angle between the two planes.
\item When $\langle \vec{g} \rangle_2$ is an ideal line $p_{01} \e{01} + p_{02} \e{02}+p_{03} \e{03}$, it's a translation by the vector $(x,y,z)=2(p_{01}, p_{02}, p_{03})$.  
\end{enumerate}
\end{enumerate}
A rotor responsible for a translation (rotation)  is called, as before, a \emph{translator} (\emph{rotator}). There are however other direct isometries in euclidean space besides these two types.  

\versions{
}{}

\begin{definition}
A \emph{screw motion} is a isometry that can be factored as a rotation around a line $\pip$ followed by a translation in the direction of $\pip$.  $\pip$ is called the \emph{axis} of the screw motion.
\end{definition}
Like the linear line complex, a screw motion has no counterpart in 2D.  In fact, $\langle \vec{g} \rangle_2$ is a non-simple bivector $\iff$ $\vec{g}$ is the rotor of a screw motion.  To show this we need to extend 2D results on rotors.

\subsection{Rotors, Exponentials and Logarithms}
\label{sec:spinexplog3d}

As in \Sec{sec:spinexplog}.  the spin group $\spin{3}{0}{1}$ is defined to consist of all elements $\vec{g}$  of the even subalgebra $\pdclplus{3}{0}{1}$ such that $\vec{g} \vec{\widetilde{g}} = 1$.  A group element is called a \emph{rotor}. In this section we seek the logarithm of a rotor $\vec{g}$.  Things are complicated by the fact that the even subalgebra  includes the pseudo-scalar $\eye$.   Dual numbers  help overcome this difiiculty.

Write $\vec{g} =  s_r  + s_d \eye+ \pip$.  Then
\[ \vec{g}\widetilde{\vec{g}} = s_r^2 + 2 s_r s_d \eye - \pip^2 = 1
\]
Suppose $\pip^2$ is real. Then $\pip$ is simple, $s_d = 0$, and $\pip^2 = s_r^2 - 1$. If $\pip^2  < 0$, then find real $\lambda$ such that $\pip_N := \lambda \pip$ satisfies $\pip_N^2 = -1$, and evaluate the formal exponential $e^{t\pip_N}$ as before (\Sec{sec:spinexplog})  to yield: 
\begin{equation}
e^{t\pip_N}  = \cos{(t)} + \sin{(t)} \pip_N
\end{equation}
We can use this formula to derive exponential and logarithmic forms for rotations  as in the 2D case (Exercise).  If $\pip^2  = 0$, $\pip$ is ideal, and the rotor is a translator, similar to the 2D case (Exercise).  This leaves the case $s_d \neq 0$.  Let $\fio = (a+b\eye)\pip$ be the axis of $\pip$ (see \Sec{sec:dualno} above).  Since the axis is euclidean, $a \neq 0$, and the inverse exists: $c+d\eye := (a+b\eye)^{-1}$.
\begin{equation} \label{eqn:axisform}
(c+d\eye)\fio = \pip
\end{equation}
Replace the real parameter $t$ in the exponential with a dual parameter $t+u\eye$, replace $\pip$ with $\fio$, substitute $\fio^2 = -1$, and apply the results above on dual analysis:
\begin{eqnarray}
e^{(t+u\eye)\fio} &=& \cos{(t+u\eye)} + \sin{(t+u\eye)} \fio \\ \label{eqn:exp3d}
&=& \cos{(t)}  - u \sin{(t)}\eye + (\sin{(t)} + u\cos{(t)}\eye)\fio
\end{eqnarray}
We seek values of $t$ and $u$ such that $\mathbf{g}$ equals the RHS of \Eq{eqn:exp3d}:
\begin{equation}
 g = s_r  + s_d \eye+ (c+d\eye)\fio =  \cos{(t)}  - u\sin{(t)}\eye + (\sin{(t)} + u\cos{(t)}\eye)\fio
\end{equation}
We solve for $t$ and $u$, doing our best to avoid numerical problems that might arise from $\cos{(t)}$ or $\sin{(t)}$ alone:
\begin{eqnarray*}
t &=& \text{tan}^{-1}{(c,s_r)}
\end{eqnarray*}
\[
u = 
\begin{cases}
\dfrac{d}{\cos{(t)}}	&\text{if $|\cos{(t)}| >  |\sin{(t)}|$} \\
\dfrac{-s_d}{\sin{(t)}}	&\text{otherwise}
\end{cases}
\]
\begin{definition}
Given a rotor $\mathbf{g} \in Sp(3,0,1)$ with non-simple bivector part, the bivector $(t+u\eye)\fio$ defined above is the \emph{logarithm} of $\mathbf{g}$.
\end{definition}
\begin{theorem} \label{thm:log}
Let $(t+u\eye)\fio$ be the logarithm of the rotor $\mathbf{g}$ $(= s_r + s_d \eye + \pip)$ with $s_d \neq 0$. Then $u \neq 0, t \neq 0$ and $\mathbf{g}$ represents a screw motion along  the axis $\fio$ consisting of rotation by angle $2t$ and translation by distance $2u$.
\end{theorem}
\begin{proof} 
$\fio$  commutes with $\pip$ and with $\mathbf{g}$ (Exercise). Hence $ \mathbf{g} \fio \widetilde{\mathbf{g}} = \fio$ is fixed by the sandwich. Write the sandwich operation on an arbitrary blade $\vec{x}$ as the composition of a translation followed by a rotation:
\begin{eqnarray*}
{\mathbf{g}} \vec{x} \widetilde{\mathbf{g}} &=& e^{(t+u\eye)\fio} \vec{x} e^{-(t+u\eye)\fio} \\
&=& e^{t\fio}(e^{u\eye\fio} \vec{x} e^{-u\eye\fio})e^{-t\fio}
\end{eqnarray*}
This makes clear the decomposition into a translation through distance $2u$  (Exercise), followed by a rotation around $\fio$ through an angle $2t$.  One sees that the translation and rotation commute by reversing the order.
\end{proof}

\versions{
%

Note that calculating the logarithm of a rotor involves two separate normalizations.  First, the rotor must be normalized to satisfy $\mathbf{g}\tilde{\mathbf{g}} = 1$.  Then the \emph{axis} of the bivector of $\mathbf{g}$ must also be extracted, a second normalization step.  This axis, multiplied by an appropriate \quot{dual angle} $t+u\eye$ can then  be exponentiated to reproduce the original rotor.   

\textbf{Translations.} Translations represent a degenerate case.  The logarithm of a translator is not unique, since 
\[ e^{u\eye\fio} = e^{u\eye(\fio + \pip)}\]
for any ideal bivector $\pip$.    This is related to the fact that a translator has no well-defined axis, since the pencil used to define the axis (\Sec{sec:dualno}) degenerates to the translator itself. But this degeneracy does not cause difficulties, since the calculation of exponential and logarithms in this case is simpler than the general case.  We choose the logarithm that goes through the origin of the coordinate system.
}{}

We have succeeded in showing that the bivector of a rotor is non-simple if and only if the associated isometry is a nondegenerate screw motion.  This result closes our discussion of   $\pdclal{3}{0}{1}$.   For a fuller discussion, see \fullversion. 

 \versions{ 
\exerc
\label{ex:exercise5}
\begin{enumerate}
\item Handle the case $s_d = 0$ ($\vec{g} =  s_r  + \pip$) from \Sec{sec:spinexplog3d} to obtain exponential forms for rotations and translations.  Define the corresponding logarithms. What does the case $s_r = s_d = 0$ represent?
\item Find the rotator corresponding to a rotation of $\dfrac{\pi}{3}$ radians around the line through the origin and the point $(1,1,1)$. [Answer:  $.5(1 + \e{12} + \e{31} + \e{23})$.]
\item From the proof of  Thm. \ref{thm:log}: 
\begin{enumerate}
\item Show that $\fio$  commutes with $\pip$ and with $\mathbf{g}$.
\item Confirm the claim of the theorem that the translation moves points a distance $2u$. [Hint: $\fio^2 = -1$ implies $\|  \eye \fio \|_\infty = 1$.]
\end{enumerate}
\end{enumerate}
}{}

\section{Case Study:  rigid body motion}
\label{sec:rbm}
The remainder of the article shows how to model euclidean rigid body motion using the Clifford algebra structures described above.  It begins by showing how to use the Clifford algebra $\pdclal{3}{0}{1}$ to represent  euclidean motions and their derivatives. Dynamics is introduced with newtonian particles, which are collected to construct rigid bodies. The inertia tensor of a rigid body is derived as a positive definite quadratic form on the space of bivectors. Equations of motion in the force-free case are derived.  In the following, we represent velocity states by $\velo$, momentum states by $\momo$, and forces by $\foro$.  \versions{}{Due to space limitations, results are compressed. For a fuller discussion, see \fullversion.} 

\subsection{Kinematics}
\label{sec:kinematics2d}

\begin{definition}
A \emph{euclidean motion} is a $C^1$ path $g: [0,1] \rightarrow \spin{3}{0}{1}$ with $g(0) = \one$. 
\end{definition}
\begin{theorem}\label{thm:bivector}
For  a euclidean motion $\vec{g}$, $\widetilde{\vec{g}}\dot{\vec{g}}$ is a bivector.
\end{theorem}
\begin{proof} 
$\widetilde{\vec{g}}\dot{\vec{g}}$ is in the even subalgebra. For a bivector $X$, $\tilde{X} = -X$; for scalars and pseudoscalars, $\tilde{X} = X$. Hence it suffices to show $\widetilde{\widetilde{\vec{g}}\dot{\vec{g}}}= -\widetilde{\vec{g}}\dot{\vec{g}}$.
\begin{eqnarray*}
\widetilde{\vec{g}} \vec{g}&=& 1 \\
\dot{(\widetilde{\vec{g}}\vec{g} )} &=& 0 \\
\dot{\widetilde{\vec{g}}}\vec{g}+\widetilde{\vec{g}}\dot{\vec{g}} &=& 0 \\
\widetilde{\dot{\vec{g}}}\vec{g}+\widetilde{\vec{g}}\dot{\vec{g}} &=& 0 \\
\widetilde{\widetilde{\vec{g}}\dot{\vec{g}}}&=& -\widetilde{\vec{g}}\dot{\vec{g}}
\end{eqnarray*}
\qed \end{proof}

Define $\vec{\velo} := \dot{\vec{g}}(0)$; by the theorem, $\vec{\velo}$ is a bivector.  We call $\vec{\velo}$ a \emph{euclidean velocity state}.  For a point $\vec{P}$, the motion $\vec{g}$ induces a path $\vec{P}(t)$, the \emph{orbit} of the point $\vec{P}$,  given by $\vec{P}(t) ={\vec{g}}(t)\vec{P} \widetilde{\vec{g}}(t)$.  
Taking derivatives of both sides and evaluating at $t=0$ yields:
\begin{eqnarray*} \label{eqn:liebr}
\dot{\vec{P}}(t) &=& \dot{{\vec{g}}}(t)\vec{P}\widetilde{\vec{g}}(t)+{\vec{g}}(t)\vec{P}\dot{\widetilde{\vec{g}}}(t)\\
\dot{\vec{P}}(t) &=&{\dot{\vec{g}}}(t)\vec{P} \widetilde{\vec{g}}(t)+{\vec{g}}(t)\vec{P} \widetilde{\dot{\vec{g}}}(t)\\
\dot{\vec{P}}(0) &=& \vec{\velo}\vec{P} - \vec{P}\vec{\velo}\\
&=& 2(\vec{\velo} \times \vec{P})
\end{eqnarray*}
The last step follows from the definition of the commutator product of bivectors.  In this formula we can think of $\vec{P}$ as a normalized euclidean point which is being acted upon by the euclidean motion $\vec{g}$.  From  \Sec{sec:enumpro3}  
we know that $\vec{\velo} \times \vec{P}$ is a ideal point, that is, a free vector.   
We sometimes use the alternative form $\vec{\velo} \times \vec{P} = (\vec{\velo} \vee \vec{P})\eye$ (Exercise). The vector field vanishes wherever $\vec{\velo} \vee \vec{P} = 0$. This only occurs if $\velo$ is a line and $\vec{P}$ lies on it. The picture is consistent with the knowledge, gained above, that in this case $e^{t\vec{\velo}}$ generates a rotation (or translation) with axis $\vec{\velo}$.   Otherwise the motion is an instantaneous screw motion around the axis of $\velo$ and no points remain fixed.  \versions{\Fig{fig:eucvf} shows how the vector field looks in case $n=2$. It's easy to see that $e^{t\vec{V}}$ in this case yields a rotation around the point $\vec{V}$.

\begin{figure}[t]
\sidecaption[t]
\includegraphics[width=.52\columnwidth]{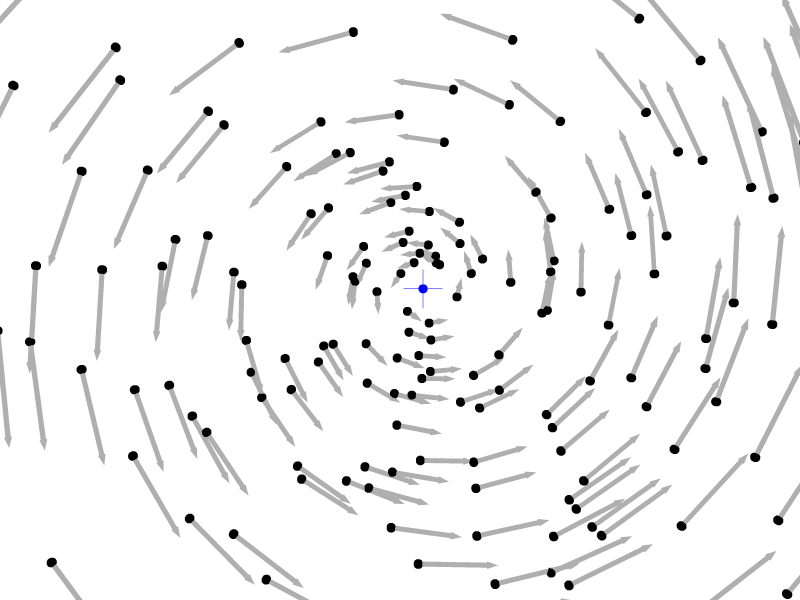}
\caption{For $n=2$ the euclidean velocity state is a point, and it acts on points. The vector field $2 \vec{V} \times \vec{P}$ in the neighborhood of $\vec{V}$.}
\label{fig:eucvf}
\end{figure}

}{\vspace{-.15in}}.

\mysubsubsection{Null plane interpretation} In the formulation $\dot{\vec{P}} =  2(\vec{\velo} \vee \vec{P})\eye$, we recognize the result as the polar point (with respect to the euclidean metric) of the null plane of $\vec{P}$ (with respect to $\velo$).  See \Fig{fig:nullplanes}. Thus, the vector field can be considered as the \emph{composition} of two simple polarities: first the null polarity on $\velo$, then the metric polarity on the euclidean quadric.   This leads to the somewhat surprising result that regardless of the metric used, the underlying null polarity remains the same.  \versions{One could say, for a given point, its null plane provide a projective \emph{ground} for kinematics, shared by all metrics; the individual metrics determine a different \emph{perpendicular} direction to the plane, giving the direction which the point moves.  This decomposition only makes itself felt in the 3D case.  In 2D, the null polarity is degenerate: $\vec{V} \vee \vec{P}$ is the joining line of $\vec{V}$ and $\vec{P}$ (a similar degeneracy occurs in 3D when $\velo$ is simple).

}{}

\begin{figure}[t]
\label{fig:nullplanes}
\begin{center}
\includegraphics[width=.9\columnwidth]{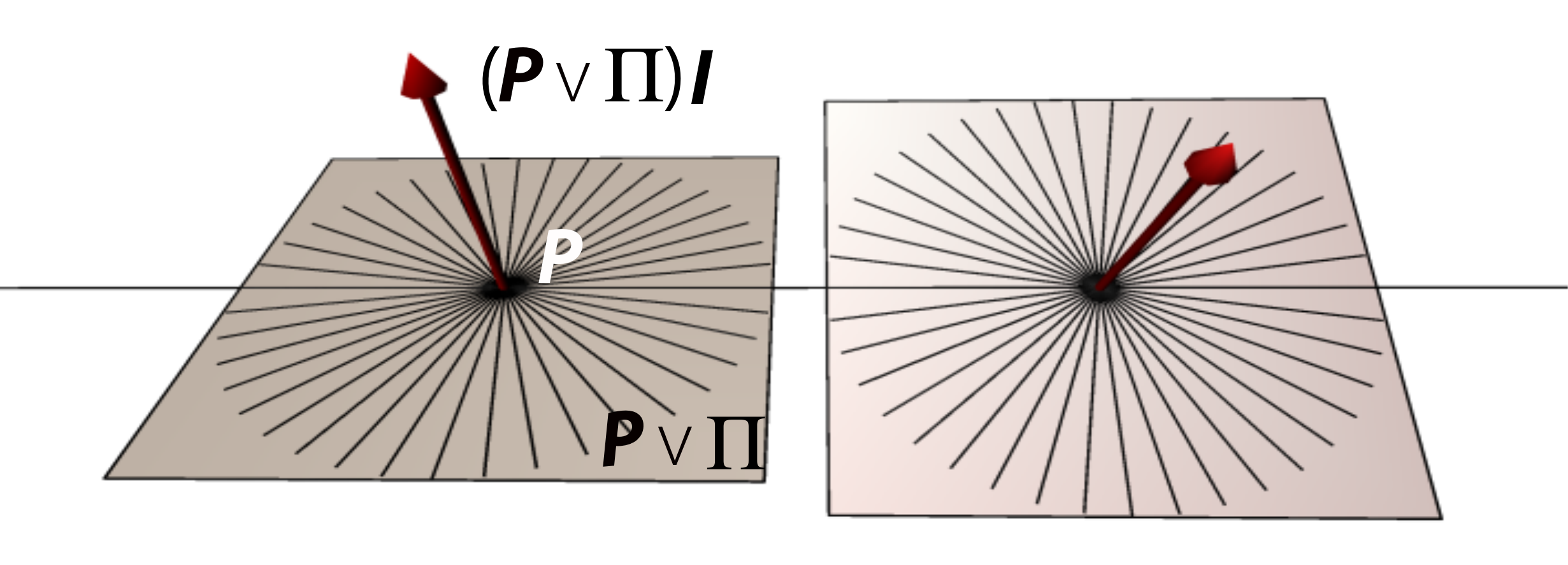}
\caption{Two examples of a point $\vec{P}$, its null plane $(\vec{P} \vee \pip)$, and the null plane's polar point.}
\end{center}
\end{figure}

\subsection{Dynamics}
\label{sec:dynamics2d}
With the introduction of forces, our discussion moves from the kinematic level to the dynamic one. We begin with a treatment of statics. We then introduce newtonian particles, build rigid bodies out of collection of such particles, and state and solve the equations of motions for these rigid bodies.

\versions{See Appendix 3 for a detailed account of how 2D statics are handled in this model.}{}

\mysubsubsection{3D Statics}
Traditional statics represents 3D a single force $F$ as a pair of 3-vectors $(V,M)$, where $ V = (v_x, v_y, v_z)$ is  the direction vector, and $M = (m_x,m_y,m_z)$ is the moment with respect to the origin (see \cite{featherstone07}, Ch. 2).   The resultant of a system of forces $F_i$ is defined to be \versions{\[ \sum_i F_i = (\sum_i{V_i} , \sum_i{M_i}) =: (V,M)\]  The forces are in equilibrium $\iff V=M=0$.  $V=0$ and $M \neq 0 \iff $ the resultant force is a \emph{force couple}. Otherwise the vectors $V$ and $M$ are orthogonal $\iff$ the system represents a single force.}{the sum of the corresponding direction vectors $V_i$ and moment vectors $M_i$.  The forces are in equilibrium if both terms of the resultant are zero.}

If $\vec{P}$ is a normalized point on the line carrying the force,  define $H(F) := \vec{P} \vee \vec{i}(V)$. We call $H(F)$ the\emph{ homogeneous form }of the force, and verify that:
\begin{eqnarray*}
H(F) &=&m_x \EE{01} + m_y \EE{02} + m_z \EE{03} + v_z \EE{12} +v_y \EE{31} + v_x \EE{23}
\end{eqnarray*}
If  $F$ is the resultant of a force system $\{F_i\}$, then $H(F) = \sum_i{H(F_i)}$.  Hence, a system of forces $\{F_i\}$ is the null force $\iff \sum_i{H(F_i)} = 0$.  Furthermore, $H(F)$ is an ideal line $\iff$ the system of forces reduces to a force-couple, and $H(F)$ is a simple euclidean bivector $\iff$ $F$ represents a single force.
Notice that the intensity of a bivector is significant, since it is proportional to the strength of the corresponding force.  For this reason we sometimes say forces are represented by \emph{weighted} bivectors.
 
 \subsubsection{Newtonian particles}
 \label{sec:newtonpart}
  The basic object of Newtonian mechanics is a particle $P$ with mass $m$ located at the point represented by a trivector $\vec{R}$.
Stated in the language of ordinary euclidean vectors, Newton's law asserts that the force $F$ acting on $P$ is: $F=m\ddot{\underline{\vec{R}}}$. 

\begin{definition}
The \emph{spear} of the particle is $\vec{\Lambda} := \vec{R} \vee \dot{\vec{R}}$.
\end{definition}
\begin{definition}
The \emph{momentum state} of the particle is $\vec{\momo} := m\vec{\Lambda}$.
\end{definition}
\begin{definition}
The \emph{velocity state} of the particle is $\vec{\pvelo} := \vec{\Lambda}\eye$. 
\end{definition}

\begin{definition}
The \emph{kinetic energy} $E$ of the particle is 
\begin{equation} \label{eqn:kinen}
E := \dfrac{m}{2} \| \dot{\vec{R}} \|_\infty^2 = -\frac{m}{2} \vec{\Lambda} \cdot \vec{\Lambda} = -\frac{1}{2}\stripeye{\vec{\pvelo} \wedge \vec{\momo} } 
\end{equation}
\end{definition}

\mysubsubsection{Remarks} Since we can assume $\vec{R}$ is normalized, $\dot{\vec{R}}$ is an ideal point. $\vec{\momo}$ is a weighted bivector whose weight is proportional to the mass and the velocity of the particle.  $\vec{\pvelo}$ is ideal, corresponding to the fact that the particle's motion is \emph{translatory}. Up to the factor $m$, $\pvelo$ is the polar line of $\momo$ with respect to the euclidean metric. It is straightforward to verify that the linear and angular momentum of the particle appear as $\vec{\momo}_o$ and $\vec{\momo}_\infty$, resp., and that the definition of kinetic energy agrees with the traditional one (Exercise). The second and third equalities in \Eq{eqn:kinen} are also left as an exercise.  

We consider only force-free systems.  \versions{}{The extension to include external forces is straightforward but lies outside the scope of this introduction.} 
\begin{theorem} If $F=0$ then $\vec{\Lambda}$, $\vec{\momo}$, $\vec{\pvelo}$, and $E$ are conserved quantities.
\end{theorem}
\begin{proof} 
$F=0$  implies $\ddot{\vec{R}} = 0$.  Then:
\begin{itemize}
\item $\dot{\vec{\Lambda}} =  ( \dot{\vec{R}} \vee \dot{\vec{R}} + \vec{R} \vee \ddot{\vec{R}}) =0~~$
\item $\dot{\vec{\momo}} = m\dot{\vec{\Lambda}} =0$
\item $\dot{\vec{\pvelo}} = (\dot{\vec{\Lambda}})\eye =0$
\item $\dot{E} = \dfrac{1}2(\stripeye{\dot{\vec{\pvelo}} \wedge \vec{\momo}} +\stripeye{\vec{\pvelo} \wedge \dot{\vec{\momo}}}) = 0$.
\end{itemize} 
\end{proof}

\mysubsubsection{Inertia tensor of a particle}
Assume the particle is \quot{governed by} a euclidean motion $\vec{g}$ with associated euclidean velocity state $\vec{\velo} := \dot{\vec{g}}(0)$.  Then $\vec{\momo}$, $\vec{\pvelo}$, and $E$ depend on $\vec{\velo}$ as follows: 
\begin{eqnarray} \label{eqn:inertiaparticle1}
\dot{\vec{R}} &=& 2(\vec{\velo} \times  \vec{R} ) \\  \label{eqn:inertiaparticle2} 
\vec{\momo} &=& 2m(\vec{R} \vee (\vec{\velo} \times  \vec{R} ))\\  \label{eqn:inertiaparticle3} 
\vec{\pvelo} &=&2(\vec{R} \vee (\vec{\velo} \times  \vec{R} ))\eye \\ \label{eqn:inertiaparticle4} 
&=&2(\vec{R} \times (\vec{\velo} \times  \vec{R} )) \\ \label{eqn:inertiaparticle5} 
E &=& -2m\stripeye{\vec{\pvelo} \wedge \vec{\momo}}\\  \label{eqn:inertiaparticle6} 
&=& -2m\stripeye{(\vec{R} \times (\vec{\velo} \times  \vec{R} )) \wedge (\vec{R} \vee (\vec{\velo} \times  \vec{R} ))}\\  \label{eqn:inertiaparticle7}
&=&  -\stripeye{\vec{\velo} \wedge \vec{\momo}}
\end{eqnarray}
\versions{The step from \Eq{eqn:inertiaparticle3} to \Eq{eqn:inertiaparticle4} follows from \Eq{eqn:ex4-5}.  The step from \Eq{eqn:inertiaparticle5} to \Eq{eqn:inertiaparticle7} is equivalent to the assertion that $2(\vec{\velo}\wedge\vec{\momo})= \vec{\pvelo} \wedge \vec{\momo}$.   
From $\vec{\momo} = m\vec{R} \vee \dot{\vec{R}}$ it is enough to show $2(\vec{\velo}\vee \vec{R}) = \vec{\pvelo} \vee \vec{R}$.   Since both sides of the equation are planes passing through $\vec{R}$, it only remains to show that the planes have the same normal vectors.  This is equivalent to 

\begin{theorem}\label{thm:simplerest}
$2(\vec{\velo} \times \vec{R}) = \vec{\pvelo} \times \vec{R}$. 
\end{theorem}
\begin{proof} 
\begin{align}
\vec{\pvelo} \times \vec{R} &= 2(\vec{R} \times (\vec{\velo} \times \vec{R}) \times \vec{R}) \\
&= 2(\vec{R} (\vec{\velo} \times \vec{R}) \vec{R}) \\
&= (\vec{R}(\vec{\velo} \vec{R} - \vec{R} \vec{\velo})\vec{R}) \\
&= (\vec{R}\vec{\velo}\vec{R}^2 - \vec{R}^2\vec{\velo}\vec{R}) \\
&= (-\vec{R}\vec{\velo} + \vec{\velo}\vec{R}) \\
&=2(\vec{\velo} \times \vec{R})
\end{align} 
Here we have used the fact that $\vec{R}\cdot (\vec{\velo} \times \vec{R})=0$, the definition of $\times$, the fact that for euclidean points $\vec{R}^2=-1$, and finally the definition of $\times$ a second time.
\qed
\end{proof} 
 
Define  Then the theorem yields immediately two corollaries:
\begin{corollary}\label{cor:restmom}
The difference $\vec{\Xi} := 2 \vec{\velo} - \vec{\pvelo}$ is a simple bivector incident with $\vec{R}$.
\end{corollary}
\begin{proof} The theorem implies $\vec{\Xi} \times \vec{R} = 0$.  Since this is the normal direction of the plane $\vec{\Xi} \vee \vec{R}$,  and $\vec{\Xi}$ is euclidean, this implies $\vec{\Xi} \vee \vec{R} = 0$. The null plane of a point, however, can only vanish if the bivector is simple and is incident with the point.  \qed 
\end{proof}
\begin{corollary}
$\vec{\Xi}$ is the conjugate line of $\vec{\pvelo}$ with respect to the null polarity $\vec{\velo}$.
\end{corollary}
\begin{proof}
This follows from the observation that the conjugate of a line $\fio$ with respect to a non-simple bivector $\vec{\velo}$ is a line lying in the pencil spanned by $\fio$ and $\vec{\velo}$.  This condition is clearly satisfied by both $\vec{\Xi}$ and $\vec{\pvelo}$.  Since there are at most two such lines in the pencil, the proof is complete. \qed \end{proof}
  
 We return to the theme of Newtonian particles below in \Sec{sec:rbm2}.

  }{The step from \Eq{eqn:inertiaparticle4} to \Eq{eqn:inertiaparticle5} is described in more detail in \fullversion.}
Define a real-valued bilinear operator $\inert$ on pairs of bivectors:
\begin{eqnarray} \label{eqn:symm1}
\inert(\vec{\velo}, \vec{\pip}) &:=&  -\frac{m}2\stripeye{((\vec{R} \vee 2(\vec{\velo} \times  \vec{R} ))\eye) \wedge (\vec{R} \vee 2(\vec{\pip} \times  \vec{R} ))} \\ \label{eqn:symm2}
&=& \frac{m}2(\vec{R} \vee 2(\vec{\velo} \times  \vec{R} )) \cdot (\vec{R} \vee 2(\vec{\pip} \times  \vec{R} ))
\end{eqnarray}
where the step from \Eq{eqn:symm1} to \Eq{eqn:symm2} can be deduced from \Sec{sec:enumpro3}.  \Eq{eqn:symm2} shows that $\inert$ is symmetric since $\cdot$ on bivectors is symmetric: $\vec{\Lambda} \cdot \vec{\Delta} = \vec{\Delta} \cdot \vec{\Lambda}$.  We call $\inert$ the \emph{inertia tensor} of the particle, since $E =\inert(\vec{\velo},\vec{\velo}) =  -\stripeye{\vec{\velo} \wedge \vec{\momo}} $. We'll construct the inertia tensor of a rigid body out of the inertia tensors of its particles below.  We overload the operator and write $\vec{\momo} = \inert(\vec{\velo})$ to indicate the polar relationship between $\vec{\momo}$ and $\vec{\velo}$. 

\subsubsection{Rigid body motion}
\label{sec:rbm2}
Begin with a finite set of mass points $P_i$; for each derive the velocity state $\vec{\pvelo}_i$, the momentum state $\vec{\momo}_i$, and the inertia tensor $\inert_i$.\footnote{We restrict ourselves to the case of a finite set of mass points, since extending this treatment to a continuous mass distribution presents no significant technical problems; summations have to be replaced by integrals. }  Such a collection of mass points is called a \emph{rigid body} when the euclidean distance  between each pair of points is constant. 

Extend the momenta and energy to the collection of particles by summation:
\begin{subequations}
\label{eqn:lwe}
\begin{eqnarray}
\vec{\momo}&:=& \sum{\vec{\momo}_i}  =  \sum\inert_i(\vec{\velo}) \label{eqn:lwe2} \\
E&:=& \sum{E_i} = \sum{\inert_i(\vec{\velo}, \vec{\velo})} 
\end{eqnarray}
\end{subequations}

Since for each single particle these quantities are conserved when $F=0$, this is also the case for the aggregate $\vec{\momo}$ and $E$ defined here.

We introduce the inertia tensor $A$ for the body:
\begin{definition}
$\inert:=\sum{\inert_i}$.
\end{definition}

Then $\vec{\momo} = \inert(\vec{\velo})$ and $E =  \inert(\vec{\velo},\vec{\velo})$, neither formula requires a summation over the particles: the shape of the rigid body has been encoded into $\inert$. \versions{We sometimes use the identity
\begin{equation} \label{eqn:dotwedge}
\inert(\vec{\velo}_1,\vec{\velo}_2) = -\stripeye{\vec{\velo}_1 \wedge \inert(\vec{\velo}_2)}~~~~\text{(Exercise)}
\end{equation}
which is a consequence that the individual inertia tensors for each particle exhibit this property.

}{}One can proceed traditionally and diagonalize the inertia tensor by finding the center of mass and moments of inertia (see \cite{arnold78}).  Due to space constraints we omit the details.  Instead, we sketch how to integrate the inertia tensor more tightly into the Clifford algebra framework.

\mysubsubsection{Clifford algebra for inertia tensor}
We define a Clifford algebra $\mathbf{C}_\inert$ based on $\pdgrassgr{4}{2}$ by attaching the positive definite quadratic form $\inert$ as the inner product.\footnote{It remains to be seen if this approach represents an improvement over the linear algebra approach which could also be maintained in this setting.}  We denote the pseudoscalar of this alternative Clifford algebra by $\eye_\inert$, and inner product of bivectors by $\langle, \rangle_\inert$.  We use the same symbols to denote bivectors in $W^*$ as 1-vectors in $\mathbf{C}_\inert$.  Bivectors in $W$ are represented by 5-vectors in $\mathbf{C}_\inert$. Multiplication by $\eye_\inert$ swaps 1-vectors and 5-vectors in $\mathbf{C}_\inert$; we use $\mathbf{J}$ (lifted to  $\mathbf{C}_\inert$) to convert  5-vectors back to 1-vectors as needed.
 The following theorem, which we present without proof, shows how to obtain $\vec{\momo}$ directly from $\eye_\inert$ in this context:
\begin{theorem} \label{thm:inert}
Given a rigid body with inertia tensor $\inert$ and velocity state $\vec{\velo}$,  the momentum state $\vec{\momo} = \inert(\vec{\velo}) = \mathbf{J}(\vec{\velo \eye_\inert })$.
\end{theorem}
Conversely, given a momentum state $\momo$, we can manipulate the formula in the theorem to deduce:
\begin{equation*}
\vec{\velo} =\inert^{-1}(\vec{\momo}) = (\mathbf{J}(\vec{\momo})\eye_\inert^{-1})
\end{equation*}
In the sequel we denote the polarity on the inertia tensor by $\inert(\vec{\velo})$ and $\inert^{-1}(\vec{\momo})$, leaving open whether the Clifford algebra approach indicated here is followed.
\versions{

\mysubsubsection{Newtonian particles, revisited} Now that we have derived the inertia tensor for a euclidean rigid body, it is instructive to return to consider the formulation of euclidean particles above (\Sec{sec:newtonpart}).  We can see that in this formulation, particles exhibit properties usually associated to rigid bodies.  
\begin{itemize}
\item $E = -\frac12\stripeye{\vec{\pvelo} \wedge \vec{\momo}}$: The kinetic energy is the result of a dual pairing between the particle's velocity state and its momentum state, considered as bivectors.
\item $\vec{\pvelo} = \frac1m \vec{\momo} \eye$: The dual pairing is given by the polarity on the euclidean metric quadric, scaled by $\frac1m$. This pairing is degenerate and only goes in one direction: from the momentum state to produce the velocity state. 
\item $E = -\stripeye{\vec{\velo} \wedge \vec{\momo}}$: the same energy is obtained by using twice the global velocity state in place of the particle's velocity state.  This follows from \Thm{thm:simplerest}.
\end{itemize}
\exerc
\begin{enumerate}
\item Verify that the linear and angular momentum of a particle appear as $\vec{\momo}_o$ and $\vec{\momo}_\infty$, resp.
\item Verify the equalities in \Eq{eqn:kinen}.
\end{enumerate}
}{\vspace{-.2in}}

\subsubsection{The Euler equations for rigid body motion}
In the absence of external forces, the motion of a rigid body is completely determined by its momentary velocity state or momentum state at a given moment.   How can one compute this motion? First we need a few facts about coordinate systems.

\versions{\textbf{Coordinate systems.} Up til now, we have been considering the behavior of the system at a single, arbitrary moment of time.  But if we want to follow a motion over time, then there will be two natural coordinate systems.  One, the \emph{body} coordinate system, is fixed to the body and moves with it as the body moves through space. The other, usually called the \emph{space} coordinate system, is the coordinate system of an unmoving observer.  Once the motion starts, these two coordinate systems diverge.}{}  The following discussion assumes we observe a system as it evolves in time.  All quantities are then potentially time dependent;  instead of writing $\vec{g}(t)$, we continue to write $\vec{g}$ and trust the reader to bear in mind the time-dependence.

We use the subscripts $X_s$ and $X_c$\footnote{From \emph{c}orpus, Latin for body.} to distinguish whether the quantity $X$ belongs to the space or the body coordinate system.  The conservation laws of the previous section are generally valid only in the space coordinate system,  for example,  $\dot{\vec{\momo}}_s = 0$. On the other hand, the inertia tensor will be constant only with respect to the body coordinate system, so, $\vec{\momo}_c = \inert( \vec{\velo}_c) $.  When we consider a euclidean motion $\vec{g}$ as being applied to the body, then 
the relation between body and space coordinate systems for \emph{any} element $\vec{X} \in \pdclal{3}{0}{1}$, with respect to a motion ${\vec{g}}$, is given by the sandwich operator:
\[
\vec{X}_s =  \vec{g}  \vec{X}_c\widetilde{\vec{g}}
\]
\begin{definition}
The \emph{velocity in the body} $\vec{\velo}_c :=\tilde{\vec{g}} \dot {\vec{g}}$, and the  \emph{velocity in space} $\vec{\velo}_s :=\vec{g} \vec{\velo}_c  \tilde{\vec{g}}$.
\end{definition}
\begin{definition}
The \emph{momentum in the body} $\vec{\momo}_c :=\inert( \vec{\velo}_c)$, and the \emph{momentum in space} $\vec{\momo}_s := \vec{g} \vec{\momo}_c \tilde{\vec{g}}$.
\end{definition}

We derive a general result for a time-dependent element (of arbitrary grade) in these two coordinate systems:
\begin{theorem} \label{thm:liebracket}
For  time-varying $\vec{X} \in \pdclal{3}{0}{1}$ subject to the motion $\vec{g}$ with velocity in the body $\vec{\velo}_c$,
\[
\dot{\vec{X}}_s = \vec{g}(\dot{\vec{X}}_c+2(\vec{\velo}_c \times \vec{X}_c )) \tilde{\vec{g}}
\]
\end{theorem}
\begin{proof}
\begin{eqnarray*}
\dot{\vec{X}_s} &=&\dot{\vec{g}}\vec{X}_c\tilde{\vec{g}}+\vec{g}\dot{\vec{X}}_c \tilde{\vec{g}}+ \vec{g}\vec{X}_c \dot{\tilde{\vec{g}}}\\
&=&\vec{g}(\tilde{\vec{g}}\dot{\vec{g}}\vec{X}_c+\dot{\vec{X}}_c+ \vec{X}_c\dot{\tilde{\vec{g}}}{\vec{g}}) \tilde{\vec{g}}\\
&=& \vec{g}({\vec{\velo}}_c\vec{X}_c+\dot{\vec{X}_c}+\vec{X}_c \widetilde{\vec{\velo}}_c) \tilde{\vec{g}}\\
&=& \vec{g}(\dot{\vec{X}}_c+ \vec{\velo}_c \vec{X}_c - \vec{X}_c\vec{\velo}_c) \tilde{\vec{g}} \\
&=& \vec{g}(\dot{\vec{X}}_c+2(\vec{\velo}_c \times \vec{X}_c )) \tilde{\vec{g}}
\end{eqnarray*}
The next-to-last equality follows from the fact that for bivectors, $\widetilde{\vec{\velo}} = - \vec{\velo}$; the last equality is the definition of the commutator product.  \qed 
\end{proof}


We'll be interested in the case $\vec{X}_c $ is a bivector. In this case,  $\vec{X}_c$ and $\vec{\velo}_c$ can be considered as Lie algebra elements, and $2(\vec{\velo}_c \times \vec{X}_c)$ is called  the \emph{Lie bracket}. It expresses the change in one ($\vec{X}$) due to an instantaneous euclidean motion represented by the other ($\vec{\velo}$).
\subsubsection{Solving for the motion}

 Since $\vec{\velo}_c =\widetilde{\vec{g}} \dot{\vec{g}}$, $\dot{\vec{g}} = \vec{g} \vec{\velo}_c$, a first-order ODE.   If we had a way of solving for $\vec{\velo}_c$,  we could solve for $\vec{g}$.   If we  had a way of solving for $\vec{\momo}_c$,  we could apply \Theorem{thm:inert} to solve for $\vec{\velo}_c$. So, how to solve for $\vec{\momo}_c$?

We apply the corollary to the case of force-free motion.  Then 
$\dot{\vec{\momo}}_s = 0$: the momentum of the rigid body in space is constant. 
By \Theorem{thm:liebracket},
\begin{equation}\label{eqn:dmzero}
0 = \dot{\vec{\momo}}_s = {\vec{g}}(\dot{\vec{\momo}}_c+2 (\vec{\velo}_c \times \vec{\momo}_c)) \widetilde{\vec{g}}
\end{equation}
The only way the RHS can be identically zero is that the expression within the parentheses is also identically zero, implying:
\[
\dot{\vec{\momo}}_c = 2 \vec{\momo}_c \times \vec{\velo}_c
\]
Use the inertia tensor to convert velocity to momentum yields a differential equation purely in terms of the momentum:
\begin{eqnarray*}
\dot{\vec{\momo}}_c &=& 2 \vec{\momo}_c \times \inert^{-1}(\vec{\momo_c})
\end{eqnarray*}
\versions{One can also express this ODE in terms of the velocity state alone:
\begin{eqnarray*}
\dot{\vec{\velo}}_c = \inert^{-1}(\dot{\vec{\momo}}_c) &=& 2\inert^{-1}(\vec{\momo}_c \times \inert^{-1}(\vec{\momo_c}) )\\
&=& 2\inert^{-1}(\inert(\vec{\velo}_c) \times \vec{\velo}_c )
\end{eqnarray*}
\vspace{-.2in}}{}

When the inner product is written out in components, one arrives at the well-known Euler equations for the motion (\cite{arnold78}, p. 143).

The complete set of equations for the motion $g$ are given by the pair of first order ODE's:
\begin{eqnarray}\label{eqn:eqnmot}
\dot{\vec{g}} &=& \vec{g}\vec{\velo}_c\\
\dot{\vec{\momo}}_c &=& 2 \vec{\momo}_c \times \vec{\velo}_c 
\end{eqnarray}
where $\vec{\velo}_c =  \inert^{-1}(\vec{\momo}_c )$. When written out in full, this gives a set of 14 first-order linear ODE's.  The solution space is 12 dimensions; the extra dimensions corresponds to the normalization $\vec{g}\widetilde{\vec{g}}=1$.   At this point the solution continues as in the traditional approach, using standard ODE solvers.  Our experience is that the cost of evaluating the Equations (\ref{eqn:eqnmot}) is no more expensive than traditional methods.   \versions{}{For an extension to account for external forces, see \fullversion.}

\versions{
\subsubsection{External forces}
\label{sec:extforce}
The external force $F$ acting on the rigid body is the sum of the external forces $F_i$ acting on the individual mass points. In traditional notation, $F_i = m_i\ddot{\mathbf{r}}$. What is the homogeneous form for $F_i$?  In analogy to the velocity state of the particle, we define the acceleration \emph{spear} of the particle $\vec{\Upsilon}_i := \vec{R}_i \vee \ddot{\vec{R}}_i$, and the force state of the particle $\vec{\foro}_i : = m_i\vec{\Upsilon}_i$. Then the total force $\vec{\foro}$ is the sum:
\begin{eqnarray*}
\vec{\foro} &=&\sum{\vec{\foro}_i} \\
\end{eqnarray*}

The following theorem will be needed below in the discussion of work (\Sec{sec:work}):
\begin{theorem}
\label{thm:force}
$\vec{\foro} = \dot{\vec{\momo}}$
\end{theorem}
\begin{proof}
Take the derivative of \Eq{eqn:lwe2}:
\begin{eqnarray*}
\dot{\vec{\momo}} &=&  \sum{m_i \frac{d}{dt}(\vec{R}_i \vee \dot{\vec{R}}_i}) \\
&=&  \sum{m_i (\dot{\vec{R}}_i \vee \dot{\vec{R}}_i + \vec{R}_i \vee \ddot{\vec{R}}_i)} \\
&=&  \sum{m_i (\vec{R}_i \vee \ddot{\vec{R}}_i)}\\
&=& \vec{\foro}
\end{eqnarray*}\qed \end{proof}

Theorem \ref{thm:force} applied to  \Eq{eqn:dmzero} yields Euler equations for motion with external forces::
\begin{eqnarray*}
\vec{\foro}_s &=& \dot{\vec{\momo}}_s =  \vec{g}(\dot{\vec{\momo}}_c+2 ( \vec{\velo}_c \times \vec{\momo}_c))\widetilde{\vec{g}} \\
\widetilde{\vec{g}} \vec{\foro}_s\vec{g}&=&\dot{\vec{\momo}}_c+2 ( \vec{\velo}_c \times \vec{\momo}_c) \\
\vec{\foro}_c&=&\dot{\vec{\momo}}_c+2 ( \vec{\velo}_c \times \vec{\momo}_c) \\
\dot{\vec{\momo}}_c&=&\vec{\foro}_c +2 (\vec{\momo}_c \times \vec{\velo}_c)
\end{eqnarray*}
Note that the forces have to converted from world to body coordinate systems.

\subsubsection{Work}
\label{sec:work}
As a final example of the projective approach, we discuss the concept of \emph{work}.  Recall that $\vec{\momo} = \inert( \vec{\velo})$, so $\dot{\vec{\momo}} = \inert(\dot{\vec{\velo}})$, and the definition of kinetic energy for a rigid body : $E = \frac{m}2 \inert(\vec{\velo}, \vec{\velo})$.

\begin{theorem}
$\dot{E} =-\stripeye{\vec{\velo} \wedge \vec{\foro}}$
\end{theorem}
\begin{proof}
\begin{eqnarray*}
\dot{E} &=& \frac{m}2\frac{d}{dt} \inert(\vec{\velo}, \vec{\velo}) \\
& = & \frac{m}2 (  \inert(\dot{\vec{\velo}}, \vec{\velo}) +\inert(\vec{\velo}, \dot{\vec{\velo}})) \\
&=&  m\inert( \vec{\velo}, \dot{\vec{\velo}})\\
&=&   -\stripeye{\vec{\velo} \wedge m\inert(\dot{\vec{\velo}})}=  -\stripeye{\vec{\velo} \wedge m\dot{\vec{\momo}} }\\
&=& -\stripeye{ \vec{\velo} \wedge \foro}
\end{eqnarray*} where we apply Leibniz rule, the symmetry of $\inert$, \Eq{eqn:dotwedge} and finally Thm. \ref{thm:force}. \qed \end{proof}

In words: the rate of change of the kinetic energy is equal to the signed magnitude of the outer product of force and velocity.  This is noteworthy in that it does not involve the metric directly.

$\dot{E}$ is sometimes called the \emph{power}.  The \emph{work} done by the force between time $t_0$ and $t$ is the defined to be the integral:
\begin{eqnarray*}
w(t) = E(t) - E(t_0) &=& \int_{t_0}^{t}\dot{E} ds \\
&=& \int_{t_0}^{t}-\stripeye{\vec{\velo} \wedge \vec{\foro}}ds \\
\end{eqnarray*}
 The integrand depends on the incidence properties of the force and the velocity, as a line/point pair.    If the two elements are incident, then $\vec{f} \wedge \vec{\velo} = 0$ and there is no work done; the further away $\vec{\velo}$ lies from $\vec{\foro}$, (in $\complexspace$!) the greater the power and hence the work done.  One can be more precise:
 \begin{theorem}
 \label{thm:forcedist}
 Suppose $\vec{\foro} $ is a single force,  $\vec{\velo} $ is a rotator, $d$ is the euclidean distance between the lines $\vec{\foro}$ and  $\vec{\velo}$, and $\alpha$ is the angle between the two direction vectors. Then
$ \stripeye{\vec{\velo} \wedge \vec{\foro}}=  -d \sin{\alpha} \| \vec{\velo} \| \| \vec{\foro} \|.$
 \end{theorem}
 \begin{proof} Left as an exercise. \qed \end{proof}
 In  words: the rate of change of the kinetic energy is proportional to the intensity of the force, the intensity of the velocity, and the euclidean distance between the force line and the velocity line.  

\textbf{Example.}
Imagine an 
ice skater who moving along the  surface of a frozen lake with negligible friction; the single force is given by gravity. Assuming gravity is in the negative $z$-direction acting on the skater located at the origin, then $\vec{\foro} = g\e{12}$ (corresponding to the intersection line of the planes $x=0$ and $y=0$ with weight gravitation constant $g$).  Consider two possible motions of the skater:
\begin{itemize}
\item The motion of the skater is a translation in the $x$-direction given by an ideal line of the form $\vec{\velo}=d\e{01}, d<0$.  $\vec{\foro} \wedge \vec{\velo} = 0$, so no work is required for the skater to skate!   
\item The skater spins around a fixed point.  Then the velocity state relative to the natural diagonalized form of the inertia tensor has null ideal part $\vec{\velo}_\infty = 0$ and the corresponding momentum state $\vec{\momo} = \inert(\vec{\velo})$ has null euclidean part $\vec{\momo}_o = 0$: it's a momentum couple: a momentum carried by an ideal line!
\item As the skater spins, she stretches her arms out, then pulls her arms close to her body.  This latter movement decreases the entries in the inertia tensor $\inert$, increasing the entries in $\inert^{-1}$; since $\vec{\velo} = \inert^{-1}(\vec{\momo})$,  her velocity increases proportionally in intensity: she spins faster.
\end{itemize}
One can see from this example the advantages of the projective approach: it unifies calculations, and handles exceptional cases, such as the translations and couples in the above example, at no extra cost.
}{
}

\mysubsubsection{Comparison}
The projective Clifford algebra approach outlined here exhibits several advantages over other approaches to rigid body motion.   The representation of kinematic and dynamic states as bivectors rather than as pairs of ordinary 3-vectors (linear and angular velocity/momentum/etc.) provides a framework free of the special cases which characterize the split approach (for example, translations are rotations around an ideal line, a force couple is a single force carried by an ideal line, etc.).   The Clifford algebra product  avoids cumbersome matrix formulations and, as seen in \Thm{thm:liebracket}, yields formulations which are valid for points, lines, and planes uniformly.  The inertia tensor of a rigid body can be represented as a separate (positive definite) Clifford algebra defined on the space of bivectors.  Finally, the treatment of Newtonian particles reveals an underlying velocity-momentum polarity in bivector space analogous to that of rigid bodies.
\vspace{-.1in}

\versions{
\section{Extension to elliptic and hyperbolic metrics}

The theoretical basis for modeling rigid body motion in noneuclidean spaces using linear line complexes was worked out in the $19^{th}$ century and the earlier part of the $20^{th}$ century (\cite{lindemann73}, \cite{heath84}, \cite{francesco02}, \cite{blaschke42}). The treatment given here for euclidean geometry and physics can be extended to elliptic and hyperbolic space by changing the defining signature to be $(4,0,0)$ for elliptic and $(3,1,0)$ for hyperbolic space.  Surprisingly little else needs to be changed.   \Fig{fig:hyperbolicMotion} shows a simulation of rigid body motion in hyperbolic space using the Clifford algebra $Cl(3,1,0)$.  A detailed description of the phenomenology lies outside the scope of this article.

\begin{figure}
\sidecaption[t]
\label{fig:hyperbolicMotion}
\includegraphics[width=.42\columnwidth]{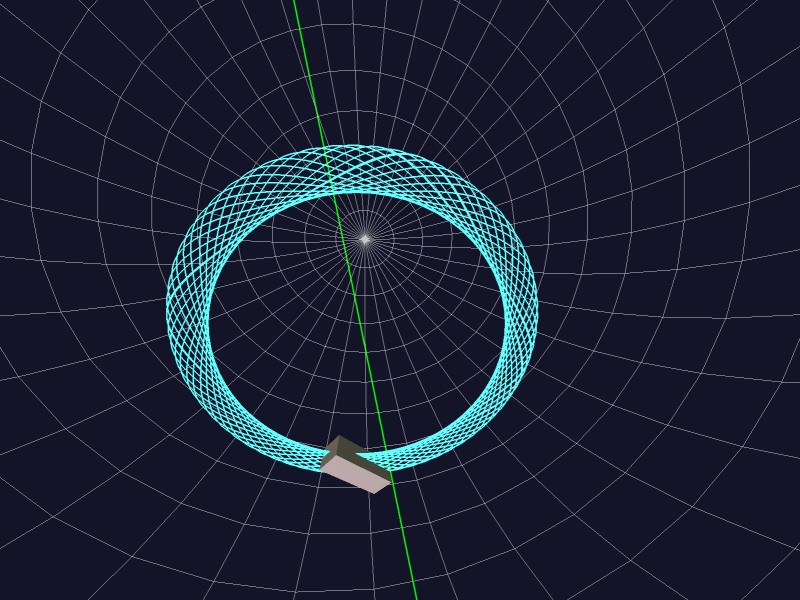}
\caption{Rigid body motion in hyperbolic space. Note that the motion is bounded but has no fixed point, instead seems to fill out a solid torus. The momentum state is a hyper-ideal line which lies completely outside hyperbolic space.  }
\end{figure}

In this way one obtains a setting in which euclidean geometry is positioned as a limiting case of  geometries with nondegenerate metric;  the connections to its neighbors sheds new light on euclidean geometry and physics.  Since Riemann, we know that the metric of space can only be learned through empirical research.  
The difference to nondegenerate metric spaces is particularly interesting in the case of rigid body motion.  Here one becomes aware of peculiarities of euclidean space.  

For example, one considers a rigid body as composed of points.  But in a nondegenerate metric, every point has a unique polar plane which kinematically and dynamically cannot be separated from it -- moves along with it, as it were.  The one-sided focus on points has no validity here.  One of Lindemann's important contributions was, in this setting, to state the underlying dynamics of a rigid body, not in terms of the instantaneous lines of translation of points, but in terms of the instantaneous axes of rotation of planes.  Another difference: the concept of \emph{center of mass} is specific to euclidean space; in elliptic space a rigid body has 4 centers.  But space constraints prevent further discussion of this interesting theme.

}{}

\section{Guide to the Literature}

The work presented here, in its conceptual basis, is derived from \cite{study03}.  
This seminal  book worked out  in  impressive detail the structure which in modern form appears as the even subalgebra $\pdclplus{3}{0}{1}$.  Study avoided using the term quaternion;  the structure he described has nonetheless become known as the \emph{dual quaternions}.  His parameter $\epsilon$ maps to the pseudoscalar $\eye$.  A full description of the correspondence between the two systems lies outside the scope of this work, nor is the full scope of \cite{study03} reflected in material presented here.
 \cite{ziegler} gives an excellent survey of the historical development that led up to Study's work. M\"{o}bius, Pl\"{u}cker, Hamilton and Klein were Study's most important predecessors; he and Ball (\cite{ball00}) had a relationship based on mutual appreciation. Weiss, a student of Study's, wrote a concise introduction \cite{weiss35} to Study's investigations which can be recommended.    For beginners, \cite{blaschke54a} provides a simpler introduction to Study's approach.

Study's contribution in the direction of mechanics was developed further in \cite{vm24} and \cite{blaschke42}. The former was hampered by the awkward matrix notation required for stating transformation laws.  These, on the other hand, are \quot{built in} to the Clifford algebra approach and simplify the approach considerably.  \cite{blaschke42} concentrated on the noneuclidean case. 

The modern legacy of this work is varied.  Some modern literature, such as \cite{featherstone07}, use \emph{spatial vectors} to model rigid body motion; these are 6D vectors equivalent to our bivectors, but developed within a linear algebra framework reminiscent of \cite{vm24}. 
While the dual quaternions have developed a following in robotics and other applied areas (\cite{mcc})\mymarginpar{reference?}, their embedding as the even subalgebra of $\pdclal{3}{0}{1}$ has not received widespread acceptance.  Modern sources using $\pdclal{3}{0}{1}$ include \cite{selig00} and \cite{selig05}. 
Much contemporary work which applies Clifford algebra methods to physics and engineering uses   the conformal model (\cite{doran03}, \cite{perwass09}, \cite{hestenes10}) rather than the homogeneous model presented here.  Our handling of rigid body motion owes much to the spirit of \cite{arnold78}.

Study himself appears to be aware of the possibility of extending his work within a more comprehensive algebraic structure.  He remarked at the end of \cite{study03} (p. 595, translated by the author): 
\begin{quote}
The elementary geometric theory, that hovers thus before us, will surpass the construction possibilities of the quaternions, to the same degree that the geometry of dynamen surpasses the addition of vectors. The accessory analytic machinery will consist of a system of compound quantities, with eight, or better yet, with \emph{sixteen} units. \{Study's italics!\}
\end{quote}
It seems likely that $\pdclal{3}{0}{1}$ is in fact the 16-dimensional algebraic realizations of Study's prophetic inkling.  
 
 \versions{
 
 \section{Further work}

\subsection{Loose ends}
There are some loose ends in the formulation presented here, that would be good to tidy up:
\begin{itemize}
\item \textbf{Orientation.}  The derivation of the reflection in \Sec{sec:eucsand2d} shows that, before dehomogenizing, the reflection actually flips the sign of the transformed coordinates.  In particular the homogeneous coordinate is -1, not 1.  One needs to first establish that this phenomenon is generally present in indirect versors (generally, odd-grade versors) and also for other types of operands (sandwich filling), not just points.  If it is generally true, it would be interesting, since such a reversal of sign can be interpreted as a change of orientation due to the reflection.  That would imply that dehomogenizing should be applied selectively, in order to retain this information.  For background on this theme see \cite{stolfi91}.  Previous work on this theme (\cite{pappas96}) used metric polarity to implement duality instead of the projective approach used here, but presumably the results obtained do not depend on this fact.
\item \textbf{Software implementations.} The past decade has seen a explosion of software which implements geometric algebras of various sorts.  Such software is generally of two kinds: implementations in so-called general purpose languages such as C++ or Java, and those based on symbolic algebra systems such as Mathematica or Maple.  Examples of the former include \cite{gaigen}, \cite{clu}, and \cite{gaalop}; and of the latter, \cite{aragon}. These generally have not been written with degenerate metrics in mind.  Consequently, either there is no support for degenerate metrics, or it hasn't been debugged.  As the discussion has made clear, we recommend that the implementation of what is called \emph{duality} to be based on projective principles, not on a metric polarity (Appendix 1).  This typically should not present a serious problem for software systems, since true duality can in fact be simulated using the metric polarity.  One needs also to establish to what extent the assumption of a non-degenerate metric has been \quot{hard-wired} into the code, for example, in the existence of $\eye^{-1}$.  The author adapted \cite{aragon} to support the degenerate metric used here;  these changes were not difficult to implement.  In correspondence with the authors of the package, it appears these changes will be integrated in the next release of the package. 
\end{itemize}

\subsection{Visualization}
One essential task is the visualization of non-simple bivectors.  As the above hopefully makes clear,  bivectors are the essential ingredient to almost every non-trivial phenomena in 3D kinematics and dynamics.  Visualization tools for non-simple bivectors are however not well-developed.  Naively drawing all lines of the linear line complex associated to the bivector is not practical.  If one takes a metric as given, then one can draw the axis of the bivector, along with some extra information to identify which element of the bivector pencil is intended.  One also need to consider ways to distinguish between axes and spears visually.  

\section{Evaluation and comparison}
We collect and present observations on the use of the homogeneous model of euclidean geometry and its Clifford algebraic implementation, which arose during the course of the work.  
\begin{itemize}
\item \textbf{Fits the phenomena.} The mathematical theory of rigid body motion maps precisely onto the projective model presented here.  Indeed, a series of renowned mathematicians in the $19^{th}$ century (Pl\"{u}cker, Klein, Clifford, Study, and others), in their research on rigid body motion,  laid the foundations of the mathematical theory  which later developed into modern Clifford algebras (\cite{klein72a}, \cite{study03}). 
\item \textbf{Duality.} Our implementation provides duality via the projective foundation of the homogeneous model, rather than attempting to obtain it using a metric.  Appendix 1 makes clear why this projectively valid method is often confused with a metric polarity.
\item \textbf{Ideal elements} play an important role in the homogeneous model of euclidean space.  
In incidence calculations involving parallel elements, they provide correct answers.  Free vectors are represented as ideal points.  An ideal line represents, in statics,  a force couple; in kinematics, the axis of a euclidean translation, and plays an important role throughout the theory presented here. 
\item \textbf{Efficient.} Before optimizing, a multiplication in the homogeneous model is 4 times cheaper than in the conformal model.  Integration of equations of motion is easier in the homogeneous model since there are fewer directions to wander off the integral surface.   The only possibility for \quot{wandering off} when solving the equations of rigid body motion is normalizing the motion $\mathbf{g}$ so that $\mathbf{g}\tilde{\mathbf{g}}=1$.
\item \textbf{Intuitive.} The geometric representation in the homogeneous model agrees with the naive sense perceptions of the human being.  Straight is straight.
\item\textbf{Insider's view} of non-euclidean spaces is naturally computed. This important visualization feature of the projective model lies outside the scope of this article; a detailed discussion can be found in \cite{gunnmn2010}.
\item \textbf{Degenerate metric.} Rather than considering this a disadvantage,  our experience leads us to appreciate the benefits it brings.  This \textbf{is} the metric of euclidean space, and it provides correct formulas and results for euclidean geometry.  While the non-invertiblity of $\eye$ is unfamiliar,  it presents no difficulties in the calculations and procedures represented here. 
\end{itemize}

}{}

\section{Conclusion}
We have established that $\pdclal{n}{0}{1}$ is a model for euclidean geometry. By using a mixture of projective, ideal and properly euclidean elements, we have avoided the problems traditionally associated with degenerate metrics.  The result provides a complete and compact representation of euclidean geometry, kinematics and rigid body dynamics.  We hope that this work will stimulate others working in these fields, to consider the homogeneous model as a practical solution to their problem domain, and to deepen and extend these initial results.  

\versions{
}{}

\versions{
\section*{Appendix 1: $\mathbf{J}$, metric polarity, and the regressive product}
\addcontentsline{toc}{section}{Appendix: \mathbf{J}}
\label{sec:J}

Because of its importance in our approach to the homogeneous model, we provide here an more stringent, dimension-independent treatment of the map $\mathbf{J}$ and its close connection to the polarity $\mathbf{\Pi}$ on the elliptic metric, concluding with reasons to prefer the use of $\mathbf{J}$ to that of $\mathbf{\Pi}$ for implementing duality in the Clifford algebra setting.

\mysubsubsection{Canonical basis} A subset $S = \{i_1,i_2,...i_k\}$ of $N = \{1,2,...,n\}$ is called a \emph{canonical \emph{k}-tuple} of $N$ if  $i_1<i_2<...<i_k$.  For each canonical $k$-tuple of $N$, define $S^\perp$ to be the canonical $(n-k)$-tuple consisting of the elements $N \setminus S$.  For each unique pair $\{S,S^\perp\}$,  swap a pair of elements of $S^\perp$ if necessary so that the concatenation $SS^\perp$, as a permutation $P$ of $N$, is even.  Call the collection of the resulting sets $\mathfrak{S}$.  For each $S \in \mathfrak{S}$, define $\vec{e}^S = \vec{e}^{i_1}...\vec{e}^{i_k}$.    We call the resulting set $\{\vec{e}^S\}$ the \emph{canonical basis} for $W$ generated by $\{\vec{e}^i\}$.

Case 1: Equip $W = \proj{\bigwedge (\R{n})}$ with the euclidean metric to form the Clifford algebra $\pclal{n}{0}{0}$ with pseudoscalar $I = \vec{e}^N$.   Then, by construction, $ \vec{e}^S\eye = \vec{e}^{(S^\perp)}$ (remember $(\vec{e}^S)^{-1} = \vec{e}^S$).
This is the polarity $\Pi: W \rightarrow W$ on the elliptic metric quadric (see \Sec{sec:ckc}). 

Case 2: Consider $W^*$, the dual algebra to $W$. Choose a basis $\{\vec{e}_1, \vec{e}_2, ... \vec{e}_n\}$ for $W^*$ so that $\vec{e}_i$ represents the same oriented subspace represented by the $(n-1)$-vector $\vec{e}^{(i^\perp)}$ of $W$.  Construct the canonical basis (as above) of $W^*$ generated by the basis $\{\vec{e}_i\}$.  Then define a map $\mathbf{J}: W \rightarrow W^*$ by $\mathbf{J}(\vec{e}^S) = \vec{e}_{S^\perp}$ and extend by linearity.  
$\mathbf{J}$ is an \quot{identity} map on the subspace structure of $\RP{n}$: it maps a $k$-blade $B \in W$ to the $(n-k)$-blade $\in W^*$ which represents the same geometric entity as $B$ does in $\RP{n}$. \emph{Proof:} By construction, $\vec{e}^{S}$ represents the join of the 1-vectors $\vec{e}_{i_j}, (i_j \in S)$ in $W$.  This is however the same subspace as the meet of the $n-k$ basis 1-vectors $\vec{e}_{i_j}, (i_j \in S^\perp)$ of $W^*$, since $\vec{e}_i$ contains $\vec{e}^j$ exactly when $j \neq i$.

\mysubsubsection{Conclusion}: Both $\mathbf{J}$ and $\Pi$ represent valid grade-reversing involutive isomorphisms.  
 The \textbf{only} difference is the target space: $\Pi: \proj{W}\rightarrow \proj{W}$, while $\mathbf{J}: \proj{W} \rightarrow \proj{W^*}$. 

\mysubsubsection{The regressive product via a metric} Given the point-based exterior algebra  with outer product $A \wedge B $, (representing the join of subspaces $A$ and $B$), the \emph{meet} operator $A \vee B$  (also known as the regressive product) is often defined $\Pi(\Pi(A) \wedge \Pi(B))$ \cite{hessob87}. That is, the euclidean metric (any nondegenerate metric suffices) is introduced 
in order to provide a solution to a projective (incidence) problem. By the above, the same result is also given via $\mathbf{J}(\mathbf{J}(A) \wedge \mathbf{J}(B)$.   (Here, the $\wedge$ denotes the outer product in $W^*$).

\mysubsubsection{The Hodge $\star$ operator}
An equivalent method for producing dual coordinates is described in \cite{pottmann01}, p. 150.  The $\star$ operator is presented as a way of generating \emph{dual} coordinates, which is an apt description of the $\mathbf{J}$ operator also.   One can also define the regressive product by $\star(\star A \wedge \star B)$.  Formally, however, since $\star$ is a map from $W$ to $W$, it is identical to the metric polarity $\Pi$.

\mysubsubsection{Comparison} The two methods yield the same result, but they have very different conceptual foundations.  As pointed out in \Sec{sec:equalrights},  the meet and join operators live the exterior algebra $W$ and its dual.  $\mathbf{J}$ provides the bridge between these two projective algebras, hence provides a projective \emph{explanation} for what is a projective operation.  A related advantage of $\mathbf{J}$ is that it also is useful when used \quot{alone} -- whereas $\Pi$ (and also $\star$) is only valid when it appears in the second power, and hence disappears.  For these reasons, we propose a differentiation of the terminology to reflect this mathematical differentiation.  Instead of referring to multiplication by $\eye$ (or $\eye^{-1}$) as the \emph{duality} operator, we propose it should be called the \emph{metric polarity} operator 
This is consistent with the mathematical literature. The purely projective term \emph{duality} would be reserved for the $\mathbf{J}$ operator. 
}{}

\versions{
\section*{Appendix 2: The euclidean distance function as a limit}
\addcontentsline{toc}{section}{Appendix 3: 2: The euclidean distance function}

Rewrite \Eq{eq:epsmetric} and manipulate:
 \begin{align}
 \cos^2{(d_\epsilon)}\langle \mathbf{X}, \mathbf{X} \rangle_\epsilon\langle \mathbf{Y}, \mathbf{Y} \rangle_\epsilon &= \langle \vec{X},\vec{Y} \rangle_\epsilon \\
 \cos^2{(d_\epsilon)}(\epsilon + x_1^2)(\epsilon+y_1^2) &= (\epsilon+ (x_1 y_1))^2 \\
 \cos^2{(d_\epsilon)}(\epsilon^2 +\epsilon ( ( y_1)^2 + (x_1)^2) + (x_1 y_1)^2) &= \epsilon^2 + 2 \epsilon(x_1 y_1) + (x_1 y_1)^2 
 \end{align}
 Now consider the limit as $\epsilon\rightarrow \infty$. It's clear from  \Eq{eq:epsmetric} that $\lim_{\epsilon\rightarrow \infty}  \cos{d_\epsilon}= 1$.  So we can replace $\cos{d_\epsilon}^2$ by $(1-d_\epsilon^2)$, and simplify the resulting expressions:
  \begin{align}
(1-d_\epsilon^2)(\epsilon^2 +\epsilon ( y_1^2 + x_1^2) + (x_1 y_1)^2) &=  \epsilon^2+ 2 \epsilon( x_1 y_1) + (x_1 y_1)^2 \\
 \epsilon ( y_1^2 + x_1^2 - 2x_1 y_1)  &=  d_\epsilon^2(\epsilon^2+\epsilon ( y_1^2 + x_1^2) + (x_1 y_1)^2)\\
 (x_1 - y_1)^2  &= d_\epsilon^2(\epsilon+( y_1^2 + x_1^2) + \frac1{\epsilon}(x_1 y_1)^2)
  \end{align}

Notice that the LHS is exactly the desired euclidean metric.  In order to make the limit converge to this, we define a new distance function $\hat{d}_\epsilon := \sqrt{\epsilon}d_\epsilon$.  This is exactly the scaling needed to prevent the distance from going to zero in the limit. One obtains:
\begin{align}
\lim_{\epsilon \rightarrow \infty} \hat{d}_\epsilon^2 &= \frac{ (x_1 - y_1)^2  }{(1+ \frac1{\epsilon}( y_1^2 + x_1^2)+ \frac1{\epsilon^2}(x_1 y_1)^2 }\\
&= (x_1 - y_1)^2 
 \end{align}
}{}

\versions{
\section*{Appendix 3: 2D Statics}
\addcontentsline{toc}{section}{Appendix 3: 2D Statics}
\label{sec:2dstatics}
We include here a detailed account of euclidean planar statics, to show how the classical separation of forces into a linear and an angular part is overcome through the use of the projective model.  

Statics is the study of systems of rigid bodies in equilibrium, or close to equilibrium.\footnote{Dynamics on the other hand is concerned with an analysis of the \emph{motion} which results when a system is not in equilibrium.} This discussion is restricted to systems of forces acting on a \emph{single} rigid body.  In two dimensional statics, forces are characterized by the following geometrical properties:
\begin{enumerate}
\item A single force  can be represented as a vector free to move along a given lin.  
One speaks of \emph{the vector associated to the force}, and 
of the \emph{carrying line} of the force.  
\item Two forces acting on lines which intersect, are equivalent to a single force acting on a line passing through the intersection of the two lines, whose vector is equal to the vector sum of the vectors of the two original forces. See Figure  \ref{fig:statics}.
\item Two forces acting along the same line are equivalent to a single force acting on this line, whose vector is the vector sum of the two vectors.
\end{enumerate}
If the vector of a force vanishes, we say it is the null force.

\begin{figure}
\centering
\includegraphics[width=.3\columnwidth]{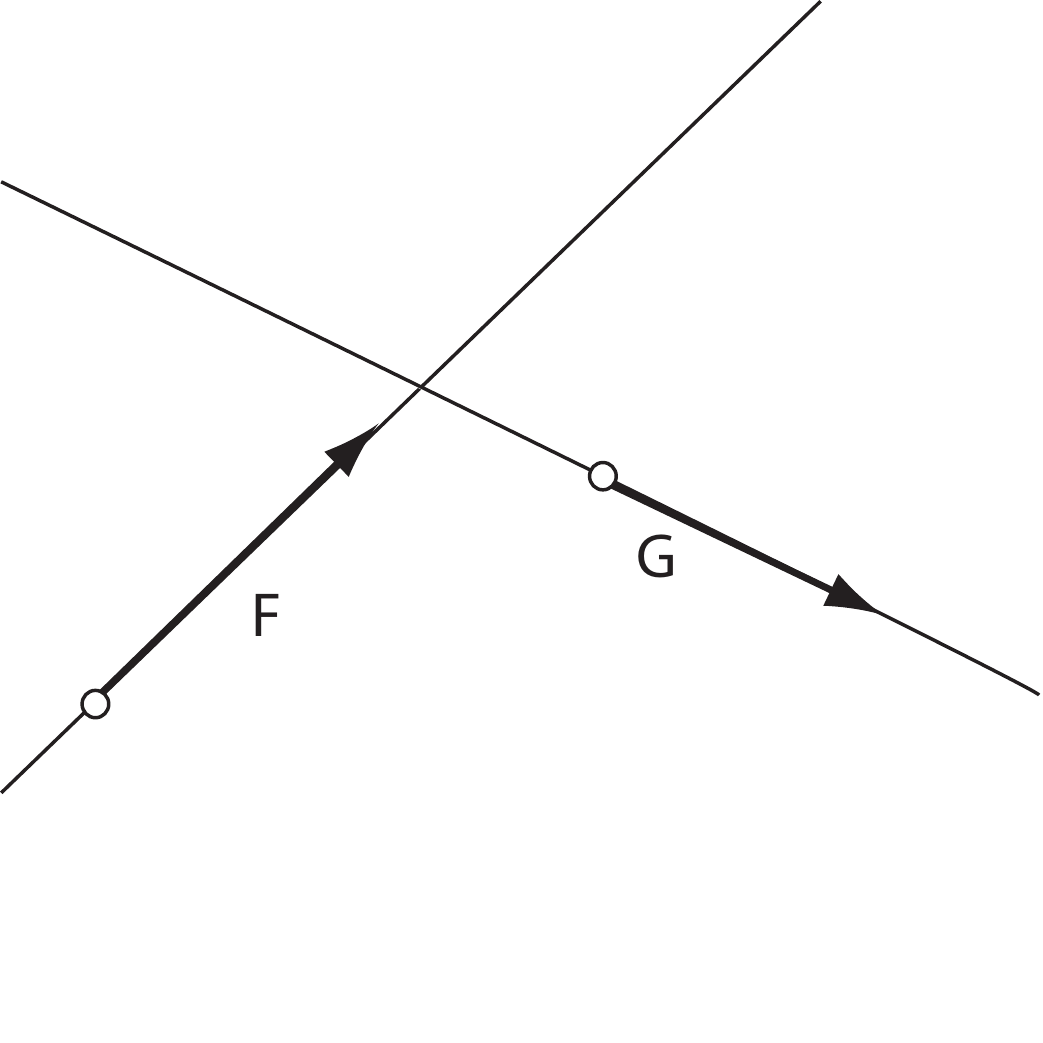}%
\hspace{.05in}
\includegraphics[width=.3\columnwidth]{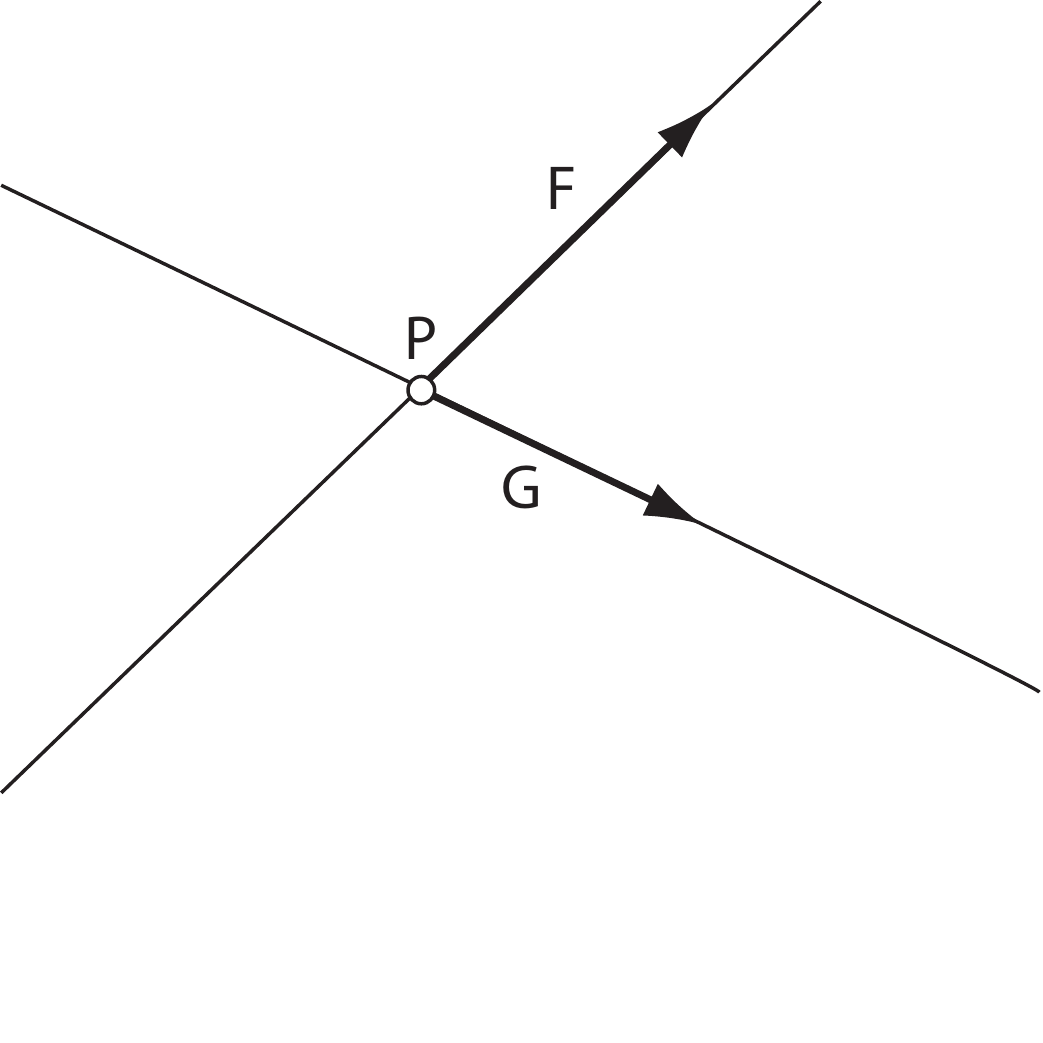}%
\hspace{.05in}
\includegraphics[width=.3\columnwidth]{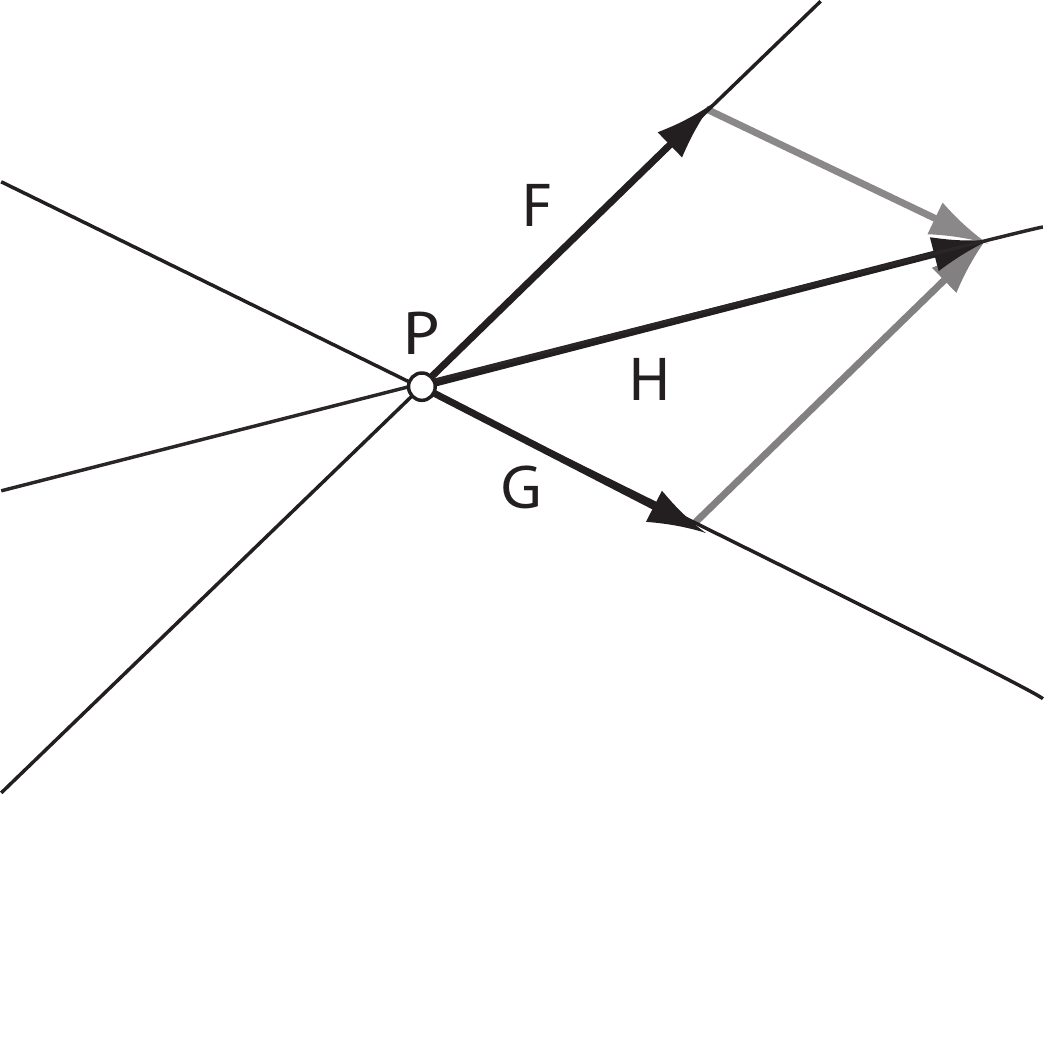}%
\vspace{-.3in}
\caption{The two forces F and G are equivalent to the force H.}
\label{fig:statics}
\end{figure}

We recall some standard results of traditional treatments of 2D statics.  Suppose $\vec{v} = (v_x, v_y)$ is the vector of a force, $P = (p_x, p_y)$ is a point on the carrying line. 
\begin{definition}
The \emph{moment}  of a force $F$ with respect to the origin is $M : = p_x v_y - p_y v_x$.
\end{definition}
The main result of planar statics takes the following form:  
\begin{theorem}
A system of forces $F_i$ is equivalent to the null force $\iff$ 
\[
\sum_i{\vec{v}_i} = \sum_i{M_i} = 0
\]
where $\vec{v}_i$ is the vector of $F_i$ and $M_i$ is the moment of $F_i$ with respect to the origin.
\end{theorem}

Let's find a homogeneous version of the above!  Given a force $F$ , with $\vec{P}$ and $\vec{Q}$ two normalized points on the line carrying the force, such that the force vector is given by $\vec{Q} - \vec{P} = (v_x, v_y)$.  Define the 1-vector $H(F) : = \vec{P} \vee \vec{Q}$.   
\begin{theorem}
A system of forces $F_i$ is equivalent to the null force $\iff$ $\sum_i{H(F_i)} = 0$.
\end{theorem}
\begin{proof}
For normalized points $\vec{P}$ and $\vec{Q}$, $\vec{P} := \EE{0} + p_x \EE{1} + p_y \EE{2}$ and 
$ \vec{Q} := \EE{0} + q_x \EE{1} + q_y \EE{2}$, 
\begin{eqnarray} \label{eqn:force}
\vec{P} \vee \vec{Q} &=& (p_x q_y - p_y q_x)\e{0} - (q_y-p_y) \e{1} + (q_x-p_x)\e{2} 
 \end{eqnarray}
Note first that the coordinates of the force vector occur (rotated as it were by 90 degrees) as the coefficients of $\e{1}$ and $\e{2}$.  A simple computation shows that the coefficient of $\e{0}$ is the moment of the force $m$:
\begin{eqnarray*}
p_x q_y - p_y q_x &=& p_x(p_y+v_y) - p_y(p_x+v_x) \\
&=& p_x v_y - p_y v_x
\end{eqnarray*}
Hence $\vec{P} \vee \vec{Q} = m \e{0} + v_y \e{1} - v_x \e{2}$. The summand in this theorem vanishes $\iff$ both summands in the previous one vanish. \qed \end{proof}
 
 In brief:  the blade determined by the force vector's head and tail,  contains both force vector \textbf{and} moment with respect to the origin.  This motivates the following definition:
\begin{definition}
 Let $F$ be a force with vector $(v_x, v_y)$ and moment $m$. Then $H(F) = m\e{0} + v_y\e{1} - v_x\e{2}$ is called the \emph{homogeneous form} of $F$.
\end{definition}
 \textbf{Remarks.} As the theorem shows us, the homogeneous form replaces two equations with a single one. 
 Furthermore, it provides a unified form which subsumes a sticky exception in euclidean statics, the so-called \emph{force couple}, which arises as the resultant of two forces with equal-and-opposite force vectors carried on parallel lines.  There is no single force equivalent to a force couple, so all results in euclidean statics general contain references to both single forces and to force couples.  However, as is easy to verify by a simple calculation, the resultant of two such forces in homogeneous form yields a resultant of the form $c \e{0}$, $c\neq0$. This corresponds to the ideal line:  just as a translation is integrated in the projective model as a rotation around a ideal point, a force couple appears in the projective model as a force carried by the ideal line, whose intensity is proportional to the distance between the two original lines, and to the magnitude of the original force vectors.
 }{}

\bibliography{GunnRef}

\begin{thebibliography}{Gun11b}

\bibitem[Arn78]{arnold78}
V.~I. Arnold.
\newblock {\em Mathematical Methods of Classical Physics}.
\newblock Springer-Verlag, New York, 1978.
\newblock Appendix 2.

\bibitem[Bal00]{ball00}
Robert Ball.
\newblock {\em A Treatise on the Theory of Screws}.
\newblock Cambridge University Press, Cambridge, 1900.

\bibitem[Bla38]{blaschke38}
Wilhelm Blaschke.
\newblock {\em Ebene Kinematik}.
\newblock Tuebner, Leibzig, 1938.

\bibitem[Bla42]{blaschke42}
Wilhelm Blaschke.
\newblock {\em Nicht-euklidische Geometrie und Mechanik}.
\newblock Teubner, Leipzig, 1942.

\bibitem[Bla54]{blaschke54a}
Wilhelm Blaschke.
\newblock {\em Analytische Geometrie}.
\newblock Birkhauser, Basel, 1954.

\bibitem[Bou89]{bourbaki89}
Nicolas Bourbaki.
\newblock {\em Elements of mathematics, Algebra I}.
\newblock Springer Verlag, 1989.

\bibitem[Con08]{conradt08}
Oliver Conradt.
\newblock {\em Mathematical Physics in Space and Counterspace}.
\newblock Verlag am Goetheanum, 2008.

\bibitem[Cox78]{cox:neg}
H.M.S. Coxeter.
\newblock {\em Non-euclidean Geometry}.
\newblock University of Toronto, Toronto, 1978.

\bibitem[Cox87]{cox:pg}
H.M.S. Coxeter.
\newblock {\em Projective Geometry}.
\newblock Springer-Verlag, New York, 1987.

\bibitem[dF02]{francesco02}
Domenico de~Francesco.
\newblock Sui moto di un corpo rigido in uno spazio di curvatura constante.
\newblock {\em Mathematische Annalen}, 55:573--584, 1902.

\bibitem[DFM09]{dfm07}
Leo Dorst, Daniel Fontijne, and Stephen Mann.
\newblock {\em Geometric Algebra for Computer Science}.
\newblock Morgan Kaufmann, San Francisco, 2009.

\bibitem[DL03]{doran03}
Chris Doran and Anthony Lasenby.
\newblock {\em Geometric Algebra for Physicists}.
\newblock Cambridge University Press, Cambridge, 2003.

\bibitem[Fea07]{featherstone07}
Roy Featherstone.
\newblock {\em Rigid Body Dynamics Algorithms}.
\newblock Springer, 2007.

\bibitem[Fon]{gaigen}
Daniel Fontijne.
\newblock Gaigen 2.5: A code generator for geometric algebra.
\newblock \url{http://sourceforge.net/projects/g25}.

\bibitem[Gun10]{gunnmn2010}
Charles Gunn.
\newblock Advances in {M}etric-neutral {V}isualization.
\newblock In {\em GraVisMa 2010}, pages 17--26, Brno, 2010. Eurographics,
  \url{http://gravisma.zcu.cz/GraVisMa-2010/GraVisMa-2010-proceedings.pdf}.

\bibitem[Gun11a]{gunn2011}
Charles Gunn.
\newblock On the homogeneous model for euclidean geometry.
\newblock In {\em Applications of Geometric Algebra in Computer Science and
  Engineering}, pages ?--? Springer, 2011.

\bibitem[Gun11b]{gunnFull2010}
Charles Gunn.
\newblock On the homogeneous model for euclidean geometry: Extended version.
\newblock {\em \url{http://arxiv.org/abs/1101.4542}}, 2011.

\bibitem[Hea84]{heath84}
R.~S. Heath.
\newblock On the dynamics of a rigid body in elliptic space.
\newblock {\em Phil. Trans. Royal Society of London}, 175:281--324, 1884.

\bibitem[Hes10]{hestenes10}
David Hestenes.
\newblock New tools for computational geometry and rejuvenation of screw
  theory.
\newblock In Eduardo~Jose Bayro-Corrochano and Gerik Scheuermann, editors, {\em
  Geometric Algebra Computing: in Engineering and Computer Science}, pages
  3--35. Springer, 2010.

\bibitem[Hil]{gaalop}
Dieter Hildebrand.
\newblock Gaalop: Geometric algebra algorithms optimizer.
\newblock \url{http://www.gaalop.de/}.

\bibitem[Hit03]{nh:pg}
Nigel Hitchin.
\newblock {\em Projective Geometry}.
\newblock
  \url{http://people.maths.ox.ac.uk/hitchin/hitchinnotes/Projective_geometry/C%
hapter_3_Exterior.pdf}, 2003.

\bibitem[HS87]{hessob87}
David Hestenes and Garret Sobczyk.
\newblock {\em Clifford Algebra to Geometric Calculus}.
\newblock Fundamental Theories of Physics. Reidel, Dordrecht, 1987.

\bibitem[JLAV]{aragon}
et.~al. JosŽ Luis Arag—n~Vera.
\newblock clifford.m.
\newblock \url{http://http://www.fata.unam.mx/aragon/software}.

\bibitem[Kle72]{klein72a}
Felix Klein.
\newblock {\"{U}}ber liniengeometrie und metrische geometrie.
\newblock {\em Mathematische Annalen}, 2:106--126, 1872.

\bibitem[Kle49]{klein:neg}
Felix Klein.
\newblock {\em Vorlesungen {\"{U}}ber Nicht-euklidische Geometrie}.
\newblock Chelsea, New York, 1949.
\newblock (Original 1926, Berlin).

\bibitem[Kow09]{kowol2009}
Gerhard Kowol.
\newblock {\em Projektive Geometrie und Cayley-Klein Geometrien der Ebene}.
\newblock Birkhauser, 2009.

\bibitem[Li08]{li08}
Hongbo Li.
\newblock {\em Invariant Algebras and Geometric Algebra}.
\newblock World Scientific, Singapore, 2008.

\bibitem[Lin73]{lindemann73}
F.~Lindemann.
\newblock {\"{U}}ber unendlich kleine {B}ewegungen und {\"{u}}ber
  {K}raftsysteme bei allgemeiner projektivischer {M}assbestimmung.
\newblock {\em Mathematische Annalen}, 7:56--144, 1873.

\bibitem[McC90]{mcc}
J.~Michael McCarthy.
\newblock {\em An Introduction to Theoretical Kinematics}.
\newblock MIT Press, Cambridge, MA, 1990.

\bibitem[Pap96]{pappas96}
Richard Pappas.
\newblock Oriented projective geometry with clifford algebra.
\newblock In {\em Clifford Algebras with Numerical and Symbolic Calculations},
  pages 233--251. Birkh{\"{a}}user, 1996.

\bibitem[Per]{clu}
Christian Perwass.
\newblock Clu project.
\newblock \url{http://www.perwass.de/clu/}.

\bibitem[Per09]{perwass09}
Christian Perwass.
\newblock {\em Geometric Algebra with Applications to Engineering}.
\newblock Springer, 2009.

\bibitem[PW01]{pottmann01}
Helmut Pottmann and Johannes Wallner.
\newblock {\em Computational Line Geometry}.
\newblock Springer, Berlin, 2001.

\bibitem[Sel00]{selig00}
Jon Selig.
\newblock Clifford algebra of points, lines, and planes.
\newblock {\em Robotica}, 18:545--556, 2000.

\bibitem[Sel05]{selig05}
Jon Selig.
\newblock {\em Geometric Fundamentals of Robotics}.
\newblock Springer, 2005.

\bibitem[Sto92]{stolfi91}
Jorge Stolfi.
\newblock {\em Oriented Projective Geometry}.
\newblock Academic, 1992.

\bibitem[Stu91]{study91}
Eduard Study.
\newblock Von den bewegungen und umlegungen.
\newblock {\em Mathematische Annalen}, 39:441--566, 1891.

\bibitem[Stu03]{study03}
Eduard Study.
\newblock {\em Geometrie der Dynamen}.
\newblock Tuebner, Leibzig, 1903.

\bibitem[TL97]{lt:gt3m}
William Thurston and Silvio Levy.
\newblock {\em Three-Dimensional Geometry and Topology: Volume 1}.
\newblock Princeton University Press, 1997.

\bibitem[vM24]{vm24}
Richard von Mises.
\newblock Die {M}otorrechnung: {E}ine {N}eue {H}ilfsmittel in der {M}echanik.
\newblock {\em Zeitschrift {f\"{u}}r Rein und Angewandte Mathematik und
  Mechanik}, pages 4(2):155--181, 1924.

\bibitem[Wei35]{weiss35}
Ernst~August Weiss.
\newblock {\em Einf{\"{u}}hrung in die Liniengeometrie und Kinematik}.
\newblock Teubner, Leibzig, 1935.

\bibitem[Whi98]{whitehead98}
A.~N. Whitehead.
\newblock {\em A Treatise on Universal Algebra}.
\newblock Cambridge University Press, 1898.

\bibitem[Wik]{wikiextalg}
Wikipedia.
\newblock \url{http://en.wikipedia.org/wiki/Exterior_algebra}.

\bibitem[Zie85]{ziegler}
Renatus Ziegler.
\newblock {\em Die Geschichte Der Geometrischen Mechanik im 19. Jahrhundert}.
\newblock Franz Steiner Verlag, Stuttgart, 1985.

\end{thebibliography}
\bibliographystyle{alpha}

\end{document}